\documentclass[10pt]{article}

\usepackage[latin1]{inputenc}
\usepackage[english,frenchb]{babel}
\usepackage{amssymb}
\usepackage{amsmath}
\usepackage{enumerate}
\usepackage{indentfirst}
\usepackage{color}
\usepackage{textcomp}
\usepackage{fancyhdr}
\usepackage{titlesec}

\newtheorem{theorem}{Theorem}[section]
\newtheorem{definition}{Definition}[section]
\newtheorem{corollary}{Corollary}[section]
\newtheorem{proposition}{Proposition}[section]
\newtheorem{lemma}{Lemma}[section]

\newtheorem{theoreme}{Théorème}[section]
\newtheorem{remark}{Remark}[section]

\newcommand{\eq}[1]{\begin{equation}#1\end{equation}}
\newcommand{\arr}[2]{\begin{array}{#1}#2\end{array}}
\newcommand{\theo}[1]{\begin{theorem}#1\end{theorem}}

\newcommand{\R}{\mathbb{R}}
\newcommand{\dv}{\hbox{div} \hspace{0.1cm}}
\newcommand{\curl}{\hbox{curl} \hspace{0.1cm}}

\newcommand{\lem}[1]{\begin{lemma}#1\end{lemma}}
\newcommand{\aligne}[1]{\begin{align*}#1\end{align*}}
\newcommand{\cent}[1]{\begin{center}#1\end{center}}

\makeatletter

\@addtoreset{equation}{section}
\makeatother

\title{Asymptotic profiles for the third grade fluids equations}
\author{Olivier Coulaud}
\date{}
\begin{document}
\selectlanguage{english}
\maketitle
\begin{abstract}
We study the long time behaviour of the solutions of the third grade fluids equations in dimension 2. Introducing scaled variables and performing several energy estimates in weighted Sobolev spaces, we describe the first order of an asymptotic expansion of these solutions. It shows in particular that, under smallness assumptions on the data, the solutions of the third grade fluids equations converge to self-similar solutions of the heat equations, which can be computed explicitly from the data.
\vspace{0.5cm}
\end{abstract}
\section{Introduction}
The study of the behaviour of the non-Newtonian fluids is a significant topic of research in mathematics, but also in physics or biology. Indeed, these fluids, the behaviour of which cannot be described with the classical Navier-Stokes equations, are found everywhere in the nature. For examples, blood, wet sand or certain kind of oils used in industry are non-Newtonian fluids. In this paper, we investigate the behaviour of a particular class of non-Newtonian fluids that is the third grade fluids, which are a particular case to the Rivlin-Ericksen fluids (see \cite{rivlinericksen55}, \cite{truesdellnoll65}). The constitutive law of such fluids is defined through the Rivlin-Ericksen tensors, given recursively by
\cent{
$
\arr{l}{A_1 = \nabla u + \left(\nabla u\right)^t,\\
A_k = \partial_t A_{k-1} + u.\nabla A_{k-1}+\left(\nabla u\right)^t A_{k-1}+A_{k-1} \nabla u,
}
$
}
where $u$ is a divergence free vector field of $\R^2$ or $\R^3$ which represents the velocity of the fluid. The most famous example of a Rivlin-Ericksen fluid is the class of the Newtonian fluids, which are modelized through the stress tensor
\cent{$\sigma = - p I + \nu A_1$,}
where $\nu> 0$ is the kinematic viscosity and $p$ is the pressure of the fluid. Introduced into the equations of conservation of momentum, this stress tensor leads to the well known Navier-Stokes equations.
\vspace{0.5cm}\\
In this article, we consider a larger class of fluids, for which the stress tensor is not linear in the Rivlin-Ericksen tensors, but a polynomial function of degree 3. As introduced by Fosdick and Rajagopal in \cite{fosdickrajagopal80}, the stress tensor that we consider is defined by
\cent{
$
\sigma = -p I + \nu A_1 + \alpha_1 A_2 + \alpha_2 A^2_1 + \beta \left|A_1\right|^2 A_1,
$
}
where $\nu> 0$ is the kinematic viscosity, $p$ is the pressure, $\alpha_1> 0$, $\alpha_2 \in \R$ and $\beta\geq 0$.
\vspace{0.5cm}\\
We assume in this article that the density of the fluid is constant in space and time and equals $1$. Actually, the value of the density is not significant, since we can replace the parameters $\nu$, $\alpha_1$, $\alpha_2$ and $\beta$ by dividing them by the density. Introduced into the equations of conservation of momentum, the tensor $\sigma$ leads to the system
\eq{\label{g3u}
\arr{l}{
\partial_t\left(u-\alpha_1\Delta u\right)-\nu \Delta u + \curl\left(u-\alpha_1 \Delta u\right)\wedge u\\
\hspace{3cm} -\left(\alpha_1+\alpha_2\right)\left(A.\Delta u+2\dv\left(L L^t\right)\right)-\beta \dv \left(\left|A	\right|^2 A\right)+\nabla p=0,\\
\dv u = 0,\\
u_{\left|t=0\right.}=u_0,
}
}
where $L = \nabla u$,  $A(u)=\nabla u +\left(\nabla u\right)^t$ and $\wedge$ denotes the classical vectorial product of $\R^3$. For matrices $A,B \in \mathcal M_d(\R)$, we define $\displaystyle A : B = \sum^d_{i,j=1} A_{i,j} B_{i,j}$ and $\left|A\right|^2 = A:A$. If the space dimension is 2, we use the convention $u=(u_1,u_2,0)$ and $\curl u = (0,0,\partial_1 u_2-\partial_2 u_1)$. Notice also that if $\alpha_1 + \alpha_2=0$ and $\beta=0$, we recover the equations of motion of second grade fluids, which are another class of non-Newtonian fluids, introduced earlier by Dunn and Fosdick in 1974 (see \cite{dunnfosdick74}, \cite{galdivandalsensauer93} or \cite{cioranescuelhacene84}). If in addition $\alpha_1=0$, then one recovers the classical Navier-Stokes equations.
\vspace{0.5cm}\\
\indent The system of equations (\ref{g3u}) has been studied in various cases, on bounded domains of $\R^d$, $d=2,3$ or in the whole space $\R^d$ (see \cite{amrouchecioranescu97}, \cite{bernard99}, \cite{breschlemoine99}, \cite{busuiociftimie04}, \cite{busuiociftimie06} and \cite{paicu08}). On a bounded domain $\Omega$ of $\R^d$ with Dirichlet boundary conditions, Amrouche and Cioranescu have shown the existence of local solutions to (\ref{g3u}) when the initial data belong to the Sobolev space $H^3(\Omega)^d$ (see \cite{amrouchecioranescu97}). In addition, these solutions are unique. For this study, the authors have assumed the restriction 
\cent{$\left|\alpha_1+\alpha_2\right|\leq \left(24 \nu \beta\right)^{1/2}$,}
which is justified by thermodynamics considerations. The proof of their result is obtained via a Galerkin method with functions belonging to the eigenspaces of the operator $\quad \curl \left(I-\alpha_1 \Delta \right)$. In dimension 3, a slightly different method has been applied by D. Bresch and J. Lemoine, who used Schauder's fixed point Theorem to extend the result of \cite{amrouchecioranescu97} to the case of initial data belonging to the Sobolev spaces $W^{2,r}(\Omega)^3$, with $r>3$.  They have shown in \cite{breschlemoine99} the local existence of unique solutions of (\ref{g3u}) in the space $C^0\left(\left[0,T\right],W^{2,r}(\Omega)^3\right)$, where $T>0$. In addition, if the data are small enough in the space $W^{2,r}(\Omega)^3$, the solutions are global in time. Notice also that the existence of such solutions holds without restrictions on the parameters of the system (\ref{g3u}).
\vspace{0.5cm}\\
\indent In the case of third grade fluids filling the whole space $\R^d$, $d=2,3$ , V. Busuioc and D. Iftimie have established the existence of global solutions with initial data belonging to $H^2(\R^d)^d$, without restrictions on the parameters or on the size of the data (see \cite{busuiociftimie04}). In this study, the authors used a Friedrichs scheme and performed a priori estimates in $H^2$ which allow to show the existence of solutions of (\ref{g3u}) in the space $L^\infty_{loc}\left(\R^+, H^2(\R^d)^d\right)$. Besides, these solutions are unique if $d=2$. Later, M. Paicu has extended the results of \cite{busuiociftimie04} to the case of initial data belonging to $H^1(\R^d)^d$, assuming additional restrictions on the parameters of the equation ; the uniqueness is not known in this space (see \cite{paicu08}). The method that he used is slightly different from the one used in \cite{busuiociftimie04}. Indeed, although M. Paicu also considered a Friedrichs scheme, the convergence of the approximate solutions to a solution of (\ref{g3u}) is done via a monotonicity method. Notice that Theorem \ref{theo1} of this article shows the existence of solutions of the equations of third grade fluids on $\R^2$ for initial data in weighted Sobolev spaces (see Section 3). 
\vspace{0.5cm}\\
\indent In what follows, we consider a third grade fluid filling the whole space $\R^2$. Actually, the equations that we consider are not exactly the system (\ref{g3u}) but the one satisfied by $w=\curl u = \partial_1 u_2-\partial_2 u_1$. In dimension 2, the vorticity equations of the third garde fluids are given by
\eq{
\label{g3w}
\arr{l}{
\partial_t\left(w-\alpha_1 \Delta w\right) - \nu \Delta w + u .\nabla \left(w-\alpha_1 \Delta w\right)\\
\hspace{3cm}-\beta \dv\left(\left|A\right|^2\nabla w\right)-\beta \dv\left(\nabla\left(\left|A\right|^2\right)\wedge A\right) =0,\\
\dv u=0,\\
w_{\left|t=0\right.}=w_0=\curl u_0.
}
}
Notice that the parameter $\alpha_2$ does no longer appear in (\ref{g3w}) and thus does not play any role in the study of these equations. Indeed, due to the divergence free property of $u$, a short computation shows that $\quad \curl \left(A.\Delta u+2\dv\left(L L^t\right)\right)=0$, or equivalently there exists $q$ such that $A.\Delta u+2\dv\left(L L^t\right) = \nabla q$. This phenomenon is very particular to the dimension 2 and does not occur in dimension 3. Notice also that the previous system is autonomous in $w$. Indeed, the vector field $u$ depends on $w$ and can be recovered from $w$ via the Biot-Savart law, which is a way to get a divergence free vector field such that $\curl u=w$. The motivation for considering the vorticity equations instead of the equations of motion comes from the fact that, due to spectral reasons, we have to study the behaviour of the solutions of (\ref{g3w}) in weighted Lebesgue spaces. Indeed, in what follows, we will consider scaled variables, which make appear a differential operator whose essential spectrum can be "pushed to the left" by taking a convenient weighted Lebesgue space. We will see that the rate of convergence of the solutions of (\ref{g3w}) is linked to the spectrum of this operator. Unfortunately, the weighted Lebesgue spaces are not suitable for the equations of motions and are not preserved by the system (\ref{g3u}). Anyway, one can obtain the asymptotic profiles of the solutions of the equations of motion (\ref{g3u}) from the study of the asymptotic behaviour of the solutions of the vorticity equations (see Corollary \ref{cor1} below). We also emphasize that the system (\ref{g3w}) allows to consider solutions whose velocity fields are not bounded in $L^2$.
\vspace{0.5cm}\\
\indent In this article, we establish the existence and uniqueness of solutions of (\ref{g3w}) in weighted Sobolev spaces, but the main aim is the study of the asymptotic behaviour of these solutions when $t$ goes to infinity. More precisely, we want to describe the first order asymptotic profiles of the solutions of (\ref{g3w}). We consider a fluid of third grade which fills $\R^2$ without forcing term applied to it. In this case, as it is expected, the solutions of (\ref{g3w}) tend to $0$ as $t$ goes to infinity. Our motivation is to show that these solutions behave like those of the Navier-Stokes equations. In our case, we will show that the solutions of (\ref{g3w}) behave asymptotically like solutions of the heat equations, up to a constant that we can compute from the initial data. The methods that we use in the present paper are based on scaled variables and energy estimates in several functions spaces. This work is inspired by several older results obtained for other fluid mechanics equations. The first and second order asymptotic profiles have been described for the Navier-Stokes equations in dimensions 2 and 3 by T. Gallay and E. Wayne (see\cite{gallaywayne02}, \cite{gallaywayne02bis}, \cite{gallaywayne05} and \cite{gallaywayne06}).  In dimension 2, they have shown in \cite{gallaywayne02} and \cite{gallaywayne05} that the first order asymptotic profiles of the Navier-Stokes equations are given up to a constant by a smooth Gaussian function called the Oseen vortex sheet. More precisely, for a solution $w$ of the vorticity Navier-Stokes equations (that is the system \ref{g3w} with $\alpha_1 = \beta = 0$), for every $2\leq p \leq + \infty$, the following property holds:
\cent{
$
\displaystyle \left\|w(t) - \frac{\int_{\R^2} w_0(x) dx}{t}G\left(\frac{.}{\sqrt{t}}\right)\right\|_{L^p} = \mathcal O (t^{-\frac{3}{2}+\frac{1}{p}})
$, when $t\rightarrow +\infty$,
}
where $G$ is the Oseen vortex sheet
\eq{
\label{oseen}
G(x)= \frac{1}{4\pi}e^{- \frac{\left|x\right|^2}{4}}.
}
The methods that they used in \cite{gallaywayne02} are very different from the ones that we develop in this article. Although they also considered scaled variables, the convergence to the asymptotic profiles is not obtained through energy estimates. Indeed, using dynamical systems arguments, they established the existence of a finite-dimensional manifold which is locally invariant by the semiflow associated to the Navier-Stokes equations. Then, they showed that, under restrictions on the size of the data, the solutions of the Navier-Stokes equations behave asymptotically like solutions on this invariant manifold. The description of the asymptotic profiles is thus obtained by the description of the dynamics of the Navier-Stokes equations on the invariant manifold. Later, the smallness assumption on the data has been removed (see \cite{gallaywayne05}). In \cite{jaffal11}, Jaffal-Mourtada describes the first order asymptotics of second grade fluids, under smallness assumptions on the initial data in weighted Sobolev spaces. She has shown that the solutions of the second grade fluids equations converge also to the Oseen vortex sheet. In this paper, we apply the methods used by Jaffal-Mourtada, namely scaled variables and energy estimates. According to these results, one can say that the fluids of second grade behave asymptotically like Newtonian fluids. In this paper, we show that, under the same smallness assumptions on the initial data, the same behaviour occurs for the third grade fluids equations. We emphasize that the rate of convergence that we obtain is better than the one obtained in \cite{jaffal11}. Actually, we show that we can choose the rate of convergence as close as wanted to the optimal one, assuming that the initial data are small enough. Since second grade fluids are a particular case of third grade fluids, we establish an improvement of the rate obtained in \cite{jaffal11}. Actually, the main difference between third and second grade fluids equations in dimension $2$ is the presence of the additional term $\quad\beta \dv\left(\left|A\right|^2 A\right)\quad$ in the third grade fluids equations. Sometimes, this term helps to obtain global estimates, like in \cite{busuiociftimie04} or \cite{paicu08}, but introduces additional difficulties when one looks for estimates in $H^3$ or in more regular Sobolev spaces (see \cite{amrouchecioranescu97}, \cite{bernard99} or \cite{busuiociftimie06}). Here, we have to establish estimates in weighted Sobolev spaces with $H^2$ regularity for the vorticity $w$, which is harder than doing estimates in $H^3$ for $u$.
\vspace{0.5cm}\\
\indent We next introduce scaled variables. In order to simplify the notations, we assume that $\nu=1$. Let $T>1$ be a positive constant which is introduced in order to avoid restrictions on the size of the parameter $\alpha_1$ and which will be made more precise later. We consider the solution $w$ of (\ref{g3w}) and define $W$ and $U$ such that $\curl U=W$ through the change of variables $\displaystyle X= \frac{x}{\sqrt{t+T}}$ and $\tau = \log(t+T)$. We set
\eq{\label{scaledv1}
\left\{\arr{l}{
\displaystyle u(t,x)=\frac{1}{\sqrt{t+T}} U\left( \log(t+T),\frac{x}{\sqrt{t+T}}\right),\\
\displaystyle w(t,x)=\frac{1}{t+T}W\left(\log (t+T),\frac{x}{\sqrt{t+T}}\right).
}\right.
}
For $\tau \geq \log(T)$, we have
\eq{\label{scaledv2}
\left\{\arr{l}{
U(\tau,X)=e^{\tau/2}u\left(e^\tau-T,e^{\tau/2}X\right),\\
W(\tau,X)=e^\tau w\left(e^\tau-T,e^{\tau/2}X\right).
}\right.
}
These variables, called scaled or self-similar variables have been introduced in order to study the long time asymptotic of solutions of parabolic equations and particularly to show the convergence to self-similar solutions (see \cite{escobedokavianmatano95}, \cite{escobedozuazua91}, \cite{galaktionov91} or \cite{kavian87}), that is to say under the form $\displaystyle \frac{1}{t+T}F\left(\frac{x}{\sqrt{t+T}}\right)$.
\vspace{0.5cm}\\
Scaled variables have been used to deal with the asymptotic behaviour of many equations, not necessarily parabolic ones (see  \cite{carpio94}, \cite{carpio96} \cite{jaffal11}, \cite{gallayraugel98} or \cite{gallayraugel00}). For instance, in \cite{gallayraugel98}, T. Gallay and G. Raugel have described the first and second order asymptotic profiles in weighted Sobolev spaces for damped wave equations, using scaled variables. In \cite{gallayraugel00}, they use scaled variables to show a stability result of hyperbolic fronts for the same equations.
\vspace{0.5cm}\\
For sake of simplicity, we set $A^{i,j}=\partial_j U_i + \partial_i U_j$. Considering self-similar variables, one can see that $W$ and its corresponding divergence free vector field $U$ satisfy the system
\eq{
\label{g3W}
\arr{l}{
\partial_\tau\left(W-\alpha_1 e^{-\tau}\Delta W\right)- \mathcal L (W)+ U.\nabla \left(W-\alpha_1 e^{-\tau}\Delta W\right)+\alpha_1 e^{-\tau}\Delta W\\
\hspace{2cm} +\alpha_1 e^{-\tau}\frac{X}{2}.\nabla \Delta W-\beta e^{-2\tau}\dv\left(\left|A\right|^2\nabla W\right)-\beta e^{-2\tau}\dv\left(\nabla\left(\left|A\right|^2\right)\wedge A\right)=0,\\
\dv U=0,\\
W_{\left|\tau=\tau_0\right.}=W_0,
}
}
where $\tau_0 = \log(T)$, $W_0(X) = e^{\tau_0} w_0 \left(e^{\tau_0/2}X\right)$ and $\mathcal L$ is the linear differential operator defined by
\cent{
$\mathcal L(W)=\Delta W + W + \frac{X}{2}.\nabla W$.
}
Notice that the initial time of the system (\ref{g3W}) is $\log (T)$. By choosing $T$ sufficiently large, one can consider $\alpha_1 e^{-\tau}$ as small as wanted. This fact allows to study the behaviour of the solutions of (\ref{g3W}) without restrictions on the size of $\alpha_1$. Formally, we see that most of the terms of the system (\ref{g3W}) tend to $0$ as time goes to infinity. The purpose of the present paper is to show that the solutions of (\ref{g3W}) asymptotically behave like solutions of
\eq{\label{g3lim}
\partial_\tau W_\infty = \mathcal{L}(W_\infty).
}
In order to describe the solutions of the system (\ref{g3lim}), we have to study the spectrum of the linear differential operator $\mathcal L$ in appropriate functions spaces. The form of the previous system and the definition of $\mathcal L$ lead to consider weighted Lebesgue spaces. For $m\in \mathbb N$, we define
\cent{
$L^2(m)= \left\{u\in L^2(\R^2) : \left(1+\left|x\right|^2\right)^{m/2} u \in L^2(\R^2)\right\}$,
}
equipped with the norm 
\cent{$\left\|u\right\|_{L^2(m)} = \displaystyle \left(\int_{\R^2} \left(1+\left|x\right|^2\right)^m \left|u(x)\right|^2 dx\right)^{1/2}$.}
The spectrum of $\mathcal L$ in $L^2(m)$ is given in \cite[Appendix A]{gallaywayne02}. It is composed of the discrete spectrum
\cent{$
\sigma_d\left(\mathcal L\right)=\left\{-\frac{k}{2} : k\in \left\{0,1,...,m-2\right\}\right\},
$}
and the continuous spectrum
\cent{$
\sigma_c (\mathcal L)=\left\{\lambda \in \mathbb C : \hbox{Re}(\lambda) \leq -\frac{m-1}{2}\right\}.
$}
In particular, the eigenvalue $0$ is simple and the Oseen vortex $G$ given by (\ref{oseen}) is an eigenfunction of $\mathcal L$ associated to $0$. Of course, $G$ is a solution of (\ref{g3lim}) and we will show that the solutions of (\ref{g3W}) behave like $G$ when the time goes to infinity. To this end, we decompose the solutions $W$ of (\ref{g3W}) as follows
\cent{$
W(\tau)=\eta G + f(\tau),
$}
where $\eta \in \R$ will be made more precise later and $f(\tau)$ is a rest which will tend to $0$ as $\tau$ goes to infinity.\\
\vspace{0.5cm}\\
In order to get a good rate of convergence for $f$, we shall "push" the continuous spectrum of $\mathcal L$ to the left by choosing an appropriate weighted Lebesgue space. For this reason, we work in $L^2(2)$, so that $
\sigma_c (\mathcal L)=\left\{\lambda \in \mathbb C : \hbox{Re}(\lambda) \leq -\frac{1}{2}\right\}
$. Since the second eigenvalue of $\mathcal L$ in $L^2(2)$ is $-\frac{1}{2}$, the best result that we expect is
\cent{
$
\displaystyle f(\tau) = \mathcal{O}(e^{- \tau/2}) \quad$ in $L^2(2), \quad$ when $\quad \tau \rightarrow+\infty.
$
}
Notice that choosing a weighted space $L^2(m)$ with $m>2$ would be useless for describing the first order asymptotics only. Indeed, if we take $m>2$, the second eigenvalue would still be $-\frac{1}{2}$ and the rate of convergence could not be better than $e^{-\tau/2}$.
\vspace{0.5cm}\\
For later use, we define the divergence free vector field $V$ such that $\curl V=G$. It is obtained by the Biot-Savart law and given by
\eq{\label{oseenV}
\displaystyle V(X)= \frac{1-e^{-\frac{\left|X\right|^2}{4}}}{2\pi\left|X\right|^2}\left(\arr{c}{-X_2\\X_1}\right).
}
In particular, for every $X\in \R^2$, one has
\cent{$
V(X).X=0, \hspace{0.2cm} V(X).\nabla G(X) = 0 \hspace{0.2cm} \hbox{and} \hspace{0.2cm} V(X).\nabla \Delta G(X)=0.
$}
Before stating the main theorem of this paper, we have to define some additional functions spaces. For $m\in \mathbb N$, we set
\aligne{
&H^1(m)=\left\{u\in L^2(m) : \partial_j u \in L^2(m); j\in\{1,2\}\right\},\\
&H^2(m)=\left\{u\in H^1(m) : \partial_j u \in H^1(m); j\in\{1,2\}\right\},
}
equipped with the norms 
\cent{$\left\|u\right\|_{H^1(m)} = \left(\left\|u\right\|^2_{L^2(m)}+\left\|\nabla u\right\|^2_{L^2(m)}\right)^{1/2}$\hspace{0.5cm} and \hspace{0.5cm}$\left\|u\right\|_{H^2(m)} = \left(\left\|u\right\|^2_{H^1(m)}+\left\|\nabla^2 u\right\|^2_{L^2(m)}\right)^{1/2}$,}
where $\displaystyle\left|\nabla u\right|^2 = \sum^2_{i=1} \left(\partial_i u\right)^2$ and $\displaystyle\left|\nabla^2 u\right|^2 = \sum^2_{i,j=1} \left(\partial_i\partial_j u\right)^2$.
\vspace{0.5cm}\\
The following theorem describes the first order asymptotic profile of $W$ in $H^2(2)$, if one assumes that the initial data $W_0$ are small enough in the weighted Sobolev space $H^2(2)$.
\begin{theorem}\label{theo1}
Let $\theta$ be a constant such that $\displaystyle 0<\theta < 1$. There exist two positive constants $\gamma_0=\gamma_0(\alpha_1,\beta)$ and $T_0=T_0(\alpha_1)\geq 1$ such that, for all $W_0\in H^2(2)$ satisfying the condition
\eq{\label{cond1}
\left\|W_0\right\|^2_{H^1}+\frac{\alpha_1}{T}\left\|\Delta W_0\right\|^2_{L^2}+\left\|\left|X\right|^2 W_0\right\|^2_{L^2}+\frac{\alpha^2_1}{T^2}\left\|\left|X\right|^2\Delta W_0\right\|^2_{L^2} \leq \gamma \left(1-\theta\right)^6,
}
for some $T \geq T_0$ and $0< \gamma \leq \gamma_0$,
\vspace{0.5cm}\\
there exist a unique global solution $W\in C^0\left(\left[\tau_0,+\infty\right),H^2(2)\right)$ of $(\ref{g3W})$ and a positive constant $C=C(\alpha_1,\beta,\theta)$ such that, for all $\tau\geq \tau_0$,
\eq{\label{inetheo1}
\left\|\left(1-\alpha_1e^{-\tau}\Delta\right)\left(W(\tau)-\eta G\right)\right\|^2_{L^2(2)}\leq C \gamma e^{-\theta \tau},}
where $\eta=\displaystyle \int_{\R^2}W_0(X) dX$, $\tau_0 =\log (T)$ and the parameters $\alpha_1$ and $\beta$ are fixed and given in (\ref{g3u}).
\end{theorem}
\begin{remark}The smallness assumption (\ref{cond1}) is not optimal. By working harder, it is possible to get $\gamma\left(1-\theta\right)^p$ with $p<6$ in the right hand side of the inequality.
\end{remark}
\begin{remark}
Notice that Theorem \ref{theo1} establishes an improvement of \cite[Theorem 1.1]{jaffal11} concerning the first order asymptotics of the second grade fluids equations. Indeed, the above theorem holds also with $\beta = 0$ and consequently describes the first order asymptotic profiles of the solutions of the second grade fluids equation. The improvement comes from the fact that one can choose $\theta$ as close as wanted to $1$, which is the optimal rate. In \cite{jaffal11}, the constant $\theta$ can not be bigger than $\frac{1}{2}$.
\end{remark}
Theorem \ref{theo1} implies the following result in the unscaled variables. In particular, it gives a description of the asymptotic profiles of the solutions of the equations of motion (\ref{g3u}).
\begin{corollary}\label{cor1}
Let $\theta$ be a constant such that $0<\theta<1$. There exist two positive constants $\gamma_0=\gamma_0(\alpha_1,\beta)$ and $T_0=T_0(\alpha_1,\beta)\geq 1$ such that, for all $w_0\in H^2(2)$ satisfying the condition
\eq{
T\left\|w_0\right\|^2_{L^2} + T^2\left\|\nabla w_0\right\|^2_{L^2}+ \frac{1}{T}\left\|\left|x\right|^2 w_0 \right\|^2_{L^2}+\alpha_1 T^3 \left\|\Delta w_0\right\|^2_{L^2}+\frac{\alpha^2_1}{T} \left\|\left|x\right|^2 \Delta w_0 \right\|^2_{L^2} \leq  \gamma \left(1-\theta\right)^6,
}
for some $T \geq T_0$ and $0< \gamma \leq \gamma_0$,
\vspace{0.5cm}\\
there exists a unique global solution $w \in C^0\left(\left[0,+\infty\right),H^2(2)\right)$ of (\ref{g3w}) such that, for all $1\leq p\leq 2$, there exists a positive constant $C=C(\alpha_1,\beta,\theta)$ such that, for all $t\geq 0$,
\cent{
$
\displaystyle \left\|\left(1-\alpha_1 \Delta\right)\left(w(t)-\frac{\eta}{{t+T}}G\left(\frac{x}{\sqrt{t+T}}\right)\right)\right\|_{L^p} \leq C \gamma\left(t+T\right)^{-1-\frac{\theta}{2}+\frac{1}{p}},
$
}
where $\displaystyle\eta=\int_{\R^2} w_0(x)dx$.
\vspace{0.5cm}\\
Moreover, for all $2< q<+\infty$, there exists a positive constant $C=C(\alpha_1,\beta,\theta,q)$ such that, for all $t\geq 0$,
\cent{
$
\displaystyle \left\|\left(1-\alpha_1 \Delta\right)\left(u(t)-\frac{\eta}{{\sqrt{t+T}}}V\left(\frac{x}{\sqrt{t+T}}\right)\right)\right\|_{L^q}\leq C\gamma \left(t+T\right)^{-\frac{1}{2}-\frac{\theta}{2}+\frac{1}{q}},
$
}
where $V$ is obtained from $G$ via the Biot-Savart law and defined by (\ref{oseenV}).
\end{corollary}
~\\
\indent Theorem \ref{theo1} describes the asymptotic behaviour of the solutions of $(\ref{g3W})$ in $H^2(2)$ at the first order. Since the solutions of the Navier-Stokes equations converge also to the Oseen vortex sheet, we can say that the fluids of third grade behave asymptotically like Newtonian fluids. Notice that the functions space $H^2(2)$ is suitable for the first order asymptotics because it "pushes" the continuous spectrum of $\mathcal L$ far enough to get $0$ as an isolated eigenvalue. If we had to describe the asymptotics of $(\ref{g3W})$ at the second order, we should work in a space where $\mathcal L$ has at least two isolated eigenvalues. Due to the forms of $\sigma_c$ and $\sigma_d$, the second order asymptotics must be studied in functions space with polynomial weight of degree at least 3, in order to get the two isolated eigenvalues $0$ and $-\frac{1}{2}$.
\vspace{0.5cm}\\
Notice also that as the system (\ref{g3w}) and our change of variables preserve the total mass. We have, for all $\tau\geq \tau_0$ and $t\geq 0$,
\cent{
$\displaystyle \eta = \int_{\R^2} w_0(x)dx=\int_{\R^2} w(t,x)dx=\int_{\R^2} W_0(X)dX=\int_{\R^2} W(\tau,X)dX$.
}
~\\
\indent The plan of this article is as follows. In Section \ref{secbiot}, we recall classical results concerning the Biot-Savart law and give several technical lemmas. In Section \ref{secapprox}, we introduce a regularized system, which is close to (\ref{g3W}) and depends on a small parameter $\varepsilon>0$. Actually, we add the regularizing term $\varepsilon \Delta^2 W$ to the system (\ref{g3W}) and show the existence of unique regular solutions $W_\varepsilon$ to this new system. In Section \ref{secenergy}, using energy estimates in various functions spaces, we show that $W_\varepsilon$ satisfies the inequality (\ref{inetheo1}) of Theorem \ref{theo1}, and thus tends to the Oseen vortex sheet $G$ when $\tau$ goes to infinity. In Section \ref{secdem}, we let $\varepsilon$ go to $0$ and show that $W_\varepsilon$ tends in a sense to a solution $W$ of (\ref{g3W}). Additionally, this solution satisfies the inequality (\ref{inetheo1}) of Theorem \ref{theo1} and consequently tends also to the Oseen vortex sheet. Finally, we establish the uniqueness of $W$, which enables us to say that every solution of (\ref{g3W}) satisfying the assumption (\ref{cond1}) converges to the Oseen Vortex sheet when $\tau$ goes to infinity.
\section{\label{secbiot}Biot-Savart law and auxiliary lemmas}
In this section, we state several technical lemmas which are useful to prove Theorem \ref{theo1}. These lemmas concern the Biot-Savart law and state several inequalities involving weighted Lebesgue norms. In what follows, we use the notation
\cent{$\left\|u\right\|=\left\|u\right\|_{L^2}$,}
and $C$ denotes a positive constant which can depend on the fixed constants $\alpha_1$ and $\beta$.
\vspace{0.5cm}\\
The first lemma will be useful in Section \ref{secenergy} to obtain estimates in Sobolev spaces of negative order. We define, for $s\in \R$, the operator $\left(-\Delta \right)^s$, given by
\cent{
$\left(-\Delta\right)^{s} u =\bar{\mathcal F}\left(\left|\xi\right|^{2s} \widehat{u}\right)$,
}
where $\widehat{u}$ (also denoted $\mathcal F (u)$) is the Fourier transform of $u$, given by
\cent{
$
\displaystyle \widehat{u}(\xi)= \int_{\R^2} u(x) e^{-i x.\xi} dx,
$
} 
and $\bar{\mathcal F}$ denotes the inverse Fourier transform
\cent{
$
\displaystyle \bar{\mathcal F}(v)(x)= \frac{1}{\left(2\pi\right)^2}\int_{\R^2} v(\xi) e^{i x.\xi} d\xi.
$
}
\lem{\label{weight-}
Let $s$ be a positive real number such that $\frac{3}{4}<s<1$, then we have the following two inequalities.
\begin{enumerate}[1.]
\item Let $g \in L^2(1)$. Then $\left(-\Delta\right)^{-s}\nabla g \in L^2(\R^2)$ and there exists $C>0$ independent of $g$ and $s$ such that
\eq{\label{weight-nabla}
\left\|\left(-\Delta\right)^{-s}\nabla g\right\| \leq \frac{C}{\left(1-s\right)^{3/2}}\left\|g\right\|_{L^2(1)}.
}
\item  Let $g \in L^2(2)$ such that $\displaystyle \int_{\R^2} g(x) dx=0$. Then $\left(-\Delta\right)^{-s} g \in L^2(\R^2)$ and there exists $C>0$ independent of $g$ and $s$ such that
\eq{\label{weight-f}
\left\|\left(-\Delta\right)^{-s} g\right\|_{L^2} \leq \frac{C}{\left(1-s\right)^{3/2}}\left\|g\right\|_{L^2(2)}.
}
\end{enumerate}
}
\textbf{Proof : }We start by proving the inequality (\ref{weight-nabla}). For $j\in \left\{1,2\right\}$, using Fourier variables, one has
\aligne{
\left\|\left(-\Delta\right)^{-s}\partial_j g\right\|^2_{L^2}&\leq C \int_{\left|\xi\right|\leq 1}\frac{1}{\left|\xi\right|^{4s-2}}\left| \widehat{g}\right|^2 d\xi+\left\|g\right\|^2_{L^2}\\
&\leq C \left(\int_{\left|\xi\right|\leq 1}\frac{1}{\left|\xi\right|^{2s}}d\xi \right)^{\frac{2s-1}{s}}\left(\int_{\left|\xi\right|\leq 1}\left| \widehat{g}\right|^{\frac{2s}{1-s}} d\xi\right)^{\frac{1-s}{s}}+\left\|g\right\|^2_{L^2}\\
&\leq \frac{C}{\left(1-s\right)}\left\| \widehat{g}\right\|^2_{L^{\frac{2s}{1-s}}}+\left\|g\right\|^2_{L^2}.
}
We now use the continuous injection of $H^1(\R^2)$ into $\displaystyle L^{\frac{2s}{1-s}}(\R^2)$. Looking at the computations of \cite[p. 723-724]{cheminxu97}, one can see that there exists a constant $C>0$ such that
\eq{\label{ineweight-}
\displaystyle\left\|u\right\|_{L^p} \leq C p\left\|u\right\|_{H^1}
, \hspace{0.2cm} \hbox{for all} \hspace{0.2cm} u\in H^1(\R^2) \hspace{0.2cm} \hbox{and} \hspace{0.2cm} 2\leq p < +\infty.
}
Notice that $Cp$ is not the optimal constant in the previous inequality. Using the inequality (\ref{ineweight-}), one has
\aligne{
\left\|\left(-\Delta\right)^{-s} \partial_j g\right\|^2_{L^2} &\leq \frac{C}{\left(1-s\right)^3}\left\| \widehat{g}\right\|^2_{H^1}+\left\|g\right\|^2_{L^2}\\
&\leq \frac{C}{\left(1-s\right)^3}\left\|g\right\|^2_{L^2(1)}.
}
We now prove the inequality (\ref{weight-f}). Since $\displaystyle \int_{\R^2} f(x)dx =0$, using Fourier variables, we get
\aligne{
\left\|\left(-\Delta\right)^{-s} g\right\|^2_{L^2}&= \left(2\pi\right)^2\int_{\R^2} \frac{1}{\left|\xi\right|^{4s}}\left| \widehat{g}(\xi) \right|^2d\xi\\
&\leq \left(2\pi\right)^2\int_{\left|\xi\right|\leq 1} \frac{1}{\left|\xi\right|^{4s}}\left| \widehat{g}(\xi) \right|^2d\xi+ \left\|g\right\|^2_{L^2}\\
&\leq \left(2\pi\right)^2\int_{\left|\xi\right|\leq 1} \frac{1}{\left|\xi\right|^{4s}}\left|\int^1_0\xi.\nabla  \widehat{g}(\sigma\xi)d\sigma\right|^2 d\xi+ \left\|g\right\|^2_{L^2}\\
&\leq C\int_{\left|\xi\right|\leq 1} \frac{1}{\left|\xi\right|^{4s-2}}\left|\int^1_0\left|\nabla  \widehat{g}(\sigma\xi)\right|d\sigma\right|^2 d\xi+ \left\|g\right\|^2_{L^2}.
}
Cauchy-Schwarz inequality and Fubini's theorem give
\aligne{
\left\|\left(-\Delta\right)^{-s} g\right\|^2_{L^2}&\leq C\int^1_0 \int_{\left|\xi\right|\leq 1} \frac{1}{\left|\xi\right|^{4s-2}}\left|\nabla  \widehat{g}(\sigma\xi)\right|^2 d\xi d\sigma + \left\|g\right\|^2_{L^2}.
}
Using Hölder inequality, we get
\aligne{
\left\|\left(-\Delta\right)^{-s} g\right\|^2_{L^2}&\leq C\int^1_0 \left(\int_{\left|\xi\right|\leq 1} \frac{1}{\left|\xi\right|^{4-\frac{2}{s}}}d\xi\right)^s\left(\int_{\left|\xi\right|\leq 1}\left|\nabla  \widehat{g}(\sigma\xi)\right|^{\frac{2}{1-s}} d\xi\right)^{1-s} d\sigma + \left\|g\right\|^2_{L^2}\\
&\leq C\left(\frac{s}{1-s}\right)^s\int^1_0 \left(\int_{\left|\xi\right|\leq 1}\left|\nabla  \widehat{g}(\sigma\xi)\right|^{\frac{2}{1-s}} d\xi\right)^{1-s} d\sigma + \left\|g\right\|^2_{L^2}.
}
The change of variables $\zeta = \sigma \xi$ yield
\aligne{
\left\|\left(-\Delta\right)^{-s} g\right\|^2_{L^2}&\leq C\left(\frac{s}{1-s}\right)^s\int^1_0 \left(\int_{\left|\zeta\right|\leq \sigma}\frac{1}{\sigma^2}\left|\nabla  \widehat{g}(\zeta)\right|^{\frac{2}{1-s}} d\zeta\right)^{1-s} d\sigma + \left\|g\right\|^2_{L^2}\\
&\leq C\left(\frac{s}{1-s}\right)^s\left(\frac{1}{2s-1}\right)\left\|\nabla  \widehat{g}\right\|^2_{L^{\frac{2}{1-s}}} + \left\|g\right\|^2_{L^2}.
}
Finally, we use again the inequality (\ref{ineweight-}) and obtain
\aligne{
\left\|\left(-\Delta\right)^{-s} g\right\|^2_{L^2}&\leq C\left(\frac{s}{1-s}\right)^s\left(\frac{1}{2s-1}\right)\left(\frac{2}{1-s}\right)^2\left\| \widehat{g}\right\|^2_{H^2} + \left\|g\right\|^2_{L^2}\\
&\leq \frac{C}{\left(1-s\right)^3}\left\|g\right\|^2_{L^2(2)},
}
which concludes the proof of this lemma.
\begin{flushright}
$\square$
\end{flushright}
\lem{\label{weight2}
\begin{enumerate}[1.]
\item Let $1\leq p <+\infty$ and $f\in L^p(\R^2)$ such that $\left|x\right|^2 f \in L^p(\R^2)$, then $\left|x\right| f \in L^p(\R^2)$ and the following inequality holds :
\eq{\label{weightLp}
\left\|\left|x\right| f\right\|_{L^p} \leq \left\| f\right\|^{1/2}_{L^p}\left\|\left|x\right|^2 f\right\|^{1/2}_{L^p}.
}
\item Let $f \in H^2(2)$, there exists $C>0$ such that
\eq{\label{weightnabla2}
\left\|\left|x\right|^2 \nabla^2 f\right\| \leq C\left(\left\|f\right\|+\left\|\left|x\right| \nabla f\right\|+\left\|\left|x\right|^2 \Delta f\right\|\right).
}
\item Let $f \in H^2(2)$, then $\left|x\right|^2 \nabla f \in L^4(\R^2)$ and there exists $C>0$ such that
\eq{\label{weightnabla2L4}
\left\|\left|x\right|^2 \nabla f\right\|_{L^4} \leq C\left\|\left|x\right|^2 \nabla f\right\|^{1/2}\left(\left\|f\right\|^{1/2}+\left\|\left|x\right|\nabla f\right\|^{1/2}+\left\|\left|x\right|^2\Delta f\right\|^{1/2}\right).
}
\end{enumerate}}
\textbf{Proof: }The inequality (\ref{weightLp}) comes directly from Hölder's inequality. To prove the inequality (\ref{weightnabla2}), we show by a simple calculation that, for every $j,k \in \left\{1,2\right\}$,
\eq{\label{weight2.0}
\left\|\left|x\right|^2 \partial_j \partial_k f\right\|^2\leq C\left(\left\|f\right\|^2+\left\| \left|x\right| \nabla f\right\|^2+\left\|\left|x\right|^2 \Delta f \right\|^2\right).
}
Indeed, we notice that
\eq{\label{weight2.1}
\left|x\right|^2 \partial_j \partial_k f=\partial_j \partial_k\left(\left|x\right|^2 f\right)-2\delta_{j,k} f -2x_j\partial_k f-2x_k\partial_j f,
}
and furthermore
\eq{\label{weight2.2}
\left\|\partial_j \partial_k\left(\left|x\right|^2 f\right)\right\|^2 \leq C\left\|\Delta\left(\left|x\right|^2 f\right)\right\|^2 \leq C\left(\left\|f\right\|^2+\left\| \left|x\right| \nabla f\right\|^2+\left\|\left|x\right|^2 \Delta f \right\|^2\right).
}
Combining (\ref{weight2.1}) and (\ref{weight2.2}) we get the inequality (\ref{weight2.0}). 
\vspace{0.5cm}\\
To obtain (\ref{weightnabla2L4}), we use Gagliardo-Niremberg's inequality as follows:
\aligne{
\left\|\left|x\right|^2 \nabla f\right\|_{L^4} &\leq C \left\|\left|x\right|^2 \nabla f\right\|^{1/2} \left\|\nabla\left(\left|x\right|^2 \nabla f\right)\right\|^{1/2}\\
& \leq C\left\|\left|x\right|^2 \nabla f\right\|^{1/2}\left(\left\|\left|x\right|\nabla f\right\|^{1/2}+\left\|\left|x\right|^2\nabla^2 f\right\|^{1/2}\right),
}
and consequently inequality (\ref{weightnabla2}) implies (\ref{weightnabla2L4}).
\begin{flushright}
$\square$
\end{flushright}
\paragraph{Biot-Savart law:}
Let $w$ be a real function defined on $\R^2$. The Bio-Savart law is a way to build a divergence free vector field $u$ such that $\curl u=w$. It is given by
\eq{\label{biotsavart}
u(x)=\displaystyle \frac{1}{2\pi}\int_{\R^2} \frac{\left(x-y\right)^\perp}{\left|x-y\right|^2}w(y) dy,
}
where $\left(x_1,x_2\right)^\perp=\left(-x_2,x_1\right)$.\\
\vspace{0.5cm}\\
The next two lemmas give estimates on the divergence free vector field $u$ obtained from $w$ via the Bio-Savart law.
\lem{\label{biots1}
Let $u$ be the divergence free vector field given by (\ref{biotsavart}).
\begin{enumerate}
\item Assume that $1<p<2<q<\infty$ and $\frac{1}{q}=\frac{1}{p}-\frac{1}{2}$. If $w\in L^p(\R^2)$, then $u \in L^q(\R^2)^2$ and there exists $C>0$ such that
\eq{\label{biotspq}
\left\|u\right\|_{L^q} \leq C\left\|w\right\|_{L^p}.
}
\item Assume that $1\leq p<2<q\leq \infty$, and define $\alpha \in (0,1)$ by the relation $\frac{1}{2}=\frac{\alpha}{p}+\frac{1-\alpha}{q}$. If $w\in L^p(\R^2)\cap L^q(\R^2)$, then $u\in L^\infty(\R^2)^2$ and there exists $C>0$ such that
\eq{
\left\| u \right\|_{L^\infty}\leq C\left\|w\right\|^\alpha_{L^p}\left\|w\right\|^{1-\alpha}_{L^q}.
}
\item Assume that $1<p<\infty$. If $w\in L^p(\R^2)$, then $\nabla u\in L^p(\R^2)^4$ and there exists $C>0$ such that
\eq{\label{biotsnabla}
\left\|\nabla u \right\|_{L^p}\leq C\left\|w\right\|_{L^p}.
}
\end{enumerate}
In addition, $\dv u=0$ and $\curl u =w$.
}
We refer to \cite{gallaywayne02} for the proof of this lemma.
\lem{\label{biots1bis}
Let $u$ be the divergence free vector field given by (\ref{biotsavart}).
\begin{enumerate}
\item If $w\in L^2(2)$, then $u\in L^4(\R^2)^2$ and there exists $C>0$ such that 
\eq{\label{biots4}
\left\|u\right\|_{L^4} \leq C\left\|w\right\|_{L^2(2)}.
}
\item If $w\in L^2(2)\cap H^1(\R^2)$, then $u\in L^\infty(\R^2)^2$ and there exists $C>0$ such that
\eq{\label{biotsinfini}
\left\|u\right\|_{L^\infty}\leq C \left\|w\right\|^{1/2}_{H^1}\left\|w\right\|^{1/2}_{L^2(2)}.
}
\item Let $s\in \R$. If \hspace{0.2cm}$\left(-\Delta\right)^{\frac{s-1}{2}} w\in L^2(\R^2)$ \hspace{0.2cm} for $s \in \R$, then \hspace{0.2cm}$\left(-\Delta\right)^{s/2}u\in L^2(\R^2)^2$\hspace{0.2cm} and there exists $C>0$ such that
\eq{
\left\|\left(-\Delta\right)^{s/2} u \right\|\leq C\left\|\left(-\Delta\right)^{\frac{s-1}{2}}w\right\|.
}
\item Let $s\in \R$. If $w\in H^s(\R^2)$, then $\nabla u\in H^s(\R^2)^4$ and there exists $C>0$ such that
\eq{\label{biotsHs}
\left\|\nabla u\right\|_{H^s} \leq C\left\|w\right\|_{H^s}.
}
\end{enumerate}
}
The proof of the two first inequalities are shown in \cite{jaffal11}. The two other inequalities are obvious when using Fourier variables. The next lemma is useful to get energy estimates in weighted Sobolev spaces for solutions of $(\ref{g3W})$. For a vector field $u$, we set
\cent{
$\quad \left|\nabla^3 u\right|^2= \displaystyle \sum^{2}_{i,j,k,l=1} \left(\partial_j \partial_k \partial_l u_i\right)^2$.
}
\lem{\label{biots2}
Let $w \in L^2(\R^2)$ and $u$ be the divergence free vector given by (\ref{biotsavart}).
\begin{enumerate}[1.]
\item If $w\in H^1(1)$, then $ \nabla^2 u\in L^2(1)$ and there exists $C>0$ such that
\eq{\label{biots3nabla2}
\left\| \nabla^2 u \right\|_{L^2(1)}\leq C\left(\left\|w\right\|_{H^1}+\left\|\left|x\right|\nabla w \right\|\right).
}
\item If $w\in H^2(1)$, then $\left|x\right|\nabla^2 u\in L^4(\R^2)$ and there exists $C>0$ such that
\eq{\label{biots3nabla2L4}
\left\|\left|x\right| \nabla^2 u \right\|_{L^4}\leq C\left(\left\|w\right\|+\left\|\left|x\right|\nabla w \right\|\right)^{1/2}\left(\left\|\nabla w\right\|+\left\|\left|x\right|\Delta w \right\|\right)^{1/2}.
}
\item If $w\in H^2(2)$, then $u\in \left|x\right|^2 \nabla^3 u\in L^2(\R^2)$ and there exists $C>0$ such that
\eq{\label{biots3nabla3}
\left\|\left|x\right|^2 \nabla^3 u \right\|\leq C \left(\left\|w\right\|
+\left\|\left|x\right| \nabla w\right\|+\left\|\left|x\right|^2 \Delta w\right\|\right).}
\item If $w\in L^2(1)$ and $\displaystyle \int_{\R^2}w(x)dx=0$, then $u \in H^1(1)$ and there exists a positive constant $C$ such that
\eq{\label{int0nabla}\left\|u\right\|+\left\|\left|x\right|\nabla u\right\| \leq C \left\|\left|x\right| w\right\|.}
\item If $w \in H^1(2)$ and $\displaystyle \int_{\R^2}w(x)dx=0$, then $\left|x\right|^2 \nabla^2 u\in L^2(\R^2)$ and there exists a positive constant $C$ such that
\eq{\label{int0nabla2}
\left\|\left|x\right|^2\nabla^2 u\right\| \leq C\left\|w\right\|_{H^1(2)}.
}
\end{enumerate}
}
\textbf{Proof: }Let us show the inequality (\ref{biots3nabla2}). Let $w$ belong to $H^1(1)$ and $u$ be the divergence free vector field obtained via the Biot-Savart law. From the inequality \ref{biotsnabla} of Lemma \ref{biots1}, we obtain
\eq{\label{proofbiots2.00}
\left\|\nabla^2 u\right\|_{L^2} \leq C \left\|\nabla w\right\|_{L^2}.}
Since the divergence of $u$ vanishes and since we are in dimension 2, it is enough to show the inequality 
\eq{\label{proofbiots2.01}
\left\|x_i \partial^2_j u_k\right\| \leq C\left(\left\|w\right\|+\left\|\left|x\right|\nabla w \right\|\right),
}
where $i,j,k \in \left\{1,2\right\}$. 
\vspace{0.5cm}\\
We omit $k$ that doesn't appear in the following calculations. One has
\aligne{
\left\|\left|x_i\right| \partial^2_j u\right\|^2 &= \left(2\pi\right)^2\displaystyle{\int_{\R^2}} \left|\partial_i \left(\xi^2_j \widehat{u}\right)\right|^2 d\xi\\
&\leq C\displaystyle{\int_{\R^2}} \left|\xi_j \widehat{u}\right|^2 d\xi+C\displaystyle{\int_{\R^2}} \left| \xi^2_j \partial_i\widehat{u}\right|^2 d\xi\\
&\leq C\left\|\nabla u\right\|^2+C\displaystyle{\int_{\R^2}} \left| \mathcal F \left(\Delta\left( x_i u\right)\right)\right|^2 d\xi\\
&\leq C\left\|\nabla u\right\|^2+C\left\|\left|x\right|\Delta u\right\|^2.
}
Using the inequality (\ref{biotsnabla}) of Lemma \ref{biots1} with $p=2$ and remarking that $\partial_1 w =\Delta u_2$ and $\partial_2 w =\Delta u_1$, we obtain (\ref{proofbiots2.01}). Combining it with the inequality (\ref{proofbiots2.00}), we get (\ref{biots3nabla2}).
\vspace{0.5cm}\\
The inequality (\ref{biots3nabla2L4}) is a direct consequence of (\ref{biots3nabla2}) and Gagliardo-Niremberg inequality. Indeed, one has
\aligne{
\left\|x_i\partial^2_j u\right\|_{L^4} &\leq C \left\|x_i\partial^2_j u\right\|^{1/2}\left\|\nabla\left(x_i\partial^2_j u\right)\right\|^{1/2}\\
&\leq C \left\|x_i\partial^2_j u\right\|^{1/2}\left(\left\|\partial^2_j u\right\|+\left\|x_i \partial^2_j \nabla u\right\|\right)^{1/2}.\\
}
Furthermore, the inequalities (\ref{biotsnabla}) and (\ref{biots3nabla2}) yield
\aligne{
\arr{ll}{
\left\|x_i\partial^2_j u\right\|_{L^4} &\leq C\left(\left\|w\right\|+\left\|\left|x\right|\nabla w\right\|\right)^{1/2}\left(\left\|\nabla w\right\|+\left\|x_i \nabla^2 w\right\|\right)^{1/2}.\\
}
}
Making the same computations than the ones we made to establish (\ref{biots3nabla2}), we obtain
\cent{
$\left\|x_i \nabla^2 w\right\| \leq C\left(\left\|\nabla w\right\|+ \left\|\left|x\right|\Delta w\right\|\right)$,
}
which gives
\aligne{
\left\|x_i\partial^2_j u\right\|_{L^4} &\leq C\left(\left\|w\right\|+\left\|\left|x\right|\nabla w\right\|\right)^{1/2}\left(\left\|\nabla w\right\|+\left\|\left|x\right| \Delta w\right\|\right)^{1/2},\\
}
and the inequality (\ref{biots3nabla2L4}) comes when summing for $i\in\{1,2\}$.
\vspace{0.5cm}\\
In order to get the inequality (\ref{biots3nabla3}), it suffices to obtain it for $\left|x\right|^2 \partial_j \partial^2_k u$, where $j,k \in \left\{1,2\right\}$. One has
\aligne{
\left\|\left|x\right|^2 \partial_j \partial^2_k u\right\|^2 &=\displaystyle \left(2\pi\right)^2\int_{\R^2}\left| \Delta \left( \xi_j \xi^2_k \widehat{u}\right) \right|^2d\xi\\
&\leq C\left(\displaystyle{\int_{\R^2}}\left| \left|\xi\right|^2\left(\xi_j+\xi_k\right)\Delta \widehat{u} \right|^2d\xi+\displaystyle{\int_{\R^2}}\left| \left(\xi_j+\xi_k\right)\widehat{u} \right|^2d\xi+\displaystyle{\int_{\R^2}}\left|\left|\xi\right|^2\nabla \widehat{u} \right|^2d\xi\right)\\
&\leq C\left(\left\|\nabla \Delta\left( \left|x\right|^2 u\right) \right\|^2+\left\|\nabla u\right\|^2+\sum^{2}_{i=1}\left\|\Delta \left(x_i u\right)\right\|^2\right)\\
&\leq C\left(\left\|\left|x\right|^2 \nabla \Delta u \right\|^2+\left\|\nabla u\right\|^2+\left\| \left|x\right| \Delta u\right\|^2\right)\\
&\leq C\left(\left\|\left|x\right|^2 \nabla^2 w \right\|^2+\left\|w\right\|^2+\left\| \left|x\right| \nabla w\right\|^2\right).\\
}
Applying the inequality (\ref{weightnabla2}), we get (\ref{biots3nabla3}). The proof of the inequality (\ref{int0nabla}) is made in two steps. It is shown in \cite{jaffal11} that
\eq{\label{proofbiots2.1}\left\|u\right\|\leq C \left\|\left|x\right| w\right\|.}
To finish the proof of the inequality (\ref{int0nabla}), we notice that
\eq{\label{proofbiots2.2}
\displaystyle \left\|\left|x\right|w\right\|^2 = \left\|\left|x\right|\partial_1 u_2\right\|^2+\left\|\left|x\right|\partial_2 u_1\right\|^2-2 \displaystyle \int_{\R^2}\left|x\right|^2\partial_1 u_2 \partial_2 u_1 dx.
}
Integrating by parts, one gets
\aligne{
-2 \displaystyle \int_{\R^2}\left|x\right|^2\partial_1 u_2 \partial_2 u_1 dx &=  \displaystyle \int_{\R^2}\left|x\right|^2 u_2 \partial_1 \partial_2 u_1 dx+2\int_{\R^2}x_1 u_2 \partial_2 u_1dx \\
&\hspace{1cm}+\int_{\R^2}\left|x\right|^2\partial_1 \partial_2 u_2  u_1 dx+2\int_{\R^2}x_2 \partial_1 u_2 u_1dx.
}
Using the divergence free property of $u$ and integrating by parts, we have
\cent{
$\displaystyle -2 \displaystyle \int_{\R^2}\left|x\right|^2\partial_1 u_2 \partial_2 u_1 dx =  \left\|\left|x\right| \partial_1 u_1 \right\|^2 + \left\|\left|x\right| \partial_2 u_2 \right\|^2+4\int_{\R^2} x_2 u_2 \partial_2 u_2 dx+4\int_{\R^2}x_1 \partial_1 u_1 u_1dx$.
}
Finally, integrating again by parts, we get
\cent{$
-2 \displaystyle \int_{\R^2}\left|x\right|^2\partial_1 u_2 \partial_2 u_1 dx =\left\|\left|x\right| \partial_1 u_1 \right\|^2 + \left\|\left|x\right| \partial_2 u_2 \right\|^2-2\left\|u\right\|^2.
$}
Thus, going back to (\ref{proofbiots2.2}), one has
\cent{$
\left\|\left|x\right| \nabla u\right\|^2 = \left\|\left|x\right| w\right\|^2+2\left\|u\right\|^2.
$}
Combining this equality with (\ref{proofbiots2.1}), we get the inequality (\ref{int0nabla}). The inequality (\ref{int0nabla2}) is obtained in the same way.
\begin{flushright}
$\square$
\end{flushright}
\section{\label{secapprox}Approximate solutions}
In this section, we introduce a "regularized" system of equations, whose solutions are more regular than the solutions of (\ref{g3w}). Actually, this new system is very close to (\ref{g3w}), and is obtained by adding the small term $\varepsilon \Delta^2 w$ to (\ref{g3w}). Here, the positive constant $\varepsilon$ is supposed to be small and is devoted to tend to $0$. Adding this term, we are able to prove the existence of solutions to the regularized system via a semi-group method. The presence of the term $u.\nabla \Delta w$ would not let us obtain solutions to (\ref{g3w}) by a semi-group method because of the too high degree of derivatives in this term compared to the linear term $\Delta w$. We introduce now the following regularized system of equations:
\eq{
\label{g3e}
\arr{l}{
\partial_t\left(w_\varepsilon-\alpha_1 \Delta w_\varepsilon\right) +\varepsilon \Delta^2 w_\varepsilon- \Delta w_\varepsilon + u_\varepsilon \nabla \left(w_\varepsilon-\alpha_1 \Delta w_\varepsilon\right)\\
\hspace{4cm}-\beta \dv\left(\left|A_\varepsilon\right|^2\nabla w_\varepsilon\right)-\beta \dv\left(\nabla\left(\left|A_\varepsilon\right|^2\right)\wedge A_\varepsilon\right) =0,\\
w_{\varepsilon\left|t=0\right.}=w_0 \in H^2(2),
}
}
where $A_\varepsilon = \nabla u_{\varepsilon}+\left(\nabla u_{\varepsilon}\right)^t$.\\
\vspace{0.5cm}\\
The aim of this section is to prove the following theorem.
\theo{\label{theo2}
Let $w_0 \in H^2(2)$. For all $\varepsilon>0$, there exists $t_\varepsilon >0$ and a unique solution $w_\varepsilon$ of the system $(\ref{g3e})$ such that
\cent{
$w_\varepsilon \in C^1\left(\left(0,t_\varepsilon\right),H^1(2)\right)\cap C^0\left(\left[0,t_\varepsilon\right),H^2(2)\right)\cap C^0\left(\left(0,t_\varepsilon\right),H^3(2)\right)$.
}
}
\textbf{Proof: }
First of all, we introduce the change of variable $\tilde{x}=\gamma x$, where $\gamma$ is a positive constant that is close to $0$ and will be made more precise later. This is made in order to not have to consider restrictions on the size of $\alpha_1$. We note $v_\varepsilon (x)=w_\varepsilon(x/\gamma)$. The system $(\ref{g3e})$ provides a new system in $v_\varepsilon$, that we will solve in $H^2(2)$.
\eq{\label{g3ve}
\arr{l}{
\partial_t\left(v_\varepsilon-\alpha_1 \gamma ^2\Delta v_\varepsilon\right) +\varepsilon \gamma^4 \Delta^2 v_\varepsilon- \gamma^2 \Delta v_\varepsilon + \gamma u_\varepsilon. \nabla \left(v_\varepsilon-\alpha_1 \gamma^2\Delta v_\varepsilon\right)\\
\hspace{3cm}-\beta  \gamma \nabla\left(\left|A_\varepsilon\right|^2\right).\nabla v_\varepsilon-\beta  \gamma^2 \left|A_\varepsilon\right|^2\Delta v_\varepsilon-\beta \dv\left(\nabla\left(\left|A_\varepsilon\right|^2\right)\wedge A_\varepsilon\right) =0,\\
\\
v_{\varepsilon\left|t=0\right.}=w_0(x/\gamma)\in H^2(2).
}
}
Although there are terms involving $u_\varepsilon$ in this system, it is actually autonomous. In fact, one recover $w_\varepsilon$ from $v_\varepsilon$ and then recover $u_\varepsilon$ via the Biot-Savart law (\ref{biotsavart}) applied to $w_\varepsilon$. We set
\cent{
$z_\varepsilon(x)=q(x)v_\varepsilon (x)$,
}
where $q(x)=\left(1+\left|x\right|^2\right)$.\\
\vspace{0.5cm}\\
To show the existence of a solution in $H^2(2)$ to the system $(\ref{g3ve})$, we are reduced to show that there exists a solution in $H^2(\R^2)$ of the system
\eq{\label{g3ze}
\arr{l}{
\partial_t\left(z_\varepsilon-\gamma^2\alpha_1 \Delta z_\varepsilon-\alpha_1 \gamma^2 q \Delta q^{-1} z_\varepsilon -2 \gamma^2 \alpha_1 q\nabla q^{-1}.\nabla z_\varepsilon\right) +\varepsilon \gamma^4\Delta^2 z_\varepsilon  =F\left(z_\varepsilon \right),\\
z_{\varepsilon\left|t=0\right.}=q w_0(x/\gamma) \in H^2(\R^2),
}
}
where
\eq{
\arr{l}{F(z_\varepsilon)= -\varepsilon \gamma^4 q \Delta ^2 \left(q^{-1} z_\varepsilon\right)+\gamma^2 q\Delta \left(q^{-1} z_\varepsilon\right)-\gamma q u_\varepsilon \nabla \left(q^{-1} z_\varepsilon -\gamma^2 \alpha_1 \Delta\left(q^{-1} z_\varepsilon\right)\right)\\
\\
\hspace{1.5cm}+\beta \gamma q \nabla \left(\left|A_\varepsilon\right|^2\right).\nabla\left(q^{-1} z_\varepsilon\right)+\beta \gamma^2 q \left|A_\varepsilon\right|^2\Delta\left(q^{-1} z_\varepsilon\right)+\beta q \dv \left(\nabla\left(\left|A_\varepsilon\right|^2\right)\wedge A_\varepsilon\right).
}
}
We define the two linear operators $B:D(B)=H^1(\R^2) \rightarrow H^{-1}(\R^2)$ and $D:D(D)=L^2(\R^2)\rightarrow H^{-1}(\R^2)$ as follows:
\cent{$
\arr{l}{
B(z)= \alpha_1 \gamma^2 \Delta z +\alpha_1\gamma^2 q\Delta q^{-1} z,\\
\\
D(z)=2\alpha_1 \gamma^2 q\nabla q^{-1} .\nabla z.
}
$
}
Via Lax-Milgram theorem, it is easy to show that $A=\left(I-B-D\right)$ is invertible. We define the bilinear form on $H^1(\R^2)$ 
\aligne{
a(u,v)= \left(u,v\right)_{L^2} + \alpha_1 \gamma^2 \left(\nabla u,\nabla v\right)_{L^2}-\alpha_1 \gamma^2 \left(q\Delta q^{-1} u, v\right)_{L^2}-2\alpha_1 \gamma^2 \left(q\nabla q^{-1}.\nabla u,v\right)_{L^2}.
}
We notice that $a$ is obviously coninuous on $H^1(\R^2)\times H^1(\R^2)$. Using the fact that $q\Delta q^{-1}$ and $q\nabla q^{-1}$ are bounded on $\R^2$, one has, for all $u,v \in H^1(\R^2)$,
\aligne{
\left|a(u,v)\right| \leq C(\alpha_1,\gamma) \left\|u\right\|_{H^1}\left\|v\right\|_{H^1},
}
where $C(\alpha_1,\gamma)$ is a positive constant depending on $\alpha_1$ and $\gamma$.\\
\vspace{0.5cm}\\
We show now that $a$ is coercive. Via an integration by parts, we get
\aligne{
a(u,u) = \left\|u\right\|^2+\alpha_1 \gamma^2 \left\|\nabla u\right\|^2-\alpha_1 \gamma^2 \int_{\R^2} q\Delta q^{-1}\left|u\right|^2 dx+\alpha_1 \gamma^2 \int_{\R^2} \dv \left(q\nabla q^{-1}\right)\left|u\right|^2 dx.
}
Due to the boundedness of $q\Delta q^{-1}$ and $\dv \left(q\nabla q^{-1}\right)$, there exists $C>0$ such that
\aligne{
a(u,u) \geq \left(1-\alpha_1 \gamma^2 C\right)\left\|u\right\|^2+\alpha_1 \gamma^2 \left\|\nabla u\right\|^2.
}
If we take $\gamma$ sufficiently small, the bilinear form $a$ is both continuous and coercive on $H^1(\R^2)$. From the Lax-Milgram theorem, we conclude that for all $f \in H^{-1}(\R^2)$ there exists $u \in H^{1}(\R^2)$ such that
\eq{\label{lax}
a(u,v)=\left\langle f,v\right\rangle_{H^{-1}\times H^1} \quad \hbox{for all}\quad v \in H^1(\R^2),
}
and consequently $\left(I-B-D\right)^{-1}$ is defined from $H^{-1}(\R^2)$ to $H^1(\R^2)$. We define $A : D(A)=H^3(\R^2)\rightarrow H^1(\R^2)$ the linear differential operator on $H^1(\R^2)$
\cent{
$A = \varepsilon \gamma^4 \left(I-B-D\right)^{-1} \Delta^2$.
}
We rewrite the system $(\ref{g3ze})$ as follows:
\eq{\label{g3z2}
\arr{l}{
\partial_t z_\varepsilon+A \left(z_\varepsilon\right)  =\widetilde{F}\left(z_\varepsilon \right),\\
z_{\varepsilon\left|t=0\right.}=q w_0(x/\gamma) \in H^2(\R^2),
}
}
where $\widetilde{F}(z_\varepsilon)=\left(I-B-D\right)^{-1}F(z_\varepsilon)$.\\
\vspace{0.5cm}\\
To finish the proof of this theorem, we show that the operator $A$ is sectorial on $H^1(\R^2)$, which is equivalent to the fact that $-A$ generates an analytic semigroup on $H^1(\R^2)$. Then, we check that $\widetilde{F}$ is locally Lipschitz from bounded sets of a Sobolev space $H^s(\R^2)$ to $H^1(\R^2)$, where $1\leq s <3$. An easy computation leads to
\aligne{
A&= \varepsilon \gamma^4\left(I-B\right)^{-1}\Delta^2-\varepsilon \gamma^4 \left(I-B-D\right)^{-1}D\left(I-B\right)^{-1}\Delta^2
\\&= I + \varepsilon \gamma^4\left(I-B\right)^{-1}\Delta^2-I-\varepsilon \gamma^4 \left(I-B-D\right)^{-1}D\left(I-B\right)^{-1}\Delta^2\\
&=J+R,
}
where
\cent{
$\arr{l}{
J = I + \varepsilon \gamma^4 \left(I-B\right)^{-1}\Delta^2,\\
R =-I-\varepsilon \gamma^4 \left(I-B-D\right)^{-1}D\left(I-B\right)^{-1}\Delta^2.
}$
}
Using the same method as the one used to invert $\left(I-B-D\right)$, one can invert $\left(I-B\right)$ and define $\left(I-B\right)^{-1}$ from $H^{-1}(\R^2)$ to $H^1(\R^2)$. Consequently, $J$ is well defined from $H^3(\R^2)$ to $H^1(\R^2)$. In the remaining of this proof, we will show that $-J$ generates an analytic semi-group on $H^1(\R^2)$ and then show that $R$ satisfies the conditions of \cite[Theorem 2.1 p. 81]{pazy83}. According to this result, it implies that $-A$ generates an analytic semi-group on $H^1(\R^2)$. In order to show that $J$ is sectorial on $H^1(\R^2)$, we associate it to a continuous and coercive bilinear form on $H^2(\R^2)\times H^2(\R^2)$. To this end, we define a $H^1$-scalar product which is suitable to $J$. Let us define, for $u,v \in H^1(\R^2)$, the bilinear form on $H^1$ given by
\cent{
$\left\langle u,v \right\rangle_{H^1} = \left(\left(1 - \alpha_1\gamma^2 q\Delta q^{-1}\right)u,v\right)_{L^2} + \alpha_1\gamma^2\left( \nabla u,\nabla v\right)_{L^2}$.
}
If $\gamma$ is sufficiently small compared to $\alpha_1$, then $\left\langle .,.\right\rangle_{H^1}$ is a scalar product on $H^1(\R^2)$. Furthermore, for $u \in H^2( \R^2)$ and $v\in H^1(\R^2)$, one has
\cent{
$\left\langle u,v \right\rangle_{H^1} = \left(\left(I-B\right) u,v\right)_{L^2}$.
}
We define, using this scalar product, the bilinear form $j$ on $H^2(\R^2)\times H^2(\R^2)$ associated to $J$ by the formula
\cent{
$
j(u,v) = \left\langle u,v\right\rangle_{H^1} + \varepsilon \gamma^4 \left(\Delta u,\Delta v\right)_{L^2}.
$}
A short computation shows that, for $u \in H^3(\R^2)$ and $v\in H^2(\R^2)$, one has
\eq{\label{jH3}
j(u,v)=\left\langle J u , v\right\rangle_{H^1}.
}
Furthermore, if $\gamma$ is small enough, using the definition of $\left\langle .,.\right\rangle_{H^1}$ and $j$, we see that there exists $C(\alpha_1, \varepsilon, \gamma)>0$ such that, for all $u,v \in H^2(\R^2)$,
\cent{
$
j(u,v) \leq C(\alpha_1, \varepsilon, \gamma)\left\|u\right\|_{H^2}\left\|v\right\|_{H^2}.
$
}
Besides, it is simple to check that, if $\gamma$ is mall enough, there exists $C(\alpha_1, \gamma, \varepsilon)>0$ such that, for all $u \in H^2(\R^2)$,
\cent{
$
j(u,u) \geq C(\alpha_1, \gamma, \varepsilon)\left\|u\right\|^2_{H^2}.
$
}
The bilinear form $j$ is thus continuous and coercive on $H^2(\R^2)$ and the operator $J$ is consequently sectorial on $H^1(\R^2)$.  Additionally, The linear operator $R$ is defined from $H^2(\R^2)$ to $H^1(\R^2)$, and one can check that there exists $C(\alpha_1, \gamma, \varepsilon)>0$ such that, for all $u\in H^3(\R^2)$,
\eq{\label{RH3}
\left\|R u\right\|_{H^1} \leq C(\alpha_1, \gamma, \varepsilon)\left\|u\right\|_{H^2}.
}
Applying the equality (\ref{jH3}) to $u \in H^3(\R^2)$, we get
\cent{
$
j(u,u)=\left\langle J u, u\right\rangle_{H^1},
$ for all $u\in H^3(\R^2)$.
}
Because $j$ is coercive on $H^2$, we obtain, via Cauchy-Schwartz inequality,
\cent{
$
\left\|u\right\|^2_{H^2}\leq C(\alpha_1, \gamma, \varepsilon)\left\| J u\right\|_{H^1}\left\| u\right\|_{H^1},
$ for all $u\in H^3(\R^2)$.
}
Going back to (\ref{RH3}), the following property holds
\cent{
$
\left\|R u\right\|^2_{H^1} \leq C(\alpha_1, \gamma, \varepsilon)\left\| J u\right\|_{H^1}\left\| u\right\|_{H^1},
$ for all $u\in H^3(\R^2)$.
}
In particular, the Young inequality yields, for all $\delta>0$,
\cent{
$
\left\|R u\right\|^2_{H^1} \leq \delta \left\| J u\right\|^2_{H^1}+C(\alpha_1, \gamma, \varepsilon)\left\| u\right\|^2_{H^1},
$ for all $u\in H^3(\R^2)$.
}
By a classical result that we can find in \cite{henry81}, $-A$ is thus the generator of an analytic semigroup on $H^1(\R^2)$.\\
\vspace{0.5cm}\\
Lastly, it is easy to check that $\widetilde{F}$ is Lipschitzian from the bounded sets of $H^2(\R^2)$ into $H^1(\R^2)$. Combining several results from \cite[chapter 3]{henry81} and \cite[section 6.3]{pazy83}, we conclude that there exists $t_\varepsilon>0$ and a unique solution $z_\varepsilon \in C^1\left(\left(0,t_\varepsilon\right),H^1(\R^2)\right)\cap C^0\left(\left[0,t_\varepsilon\right),H^2(\R^2)\right)\cap C^0\left(\left(0,t_\varepsilon\right),H^3(\R^2)\right)$ of the system $(\ref{g3ze})$. Thus, there exists a unique solution $w_\varepsilon \in C^1\left(\left(0,t_\varepsilon\right),H^1(2)\right)\cap C^0\left(\left[0,t_\varepsilon\right),H^2(2)\right)\cap C^0\left(\left(0,t_\varepsilon\right),H^3(2)\right)$ to the system $(\ref{g3e})$.
\begin{flushright}
$\square$
\end{flushright}
\section{\label{secenergy}Energy estimates}
In this section, we perform energy estimates on the regularized solutions of the third grade fluids equations in the weighted space $H^2(2)$. These estimates are independent of $\varepsilon$ and allows us, in Section \ref{secdem}, to pass to the limit when $\varepsilon$ tends to $0$. Thus, we consider the solution $w_\varepsilon(t,x)$ of (\ref{g3e}). Let $T$, $T\geq 1$ be a fixed positive constant and $\tau_0 = \log(T)$. We define $W_\varepsilon(\tau,X)$, obtained from $w_\varepsilon$ by the change of variables (\ref{scaledv1}) and (\ref{scaledv2}). A short computation shows that $W_\varepsilon$ satisfies the system
\eq{
\label{g3We}
\arr{l}{
\displaystyle \partial_\tau\left(W_\varepsilon-\alpha_1 e^{-\tau}\Delta W_\varepsilon\right)+\varepsilon e^{-\tau}\Delta^2 W_\varepsilon- \mathcal L (W_\varepsilon)+ U_\varepsilon.\nabla \left(W_\varepsilon-\alpha_1 e^{-\tau}\Delta W_\varepsilon\right)+\alpha_1 e^{-\tau}\Delta W_\varepsilon\\
\hspace{1.6cm}\displaystyle +\alpha_1 e^{-\tau}\frac{X}{2}.\nabla \Delta W_\varepsilon -\beta e^{-2\tau}\dv\left(\left|A_\varepsilon\right|^2\nabla W_\varepsilon\right)-\beta e^{-2\tau}\dv\left(\nabla\left(\left|A_\varepsilon\right|^2\right)\wedge A_\varepsilon\right)=0,\\
\dv U_\varepsilon =0,\\
W_{\varepsilon \left|\tau=\tau_0\right.}=W_0,
}
}
where $\tau_0=\log(T)$, $U_\varepsilon$ is obtained from $W_\varepsilon$ via the Biot-Savart law (\ref{biotsavart}), $A_\varepsilon = \nabla U_\varepsilon + \left(\nabla U_\varepsilon\right)^t$ and we recall that
\cent{
$
\mathcal L(W_\varepsilon) = \Delta W_\varepsilon + W_\varepsilon + \frac{X}{2}.\nabla W_\varepsilon.
$}
By theorem \ref{theo2}, it is clear that there exists $\tau_\varepsilon>\tau_0$ such that
\cent{
$
W_\varepsilon\in C^1\left(\left(\tau_0,\tau_\varepsilon\right),H^1(2)\right)\cap C^0\left(\left(\tau_0,\tau_\varepsilon\right),H^3(2)\right).
$
}
We assume also that the initial datum $W_0 \in H^2(2)$ satisfies the assumption (\ref{cond1}) of Theorem \ref{theo1}, for some $\gamma >0$.
Let $\eta=\displaystyle \int_{\R^2} W_0(X) dX$, we write the following decompositions
\eq{\label{decomp}\arr{l}{W_\varepsilon=\eta G+f_\varepsilon,\\
U_\varepsilon=\eta V+K_\varepsilon,}
}
where $G$ is the Oseen vortex sheet defined by (\ref{oseen}) and $V$ is the divergence free vector field obtained from $G$ via the Biot-Savart law (\ref{biotsavart}). Using the fact that $\mathcal L (G)=0$, one has the equality
\eq{
\label{g3fe}
\arr{l}{
\partial_\tau\left(f_\varepsilon-\alpha_1 e^{-\tau}\Delta f_\varepsilon\right)+\varepsilon e^{-\tau} \Delta^2 f_\varepsilon-\mathcal{L}(f_\varepsilon)+ K_\varepsilon.\nabla\left(f_\varepsilon-\alpha_1 e^{-\tau} \Delta f_\varepsilon\right)\\
\\
\hspace{1cm} +\eta V.\nabla\left(f_\varepsilon-\alpha_1 e^{-\tau} \Delta f_\varepsilon\right)+\eta K_\varepsilon.\nabla\left(G-\alpha_1 e^{-\tau} \Delta G\right)+\alpha_1 e^{-\tau} \Delta f_\varepsilon\\
\\
\hspace{2cm}+\alpha_1 e^{-\tau}\frac{X}{2}.\nabla \Delta f_\varepsilon+\eta\alpha_1 e^{-\tau} \Delta G+\eta\alpha_1 e^{-\tau}\frac{X}{2}.\nabla \Delta G+\eta\varepsilon e^{-\tau} \Delta^2 G\\
\\
\hspace{2.5cm}-\beta e^{-2\tau}\dv\left(\left|A_\varepsilon\right|^2\nabla f_\varepsilon+\eta\left|A_\varepsilon\right|^2\nabla G\right)-\beta e^{-2\tau}\dv\left(\nabla\left(\left|A_\varepsilon\right|^2\right)\wedge A_\varepsilon\right)=0.
}
}
Let $M=M(\alpha_1,\beta)>2$ be a positive constant which will be made more precise later. Let $\tau^*_\varepsilon \in \left(\tau_0,\tau_\varepsilon\right]$ be the largest time (depending on $M$) such that, for all $\tau \in \left[\tau_0,\tau^*_\varepsilon\right)$, the following inequality holds
\eq{\label{cond1bis}
\arr{l}{\displaystyle \left\|W_\varepsilon(\tau)\right\|^2_{H^1}+\alpha_1 e^{-\tau}\left\|\Delta W_\varepsilon(\tau)\right\|^2_{L^2}+\left\|\left|X\right|^2 W_\varepsilon(\tau)\right\|^2_{L^2}+\alpha^2_1 e^{-2\tau}\left\|\left|X\right|^2\Delta W_\varepsilon(\tau)\right\|^2_{L^2}\leq   M\gamma \left(1-\theta\right)^6.}
}
To simplify the notations in the following computations, we assume that $0<\gamma \leq 1$ and we take $T$ sufficiently large so that $\displaystyle \frac{\alpha_1}{T}= \alpha_1 e^{-\tau_0} \leq 1$.
\vspace{0.5cm}\\
Since $W_\varepsilon \in C^0\left(\left[\tau_0,\tau_\varepsilon\right),H^2(2)\right)$ and the condition (\ref{cond1}) holds, $\tau^*_\varepsilon$ is well defined. Furthermore, there exists a positive constant $C$ independent of $W_0$ such that, for all $\tau \in \left[\tau_0,\tau^*_\varepsilon\right)$,
\eq{\label{cond2}
\eta^2+\left\|f_\varepsilon\right\|^2_{H^1}+\alpha_1 e^{-\tau}\left\|\Delta f_\varepsilon\right\|^2_{L^2}+\left\|\left|X\right|^2 f_\varepsilon\right\|^2_{L^2}+\alpha^2_1 e^{-2\tau}\left\|\left|X\right|^2 \Delta f_\varepsilon \right\|^2_{L^2}\leq C M \gamma \left(1-\theta\right)^6.
}
Indeed, using Cauchy-Schwartz inequality, we get
\aligne{
\eta &= \int_{\R^2} W_0(X) dX \\
&= \int_{\R^2} \frac{1+\left|X\right|^2}{1+\left|X\right|^2} W_0(X) dX\\
&\leq \left(\int_{\R^2} \frac{1}{\left(1+\left|X\right|^2\right)^2} dX\right)^{1/2}\left(\int_{\R^2} \left(1+\left|X\right|^2\right)^2 \left|W_0(X)\right|^2 dX\right)^{1/2}\\
&\leq C\left\|W_0\right\|_{L^2(2)}.
} 
Considering the decomposition (\ref{decomp}) and the smoothness of $G$, we obtain the inequality (\ref{cond2}).
\vspace{0.5cm}\\
To simplify the notations, in this section we write $f$ instead of $f_\varepsilon$, $W$ instead of $W_\varepsilon$, $U$ instead of $U_\varepsilon$ and $K$ instead of $K_\varepsilon$.
\vspace{0.5cm}\\
\indent The aim of this section is to show that the inequality (\ref{inetheo1}) of Theorem \ref{theo1} holds for the regularized solutions of the system (\ref{g3We}), provided that the condition (\ref{cond1}) is satisfied by $W_0$. To this end, we consider a fixed constant $\theta$ such that $0<\theta < 1$ which is twice the rate of convergence of $W$ to $\eta G$ in $H^2(2)$. In fact, we will show that, under the assumption (\ref{cond1}), the decaying of $f$ to $0$ in $H^2(2)$ is equivalent to $\displaystyle e^{-\frac{\theta\tau}{2}}$. As it is explained in the introduction of this paper, the spectrum of $ \mathcal L$ in $L^2(m)$ does not allow the rate of convergence to be better than $\displaystyle e^{-\frac{\tau}{2}}$.
\vspace{0.5cm}\\
In order to get the inequality (\ref{inetheo1}), we construct in this section an energy functional $E=E(\tau)$ such that, for every $\tau \in \left[\tau_0,\tau^*_\varepsilon\right)$,
\cent{
$E(\tau) \sim \left\|f(\tau)\right\|^2_{H^2(2)}$,
}
and there exists a positive constant $C=C(\alpha_1,\beta,\theta)$ such that, for all $\tau \in \left[\tau_0,\tau^*_\varepsilon\right)$,
\eq{\label{energyfinal}
\partial_ \tau E(\tau) + \theta E(\tau) \leq C \gamma e^{-\tau}.
}
This inequality will enable us to show that $\tau^*_\varepsilon = +\infty$ and obtain, by the application of Gronwall Lemma,
\cent{$
E(\tau) \leq C \gamma e^{-\theta \tau },
$ for all $\tau \in \left[\tau_0,+\infty\right)$.}
This functional is built as the sum of several intermediate energy functionals in various functions spaces, for which we perform convenient estimates.
\subsection{Estimates in $\dot{H}^{-\frac{1+\theta}{2}}$}
We start by performing an estimate of the solution of $(\ref{g3fe})$ in the homogeneous Sobolev space $\dot{H}^{-\frac{1+\theta}{2}}(\R^2)$. Combined with the other estimates, it will give us an estimate in the classical Sobolev space ${H}^{-\frac{1+\theta}{2}}(\R^2)$. 
The motivation to do this comes from the fact that the $H^1-$estimate that we will perform later (see Lemma \ref{lemE2}) makes the term $\left\| u\right\|^2_{L^2}$ appear on the right hand side of our $H^1-$energy inequality. In order to absorb this term, we look for an estimate in a Sobolev space of negative order. To this end, due to Lemma \ref{weight-} and the fact that $\displaystyle \int_{\R^2} f(X) dX=0$, for $\frac{3}{4}\leq s <1$, one can apply the operator $\left(-\Delta \right)^{-s}$ to the equality (\ref{g3fe}) and take the inner $L^2-$product of it with $\left(-\Delta \right)^{-s}f$. Through the computations that we will perform below, one can see that, in order to get the estimate (\ref{energyfinal}), we have to choose at least $s=\frac{1+\theta}{2}$. Actually, since we have to absorb terms coming from the non-linear part of (\ref{g3fe}), it is more convenient to take $\frac{1+\theta}{2}< s <1$, for instance $s=\frac{3+\theta}{4}$. In \cite{jaffal11}, the considered operator was $\left(-\Delta\right)^{-3/4}$, which implied the restriction $0< \theta <\frac{1}{2}$.
\vspace{0.5cm}\\
The next lemma summarizes the computations needed when applying $\left(-\Delta \right)^{-s}$ to (\ref{g3fe}) and taking the $L^2-$scalar product of it with $\left(-\Delta \right)^{-s}f$.
\lem{\label{weight-2}
Let $f\in H^3(2)$ such that $\displaystyle \int_{\R^2} f(X) dX = 0$, then, for all $\frac{1}{2}\leq s< 1$ the three following equalities hold.
\eq{\label{eqweight-2}\arr{l}{\left(\left(-\Delta\right)^{-s}\left(\frac{X}{2}.\nabla f\right),\left(-\Delta \right)^{-s} f\right)_{L^2}=-\left(s+\frac{1}{2}\right)\left\|\left(-\Delta\right)^{-s} f\right\|^2_{L^2},\\
\\
\left(\left(-\Delta\right)^{-s}\left(\mathcal L (f)\right),\left(-\Delta \right)^{-s} f\right)_{L^2}=-\left\|\left(-\Delta\right)^{\frac{1}{2}-s} f\right\|^2_{L^2}-\left(s-\frac{1}{2}\right)\left\|\left(-\Delta\right)^{-s} f\right\|^2_{L^2},\\
\\
\left(\left(-\Delta\right)^{-s}\left(\frac{X}{2}.\nabla \Delta f\right),\left(-\Delta \right)^{-s} f\right)_{L^2}=\left(s+1\right)\left\|\left(-\Delta\right)^{\frac{1}{2}-s} f\right\|^2_{L^2}.
}}
}
\textbf{Proof: }Using Fourier variables, it is easy to see that
\cent{
$\widehat{\frac{X}{2}.\nabla f}=- \widehat{f}-\frac{\xi}{2}.\nabla  \widehat{f}, \quad$ and $\quad
\widehat{\frac{X}{2}.\nabla \Delta f}=2\left|\xi\right|^2 \widehat{f}+\frac{\xi\left|\xi\right|^2}{2}\nabla  \widehat{f}.
$
}
The proof of this lemma is then obtained through the Plancherel formula and direct computations.
\begin{flushright}
$\square$
\end{flushright}
In order to obtain a priori estimates of $f$ in $\dot{H}^{-\frac{1+\theta}{2}}(\R^2)$, we define the functional
\cent{
$\displaystyle E_1(\tau)=\frac{1}{2} \left(\left\|\left(-\Delta \right)^{-\frac{3+\theta}{4}} f\right\|^2+\alpha_1 e^{-\tau} \left\|\left(-\Delta \right)^{-\frac{1+\theta}{4}} f\right\|^2\right)$.
}
The estimate in $\dot{H}^{-\frac{1+\theta}{2}}$ of $f$ under the condition $(\ref{cond1bis})$ is given in the next lemma.
\lem{\label{lemE1}Let $W\in C^1\left(\left(\tau_0,\tau_\varepsilon\right),H^1(2)\right)\cap C^0\left(\left(\tau_0,\tau_\varepsilon\right),H^3(2)\right)$ be the solution of $(\ref{g3We})$ satisfying the inequality (\ref{cond1bis}) for some $\gamma>0$. There exist $\gamma_0 >0$ and $T_0\geq 1$ such that if $T\geq T_0$ and $\gamma\leq \gamma_0$, then, for all $\tau\in \left[\tau_0,\tau^*_\varepsilon\right)$, $E_1$ satisfies the inequality
\eq{\label{E1f}
\arr{l}{\displaystyle \partial_\tau E_1+\theta E_1 +\left(1+\frac{1-\theta}{4}\alpha_1 e^{-\tau}\right)\left\|\left(-\Delta\right)^{-\frac{1+\theta}{4}}f\right\|^2\leq C M^3\gamma\left(1-\theta\right)^2 e^{-2\tau}\\
\hspace{5cm}\displaystyle+ CM\gamma \left(1-\theta\right)^2\left(\left\|f\right\|^2_{L^2(2)}+\left\|\nabla f\right\|^2+\alpha^2_1 e^{-2\tau}\left\|\Delta f\right\|^2_{L^2(1)}\right),
}
}
where $\theta$, $0<\theta<1$ is the fixed constant introduced  at the beginning of Section \ref{secenergy}.
}
\textbf{Proof: }Since $\displaystyle \int_{\R^2} f(X) dX =0$, according to Lemma \ref{weight-}, $\left(-\Delta\right)^{-\frac{3+\theta}{4}}f$ is well defined. Thus, we apply $\left(-\Delta\right)^{-\frac{3+\theta}{4}}$ to the equality (\ref{g3fe}) and we get
\eq{\label{eqE10}\arr{l}{
\partial_\tau\left(\left(-\Delta\right)^{-\frac{3+\theta}{4}} f +\alpha_1 e^{-\tau}\left(-\Delta\right)^{\frac{1-\theta}{4}} f\right)+\varepsilon e^{-\tau}\left(-\Delta\right)^{\frac{5-\theta}{4}} f-\left(-\Delta\right)^{-\frac{3+\theta}{4}}\left(\mathcal L (f)\right)\\
\hspace{3cm}-\alpha_1 e^{-\tau}\left(-\Delta\right)^{\frac{1-\theta}{4}} f+\alpha_1 e^{-\tau} \left(-\Delta\right)^{-\frac{3+\theta}{4}}\left(\frac{X}{2}.\nabla \Delta f\right)=H\left(\tau,G,f,W\right),
}}
where
\cent{
$
\arr{l}{
H\left(\tau,G,f,W\right)=\left(-\Delta\right)^{-\frac{3+\theta}{4}}\Biggl(-K.\nabla\left(f-\alpha_1 e^{-\tau} \Delta f\right)-\eta V.\nabla\left(f-\alpha_1 e^{-\tau} \Delta f\right)\\
\hspace{4cm}-\eta K.\nabla\left(G-\alpha_1 e^{-\tau} \Delta G\right)-\eta\alpha_1 e^{-\tau} \Delta G-\eta\alpha_1 e^{-\tau}\frac{X}{2}.\nabla \Delta G\\
\hspace{7cm}-\eta\varepsilon e^{-\tau} \Delta^2 G+\beta e^{-2\tau}\curl\dv\left(\left|A\right|^2A\right)\Biggr).
}
$
}
Taking the $L^2-$scalar product of (\ref{eqE10}) with $\left(-\Delta\right)^{-\frac{3+\theta}{4}}$ and taking into account the equalities
\cent{
$
\left(-\left(-\Delta \right)^{-\frac{3+\theta}{4}} \left(\mathcal L(f)\right),\left(-\Delta \right)^{-\frac{3+\theta}{4}}f\right)_{L^2}=\left\|\left(-\Delta \right)^{-\frac{1+\theta}{4}}f\right\|^2+\left(\frac{1+\theta}{4}\right)\left\|\left(-\Delta \right)^{-\frac{3+\theta}{4}} f\right\|^2,
$
}
and
\cent{
$
\left(\alpha_1 e^{-\tau}\left(-\Delta \right)^{-\frac{3+\theta}{4}} \left(\frac{X}{2}.\nabla \Delta f\right),\left(-\Delta \right)^{-\frac{3+\theta}{4}}f\right)_{L^2}=\left(\frac{7+\theta}{4}\right)\alpha_1 e^{-\tau}\left\|\left(-\Delta \right)^{-\frac{1+\theta}{4}} f\right\|^2,
$
} 
given by Lemma \ref{weight-2}, we obtain
\eq{\label{eqE1}\arr{l}{
\frac{1}{2}\partial_\tau\left(\left\|\left(-\Delta\right)^{-\frac{3+\theta}{4}}f\right\|^2+\alpha_1 e^{-\tau}\left\|\left(-\Delta\right)^{-\frac{1+\theta}{4}}f\right\|^2\right)+\varepsilon e^{-\tau} \left\|\left(-\Delta\right)^{-\frac{1-\theta}{4}}f\right\|^2\\
\hspace{2.5cm}+\left(\frac{1+\theta}{4}\right)\left\|\left(-\Delta\right)^{-\frac{3+\theta}{4}}f\right\|^2+\left(1+\left(\frac{1+\theta}{4}\right)\alpha_1 e^{-\tau}\right)\left\|\left(-\Delta\right)^{-\frac{1+\theta}{4}}f\right\|^2\\
\hspace{9cm}=\left(H\left(\tau,G,f\right),\left(-\Delta\right)^{-\frac{3+\theta}{4}} f\right)_{L^2}.}
}
Now, it remains to estimate the right hand side of (\ref{eqE1}), that we write as
\cent{$
\left(H\left(\tau,G,f\right),\left(-\Delta\right)^{-\frac{3+\theta}{4}} f\right)_{L^2}=I_1+I_2+I_3+I_4+I_5,
$}
where 
\cent{$
\arr{l}{
I_1=\left(\left(-\Delta\right)^{-\frac{3+\theta}{4}}\left(-K.\nabla\left(f-\alpha_1 e^{-\tau} \Delta f\right)\right),\left(-\Delta\right)^{-\frac{3+\theta}{4}} f\right)_{L^2},\\
I_2=\left(\left(-\Delta\right)^{-\frac{3+\theta}{4}}\left(-\eta V.\nabla\left(f-\alpha_1 e^{-\tau} \Delta f\right)\right),\left(-\Delta\right)^{-\frac{3+\theta}{4}} f\right)_{L^2},\\
I_3=\left(\left(-\Delta\right)^{-\frac{3+\theta}{4}}\left(-\eta K.\nabla\left(G-\alpha_1 e^{-\tau} \Delta G\right)\right),\left(-\Delta\right)^{-\frac{3+\theta}{4}} f\right)_{L^2},\\
I_4=\left(\left(-\Delta\right)^{-\frac{3+\theta}{4}}\left(-\eta\alpha_1 e^{-\tau} \Delta G-\eta\alpha_1 e^{-\tau}\frac{X}{2}.\nabla \Delta G-\eta\varepsilon e^{-\tau} \Delta^2 G\right),\left(-\Delta\right)^{-\frac{3+\theta}{4}} f\right)_{L^2}\\
\hspace{9cm}= \left(\left(-\Delta\right)^{-\frac{3+\theta}{4}} J,\left(-\Delta\right)^{-\frac{3+\theta}{4}} f\right)_{L^2},\\
I_5=\left(\left(-\Delta\right)^{-\frac{3+\theta}{4}}\left(\beta e^{-2\tau}\curl\dv\left(\left|A\right|^2A\right)\right),\left(-\Delta\right)^{-\frac{3+\theta}{4}} f\right)_{L^2}.
}
$}
The remaining of the proof of this lemma is devoted to the estimate of these terms. We recall that $\curl K=f$, $\curl V= G$ and $\curl U=W$. Since the divergence of $K$ vanishes, we obtain
\aligne{
I_1&=\left(\left(-\Delta\right)^{-\frac{3+\theta}{4}}\left(-\dv \left(K\left(f-\alpha_1 e^{-\tau} \Delta f\right)\right)\right),\left(-\Delta\right)^{-\frac{3+\theta}{4}} f\right)_{L^2}\\
&\leq \left\|\left(-\Delta\right)^{-\frac{3+\theta}{4}}\nabla \left(K\left(f-\alpha_1 e^{-\tau} \Delta f\right)\right)\right\|\left\|\left(-\Delta\right)^{-\frac{3+\theta}{4}} f\right\|.
}
Using the inequalities (\ref{weight-nabla}) of Lemma \ref{weight-} and (\ref{biotsinfini}) of Lemma \ref{biots1bis},  together with the Young and Hölder inequalities and the property (\ref{cond2}), we get
\eq{\label{E1I1}\arr{ll}{
I_1 &\displaystyle \leq \frac{C}{\left(1-\theta\right)^{3/2}}\left\|K\left(f-\alpha_1 e^{-\tau} \Delta f\right)\right\|_{L^2(1)}\left\|\left(-\Delta\right)^{-\frac{3+\theta}{4}} f\right\|\\
&\displaystyle\leq \frac{C}{\left(1-\theta\right)^{3/2}}\left\|K\right\|_{L^\infty}\left\|f-\alpha_1 e^{-\tau} \Delta f\right\|_{L^2(1)}\left\|\left(-\Delta\right)^{-\frac{3+\theta}{4}} f\right\|\\
&\displaystyle\leq \mu\left\|\left(-\Delta\right)^{-\frac{3+\theta}{4}} f\right\|^2 + \frac{C}{\mu \left(1-\theta\right)^{3}}\left\|f\right\|_{L^2(2)}\left\|f\right\|_{H^1}\left\|f-\alpha_1 e^{-\tau} \Delta f\right\|^2_{L^2(1)}\\
&\displaystyle\leq \mu\left\|\left(-\Delta\right)^{-\frac{3+\theta}{4}} f\right\|^2 + \frac{C M\gamma \left(1-\theta\right)^3}{\mu}\left(\left\|f\right\|^2_{L^2(1)}+\alpha^2_1 e^{-2\tau} \left\|\Delta f\right\|^2_{L^2(1)}\right),
}
}
where $\mu$ is a positive constant which is made more precise later.\\
\vspace{0.5cm}\\
Similar computations and the inequality (\ref{cond2}) give similar estimates for $I_2$. One has
\eq{\label{E1I2}
\displaystyle I_2 \leq \mu\left\|\left(-\Delta\right)^{-\frac{3+\theta}{4}} f\right\|^2 + \frac{CM\gamma \left(1-\theta\right)^3}{\mu}\left(\left\|f\right\|^2_{L^2(1)}+\alpha^2_1 e^{-2\tau} \left\|\Delta f\right\|^2_{L^2(1)}\right).
}
Likewise, we estimate the term $I_3$. Indeed, the same computations and the smoothness of $G$ yield
\aligne{
I_3 &\leq \mu\left\|\left(-\Delta\right)^{-\frac{3+\theta}{4}} f\right\|^2 + \frac{C\eta^2}{\mu \left(1-\theta\right)^{3}}\left\|f\right\|_{L^2(2)}\left\|f\right\|_{H^1}\left\|G-\alpha_1 e^{-\tau} \Delta G\right\|^2_{L^2(1)}\\
&\leq \mu\left\|\left(-\Delta\right)^{-\frac{3+\theta}{4}} f\right\|^2 + \frac{CM\gamma \left(1-\theta\right)^3}{\mu}\left(\left\|f\right\|^2_{L^2(2)}+\left\|f\right\|^2_{H^1}\right)\left(1+\alpha_1 e^{-\tau}\right).
}
Taking $T_0$ sufficiently large so that $\alpha_1 e^{-\tau} \leq 1$, we get
\eq{\label{E1I3}
\displaystyle I_3 \leq \mu\left\|\left(-\Delta\right)^{-\frac{3+\theta}{4}} f\right\|^2 + \frac{CM\gamma \left(1-\theta\right)^3}{\mu}\left(\left\|f\right\|^2_{L^2(2)}+ \left\|\nabla f\right\|^2\right).
}
Estimating $I_4$ is simple, because of the smoothness of $G$. We remark that $\displaystyle \int_{\R^2} J(X) dX=0$. Thus we can apply the inequality (\ref{weight-f}) of Lemma \ref{weight-} to obtain
\cent{
$\displaystyle \left\|\left(-\Delta\right)^{-\frac{3+\theta}{4}}J\right\| \leq \frac{C\left|\eta\right|e^{-\tau}\left\|G\right\|_{H^4(3)}}{\left(1-\theta\right)^{3/2}}$.
}
Using the above inequality and the smoothness of $G$, we can write
\eq{\label{E1I4}
\arr{ll}{
I_4 &\displaystyle \leq \frac{C\eta e^{-\tau}}{\left(1-\theta\right)^{3/2}}\left\|\left(-\Delta\right)^{-\frac{3+\theta}{4}} f\right\|\\
&\displaystyle \leq \mu \left\|\left(-\Delta\right)^{-\frac{3+\theta}{4}} f\right\|^2+ \frac{CM\gamma\left(1-\theta\right)^3 e^{-2\tau}}{\mu}.
}
}
It remains to estimate the term $I_5$. The inequality (\ref{weight-nabla}) of Lemma \ref{weight-} implies
\aligne{
I_5 &\leq \frac{C\beta e^{-2\tau}}{\left(1-\theta\right)^{3/2}}\left\|\nabla\left(\left|A\right|^2A\right)\right\|_{L^2(1)}\left\|\left(-\Delta \right)^{-\frac{3+\theta}{4}} f\right\|\\
&\leq \mu \left\|\left(-\Delta \right)^{-\frac{3+\theta}{4}} f\right\|^2+\frac{C\beta^2 e^{-4\tau}}{\mu\left(1-\theta\right)^{3}}\left\|\nabla\left(\left|A\right|^2A\right)\right\|^2_{L^2(1)}.
}
A short computation leads to
\cent{
$\displaystyle \left\|\nabla\left(\left|A\right|^2A\right)\right\|^2_{L^2(1)} \leq C\left\| \left|\nabla U\right|^2\nabla^2 U\right\|^2_{L^2(1)}.
$
}
Using Hölder inequalities, the inequality (\ref{biotsinfini}) of Lemma \ref{biots1bis} and the inequality (\ref{biots3nabla2}) of Lemma \ref{biots2}, we get
\aligne{
\left\|\nabla\left(\left|A\right|^2A\right)\right\|^2_{L^2(1)} &\leq C\left\|\nabla U\right\|^4_{L^\infty}\left\|\nabla^2 U\right\|^2_{L^2(1)}\\
&\leq \left\|W\right\|^2_{L^2(2)}\left\|W\right\|^2_{H^1}\left(\left\| W\right\|^2_{L^2}+\left\|\nabla W\right\|^2_{L^2(1)}\right).\\
}
Finally, taking into account the inequality (\ref{cond1bis}), we get
\eq{\label{E1I5}
\arr{ll}{I_5&\displaystyle \leq \mu \left\|\left(-\Delta \right)^{-\frac{3+\theta}{4}} f\right\|^2+\frac{C M^3 \beta^2  \gamma^3 \left(1-\theta\right)^{18} e^{-4\tau}}{\mu}\left(1+\frac{e^\tau}{\alpha_1}\right)\\
&\displaystyle \leq \mu \left\|\left(-\Delta \right)^{-\frac{3+\theta}{4}} f\right\|^2+\frac{C M^3 \gamma^3\left(1-\theta\right)^{18} e^{-3\tau}}{\mu}.}}
The equality (\ref{eqE1}) and the inequalities (\ref{E1I1}), (\ref{E1I2}), (\ref{E1I3}), (\ref{E1I4}), (\ref{E1I5}) imply that
\eq{
\arr{l}{
\displaystyle \frac{1}{2}\partial_\tau\left(\left\|\left(-\Delta\right)^{-\frac{3+\theta}{4}}f\right\|^2+\alpha_1 e^{-\tau}\left\|\left(-\Delta\right)^{-\frac{1+\theta}{4}}f\right\|^2\right)+\left(\frac{1-20\mu+\theta}{4}\right)\left\|\left(-\Delta\right)^{-\frac{3+\theta}{4}}f\right\|^2\\
\displaystyle\hspace{8cm}+\left(1+\left(\frac{1+\theta}{4}\right)\alpha_1 e^{-\tau}\right)\left\|\left(-\Delta\right)^{-\frac{1+\theta}{4}}f\right\|^2\\
\\
\hspace{0.5cm}\displaystyle\leq \frac{C M^3\gamma\left(1-\theta\right)^3 e^{-2\tau}}{\mu}+ \frac{CM\gamma \left(1-\theta\right)^3}{\mu}\left(\left\|f\right\|^2_{L^2(2)}+\left\|\nabla f\right\|^2+\alpha^2_1 e^{-2\tau}\left\|\Delta f\right\|^2_{L^2(1)}\right).
}}
Setting $\displaystyle \mu =\frac{1-\theta}{20}$, we finally get
\eq{
\arr{l}{
\displaystyle \partial_\tau E_1+\theta E_1 +\left(1+\frac{1-\theta}{4}\alpha_1 e^{-\tau}\right)\left\|\left(-\Delta\right)^{-\frac{1+\theta}{4}}f\right\|^2\leq C M^3\gamma\left(1-\theta\right)^2 e^{-2\tau}\\
\hspace{5cm}\displaystyle+ CM\gamma \left(1-\theta\right)^2\left(\left\|f\right\|^2_{L^2(2)}+\left\|\nabla f\right\|^2+\alpha^2_1 e^{-2\tau}\left\|\Delta f\right\|^2_{L^2(1)}\right).
}}
\begin{flushright}
$\square$
\end{flushright}
\subsection{Estimates in $H^1(\R^2)$}
We next establish the $H^1-$estimate of $f$. As explained earlier, we get it by performing the $L^2-$scalar product of (\ref{g3fe}) with $f$. In this section, we will see how useful the lemma \ref{lemE1} is for absorbing bad terms which appear in the computations made below. One defines the functional
\cent{$
\displaystyle E_2(\tau)=\frac{1}{2}\left(\left\|f\right\|^2+\alpha_1 e^{-\tau}\left\|\nabla f\right\|^2\right).
$}
The $H^1-$estimate of $f$ is given by the following lemma.
\lem{\label{lemE2}
Let $W\in C^1\left(\left(\tau_0,\tau_\varepsilon\right),H^1(2)\right)\cap C^0\left(\left(\tau_0,\tau_\varepsilon\right),H^3(2)\right)$ be the solution of $(\ref{g3We})$ satisfying the inequality (\ref{cond1bis}) for some $\gamma>0$. There exist $\gamma_0 >0$ and $T_0\geq 1$ such that if $T\geq T_0$ and $\gamma\leq \gamma_0$, then, for all $\tau\in \left[\tau_0, \tau^*_\varepsilon\right)$, $E_2$ satisfies the inequality
\eq{\label{E2}
\arr{l}{\partial_\tau E_2+E_2+\frac{1}{2}\left\|\nabla f\right\|^2+\frac{\beta}{2} e^{-2\tau} \left\|\left|A\right|\nabla f\right\|^2\\
\hspace{2cm}\leq \left\|f\right\|^2+CM\gamma\left(1-\theta\right)^6\left(\left\|f\right\|^2+\left\|\left|X\right|^2 f\right\|^2\right)+CM^{2}\gamma \left(1-\theta\right)^6 e^{-\tau},
}
}
where $\theta$, $0<\theta<1$ is the fixed constant introduced  at the beginning of Section \ref{secenergy}.
}
\textbf{Proof: }Taking the $L^2-$inner product of (\ref{g3fe}) with $f$, performing several integrations by parts and taking into account the equalities
\cent{$
\displaystyle \left(-\mathcal L (f) ,f\right)_{L^2}=\left\|\nabla f\right\|^2 - \frac{1}{2}\left\|f\right\|^2
$,}
and
\cent{$
\displaystyle \alpha_1 e^{-\tau}\left(\frac{X}{2}.\nabla \Delta f,f\right)_{L^2}=\alpha_1 e^{-\tau}\left\|\nabla f\right\|^2
$,}
we obtain the equality
\eq{\label{eqE2.1}
\arr{l}{\partial_\tau E_2+E_2+\varepsilon e^{-\tau} \left\|\Delta f\right\|^2+\left(1-\alpha_1 e^{-\tau}\right)\left\|\nabla f\right\|^2+\beta e^{-2\tau} \left\|\left|A\right|\nabla f\right\|^2\\
\hspace{8cm}= \left\|f\right\|^2+I_1+I_2+I_3+I_4+I_5,
}}
where
\aligne{
&I_1=-\left(K.\nabla\left(f-\alpha_1 e^{-\tau} \Delta f\right),f\right)_{L^2},\\
&I_2=-\eta\left(K.\nabla\left(G-\alpha_1 e^{-\tau} \Delta G\right),f\right)_{L^2},\\
&I_3=-\eta\left(V.\nabla\left(f-\alpha_1 e^{-\tau} \Delta f\right),f\right)_{L^2},\\
&I_4=-\eta\alpha_1 e^{-\tau} \left(\frac{\varepsilon}{\alpha_1}\Delta^2 G+\Delta G+\frac{X}{2}.\nabla \Delta G,f\right)_{L^2}+\eta\beta e^{-2\tau}\left(\dv\left(\left|A\right|^2\nabla G\right),f\right)_{L^2},\\
&I_5=\left(\beta e^{-2\tau}\dv\left(\nabla\left(\left|A\right|^2\right)\wedge A\right),f\right)_{L^2}.
}
We notice that, since $K$ is divergence free, $\left(K.\nabla f,f\right)_{L^2}=0$. Integrating by parts and using the inequality (\ref{biotsinfini}) of lemma \ref{biots1bis} and the inequality (\ref{cond2}), we obtain
\eq{\label{E2I1}
\arr{ll}{I_1 &\displaystyle =  -\alpha_1 e^{-\tau} \left(K  \Delta f,\nabla f\right)_{L^2}\\
&\displaystyle \leq C\alpha_1 e^{-\tau}\left\|K\right\|_{L^\infty}\left\| \Delta f\right\|\left\|\nabla f\right\|\\
&\displaystyle \leq C\alpha_1 e^{-\tau}\left\|f\right\|^{1/2}_{H^1}\left\|f\right\|^{1/2}_{L^2(2)}\left\| \Delta f\right\|\left\|\nabla f\right\|\\
&\displaystyle \leq C \sqrt{\alpha_1} M \gamma \left(1-\theta\right)^6 e^{-\tau/2}\left\|\nabla f\right\|\\
&\displaystyle \leq \mu \left\|\nabla f\right\|^2+\frac{C M^2 \gamma^2\left(1-\theta\right)^{12}}{\mu} e^{-\tau},
}
}
where $\mu>0$ will be made more precise later.\\
\vspace{0.5cm}\\
By the same method, using the inequality (\ref{biots4}) of the lemma \ref{biots1bis} and the smoothness of $G$, one has
\eq{\label{E2I2}
\arr{ll}{
I_2&\displaystyle =\eta\left(K\left(G-\alpha_1 e^{-\tau} \Delta G\right),\nabla f\right)_{L^2}\\
&\displaystyle \leq \left|\eta\right|\left\|K\right\|_{L^4} \left\|G-\alpha_1 e^{-\tau} \Delta G\right\|_{L^4}\left\|\nabla f\right\|\\
&\displaystyle \leq C\left(1+\alpha_1 e^{-\tau}\right)\left|\eta\right|\left\|f\right\|_{L^2(2)}\left\|\nabla f\right\|\\
&\displaystyle \leq \mu\left\|\nabla f\right\|^2+\frac{C M \gamma \left(1-\theta\right)^6}{\mu}\left(\left\|f\right\|^2+\left\|\left|X\right|^2 f\right\|^2\right).
}
}
The same method gives
\eq{\label{E2I3}
\arr{ll}{
I_3 &\displaystyle=-\alpha_1 e^{-\tau} \eta\left(V\Delta f,\nabla f\right)_{L^2}\\
&\displaystyle\leq \alpha_1 e^{-\tau} \left|\eta\right|\left\|V\right\|_{L^\infty}\left\|\Delta f\right\|\left\|\nabla f\right\|\\
&\displaystyle\leq C \sqrt{\alpha_1} M \gamma \left(1-\theta\right)^6 e^{-\tau/2}\left\|\nabla f\right\|\\
&\displaystyle\leq \mu \left\|\nabla f\right\|^2+\frac{C M^2 \gamma^2 \left(1-\theta\right)^{12}}{\mu} e^{-\tau}.
}
}
Because of the regularity of $G$, the estimate of $I_4$ is simple. Indeed, an integration by parts and Hölder inequalities yield
\aligne{
I_4 & \leq C\left|\eta\right| \left(\varepsilon+\alpha_1\right)e^{-\tau}\left\|f\right\|-\eta \beta e^{-2\tau}\left(\left|A\right|^2 \nabla G,\nabla f\right)_{L^2}\\
& \leq C\left|\eta\right| \left(\varepsilon+\alpha_1\right)e^{-\tau}\left\|f\right\|+C\left|\eta\right| \beta e^{-2\tau}\left\|\nabla G\right\|_{L^\infty}\left\|\nabla U\right\|^2_{L^3}\left\|\nabla f\right\|_{L^3}}
Then, by the inequality (\ref{biotsnabla}), the continuous injection of $H^1(\R^2)$ into $L^3(\R^2)$ and the inequalities (\ref{cond1bis}) and (\ref{cond2}), one obtains
\eq{\label{E2I4}
\arr{ll}{
I_4 &\displaystyle \leq C\left|\eta\right| \left(\varepsilon+\alpha_1\right)e^{-\tau}\left\|f\right\|+C\left|\eta\right| \beta e^{-2\tau}\left\|W\right\|^2_{L^3}\left\|\nabla f\right\|_{L^3}\\
&\displaystyle \leq C\left(\varepsilon+\alpha_1\right) M \gamma \left(1-\theta\right)^6 e^{-\tau} + C\left|\eta\right| \beta e^{-2\tau}\left\|W\right\|^2_{H^1}\left(\left\|\nabla f\right\|+\left\|\Delta f\right\|\right)\\
&\displaystyle\leq C\left(\varepsilon+\alpha_1\right) M \gamma \left(1-\theta\right)^6 e^{-\tau} +C M^{3/2} \gamma^{3/2} \left(1-\theta\right)^{9} e^{-\frac{3\tau}{2}}\\
&\displaystyle\leq C M^{3/2}  \gamma \left(1-\theta\right)^6 e^{-\tau}.
}
}
Finally, using the same arguments, due to the inequality (\ref{biotsnabla}) and the continuous injection of $H^1(\R^2)$ into $L^4(\R^2)$, one has
\eq{\label{E2I5}
\arr{ll}{
I_5 &\displaystyle = -\beta e^{-2\tau}\left(\nabla\left(\left|A\right|^2\right)\wedge A,\nabla f\right)_{L^2}\\
&\displaystyle \leq C\beta e^{-2\tau}\left\|\nabla U\right\|_{L^4}\left\|\nabla^2 U\right\|_{L^4}\left\|\left|A\right|\nabla f\right\|\\
&\displaystyle \leq C\beta e^{-2\tau}\left\|W\right\|_{H^1}\left\|W\right\|_{H^2}\left\|\left|A\right|\nabla f\right\|\\
&\displaystyle \leq C\beta M \gamma \left(1-\theta\right)^6 e^{-2\tau}\frac{e^{\tau/2}}{\sqrt{\alpha_1}}\left\|\left|A\right|\nabla f\right\|\\
&\displaystyle \leq \frac{\beta}{2} e^{-2\tau}\left\|\left|A\right|\nabla f\right\|^2+CM^2 \gamma^2 \left(1-\theta\right)^{12} e^{-\tau}.
}
}
Taking into account the inequalities (\ref{E2I1}), (\ref{E2I2}), (\ref{E2I3}), (\ref{E2I4}) and (\ref{E2I5}) and assuming that $\gamma \leq 1$, we deduce from (\ref{eqE2.1}) that
\eq{\arr{l}{
\displaystyle \partial_\tau E_2+E_2+\left(1-3\mu-\alpha_1 e^{-\tau}\right)\left\|\nabla f\right\|^2+\frac{\beta}{2} e^{-2\tau} \left\|\left|A\right|\nabla f\right\|^2\leq \\
\displaystyle \hspace{3cm}\left\|f\right\|^2+\frac{CM\gamma\left(1-\theta\right)^6}{\mu}\left(\left\|f\right\|^2+\left\|\left|X\right|^2 f\right\|^2\right)+CM^{2}\gamma \left(1-\theta\right)^6 e^{-\tau},
}
}
If we choose for instance $\mu=\frac{1}{12}$ and $T_0$ large enough to have $\alpha_1 e^{-\tau} \leq \frac{1}{4}$, we get
\eq{
\arr{l}{\partial_\tau E_2+E_2+\frac{1}{2}\left\|\nabla f\right\|^2+\frac{\beta}{2} e^{-2\tau} \left\|\left|A\right|\nabla f\right\|^2\leq \\
\hspace{3cm}\left\|f\right\|^2+CM\gamma\left(1-\theta\right)^6\left(\left\|f\right\|^2+\left\|\left|X\right|^2 f\right\|^2\right)+CM^{2}\gamma \left(1-\theta\right)^6 e^{-\tau}.
}
}
\begin{flushright}
$\square$
\end{flushright}
To achieve the $H^1-$estimate of $f$, we have to combine the inequalities (\ref{E1f}) and (\ref{E2}). We can interpolate $\left\|f\right\|^2$ between $\left\|\left(-\Delta\right)^{-\frac{1+\theta}{4}}f\right\|^2$ and $\left\|\nabla f\right\|^2$. Indeed, via Hölder and Young inequalities, we get
\aligne{
\left\|f\right\|^2 &= \left(2\pi\right)^2\int_{\R^2} \frac{1}{\left|\xi\right|^{\frac{2\left(1+\theta\right)}{3+\theta}}} \left|\xi\right|^{\frac{2\left(1+\theta\right)}{3+\theta}}\left| \widehat{f}\right|^{\frac{2\left(1+\theta\right)}{3+\theta}}\left| \widehat{f}\right|^{\frac{4}{3+\theta}} d\xi\\
&\leq \left(2\pi\right)^2\left(\int_{\R^2} \frac{1}{\left|\xi\right|^{1+\theta}} \left| \widehat{f}\right|^{2} d\xi\right)^{\frac{2}{3+\theta}}\left(\int_{\R^2} \left|\xi\right|^2 \left| \widehat{f}\right|^{2} d\xi\right)^{\frac{1+\theta}{3+\theta}}\\
&\leq \left\|\left(-\Delta\right)^{-\frac{1+\theta}{4}} f\right\|^{\frac{4}{3+\theta}}\left\|\nabla f\right\|^{\frac{2+2\theta}{3+\theta}}\\
&\leq \left(\frac{1+\theta}{3+\theta}\right)\frac{3}{8} \left\|\nabla f\right\|^2+\left(\frac{2}{3+\theta}\right)\left(\frac{8}{3}\right)^{\frac{1+\theta}{2}}\left\|\left(-\Delta\right)^{-\frac{1+\theta}{4}} f\right\|^2.
}
Since, $0<\theta<1$, we obtain
\eq{\label{interpolation}
\left\|f\right\|^2 \leq \frac{1}{4}\left\|\nabla f\right\|^2+5\left\|\left(-\Delta\right)^{-\frac{1+\theta}{4}} f\right\|^2.
}
Thus, we have
\eq{\label{E2bis}
\arr{l}{\displaystyle \partial_\tau E_2+E_2+\frac{1}{4}\left\|\nabla f\right\|^2+\frac{\beta}{2} e^{-2\tau} \left\|\left|A\right|\nabla f\right\|^2\leq \\
\displaystyle \hspace{1cm }5\left\|\left(-\Delta\right)^{-\frac{1+\theta}{4}} f\right\|^2+CM\gamma\left(1-\theta\right)^6\left(\left\|f\right\|^2+\left\|\left|X\right|^2 f\right\|^2\right)+CM^{2}\gamma \left(1-\theta\right)^6  e^{-\tau}.
}
}
We define $E_3=6E_1+E_2$. The inequalities $(\ref{E1f})$ and $(\ref{E2bis})$ give
\eq{
\arr{l}{
\displaystyle\partial_\tau E_3+\theta E_3+\left(1+ \frac{3}{2}\left(1-\theta\right)\alpha_1 e^{-\tau}\right)\left\|\left(-\Delta\right)^{-\frac{1+\theta}{4}}f\right\|^2+\frac{1}{4}\left\|\nabla f\right\|^2 \leq CM^{3}\gamma \left(1-\theta\right)^2 e^{-\tau} \\
\displaystyle\hspace{0.3cm}+CM\gamma\left(1-\theta\right)^2\left(\left\|f\right\|^2+\left\|\nabla f\right\|^2+\alpha^2_1 e^{-2\tau} \left\|\Delta f\right\|^2+\left\|\left|X\right|^2 f\right\|^2+\alpha^2_1 e^{-2\tau} \left\|\left|X\right|^2 \Delta f\right\|^2\right)
.}
}
Interpolating again $\left\|f\right\|^2$ between $\left\|\nabla f\right\|^2$ and $\left\|\left(-\Delta\right)^{-\frac{1+\theta}{4}} f\right\|^2$ and taking $\gamma$ sufficiently small, we obtain
\eq{\label{E3}
\arr{l}{
\displaystyle\partial_\tau E_3+\theta E_3+\frac{1}{2}\left\|\left(-\Delta\right)^{-\frac{1+\theta}{4}}f\right\|^2+\frac{1}{8}\left\|\nabla f\right\|^2\leq CM^{3}\gamma \left(1-\theta\right)^2 e^{-\tau}\\
\\
\hspace{2cm}\displaystyle  +CM\gamma\left(1-\theta\right)^2\left(\alpha^2_1 e^{-2\tau} \left\|\Delta f\right\|^2+\left\|\left|X\right|^2 f\right\|^2+\alpha^2_1 e^{-2\tau} \left\|\left|X\right|^2 \Delta f\right\|^2\right)
.}
}
\subsection{Estimates in $H^2(\R^2)$}
We now perform the $H^2-$estimate of $f$. This is done with the same method as for the $H^1-$estimate in the previous section. Indeed, we perform the $L^2-$product between (\ref{g3fe}) and $-\Delta f$ and, after some computations, we see that the inequality (\ref{cond1bis}) enables us to absorb all the terms involving the $H^2-$norm of $f$. Combined with (\ref{E3}), we get an estimate in $H^2$, where only terms with weighted norms remain. More precisely, we introduce the following functional.
\cent{$\displaystyle E_4(\tau)=\frac{1}{2}\left(\left\|\nabla f\right\|^2+\alpha_1 e^{-\tau}\left\|\Delta f\right\|^2\right)$.}
The $H^2-$estimate of $f$ is given by the lemma below.
\lem{\label{lemE4}
Let $W\in C^1\left(\left(\tau_0,\tau_\varepsilon\right),H^1(2)\right)\cap C^0\left(\left(\tau_0,\tau_\varepsilon\right),H^3(2)\right)$ be the solution of $(\ref{g3We})$ satisfying the inequality (\ref{cond1bis}) for some $\gamma>0$. There exist $\gamma_0 >0$ and $T_0 \geq 1$ such that if $T\geq T_0$ and $\gamma\leq \gamma_0$, then for all $ \tau\in\left[\tau_0,  \tau^*_\varepsilon\right)$, $E_4$ satisfies the inequality
\eq{\label{E4}
\arr{l}{\displaystyle\partial_\tau  E_4 + E_4+\frac{1}{2}\left\|\Delta f\right\|^2+\frac{\beta}{2}e^{-2\tau}\left\|\left|A\right|\Delta f\right\|^2 \leq \frac{3}{2} \left\|\nabla f\right\|^2+ C  M^2\gamma \left(1-\theta\right)^6 e^{-2\tau}\\
\displaystyle \hspace{6cm}+ C M^2\gamma \left(1-\theta\right)^6\left(\left\|f\right\|^2+\left\|\nabla f\right\|^2+\left\|\left|X\right|^2f\right\|^2\right),
}
}
where $\theta$, $0<\theta<1$ is the fixed constant introduced  at the beginning of Section \ref{secenergy}.
}
\textbf{Proof: }We take the $L^2-$product of $(\ref{g3fe})$ with $-\Delta f$. Doing several integrations by parts, it is easy to see that
\cent{$
\displaystyle \left(-\mathcal L(f),-\Delta f\right)_{L^2}=\left\|\Delta f\right\|^2-\left\|\nabla f\right\|^2
$,}
and
\cent{$\displaystyle -\left(\alpha_1 e^{-\tau}\frac{X}{2}.\nabla \Delta f, \Delta f\right)_{L^2}= \frac{1}{2}\alpha_1 e^{-\tau}\left\|\Delta f\right\|^2$.}
Furthermore, one also has
\cent{$
\displaystyle \beta e^{-2\tau}\left(\dv\left(\left|A\right|^2 \nabla f \right),\Delta f\right)=\beta e^{-2\tau} \left\|\left|A\right|\Delta f\right\|^2+\beta e^{-2\tau} \sum^{2}_{j=1} \int_{\R^2} A:\partial_j A \partial_j f\Delta f dX
$.}
Using Hölder inequalities, the inequality (\ref{biotsnabla}) of lemma \ref{biots1}, the continuous injections of $H^1(\R^2)$ into $L^4(\R^2)$ and the inequality (\ref{cond1bis}), we get
\aligne{
\beta e^{-2\tau} \sum^{2}_{j=1} \int_{\R^2} A:\partial_j A \partial_j f\Delta f dX&\leq C\beta e^{-2\tau}\left\|\left|A\right|\Delta f\right\|\left\|\nabla A \nabla f \right\|\\
&\leq C\beta e^{-2\tau}\left\|\left|A\right|\Delta f\right\|\left\|\nabla^2 U\right\|_{L^4}\left\| \nabla f \right\|_{L^4}\\
&\leq C\beta e^{-2\tau}\left\|\left|A\right|\Delta f\right\|\left\| \nabla W\right\|_{H^1}\left\| \nabla f \right\|_{H^1}\\
&\leq \mu_1 \beta e^{-2\tau}\left\|\left|A\right|\Delta f\right\|^2\\
&\hspace{2cm}+ \frac{C\beta}{\mu_1} e^{-2\tau}\left\|W\right\|^2_{H^2}\left(\left\|\nabla f\right\|^2+ \left\|\Delta f\right\|^2\right)\\
&\leq \mu_1 \beta e^{-2\tau}\left\|\left|A\right|\Delta f\right\|^2\\
&\hspace{2cm}+\frac{C M \gamma \left(1-\theta\right)^6}{\mu_1} e^{-\tau}\left(\left\|\nabla f\right\|^2+\left\|\Delta f\right\|^2\right),
}
where $\mu_1>0$ will be chosen later.\\
\vspace{0.5cm}\\
Consequently, we get
\eq{\label{E4i}
\arr{l}{
\displaystyle \partial_\tau E_4 + \varepsilon e^{-\tau}\left\|\nabla \Delta f\right\|^2+\left(1-\frac{\alpha_1}{2} e^{-\tau}\right)\left\|\Delta f\right\|^2+\beta\left(1-\mu_1\right) e^{-2\tau}\left\|\left|A\right|\Delta f\right\|^2 \\
\hspace{2cm}\leq  \displaystyle \left\|\nabla f\right\|^2+\frac{C M \gamma \left(1-\theta\right)^6}{\mu_1} e^{-\tau}\left(\left\|\nabla f\right\|^2+\left\|\Delta f\right\|^2\right) +I_1+I_2+I_3+I_4+I_5,
}
}
where
\aligne{
&I_1=\left(U.\nabla\left(f-\alpha_1 e^{-t} \Delta f\right),\Delta f\right)_{L^2},\\
&I_2=\eta\left(K.\nabla\left(G-\alpha_1 e^{-t} \Delta G\right),\Delta f\right)_{L^2},\\
&I_3=\eta\alpha_1 e^{-t} \left(\frac{\varepsilon}{\alpha_1} \Delta^2 G+ \Delta G+\frac{X}{2}.\nabla \Delta G,\Delta f\right)_{L^2}+\eta\beta e^{-2t}\left(\dv\left(\left|A\right|^2\nabla G\right),\Delta f\right)_{L^2},\\
&I_4=\beta e^{-2t}\left(\dv\left(\nabla\left(\left|A\right|^2\right)\wedge A\right),\Delta f\right)_{L^2}.
}
Integrating by parts and using the divergence free property of $K$, one can show that
\cent{
$I_1=-\displaystyle \sum^2_{j,k=1} \int_{\R^2} \partial_k U_j \partial_j f\partial_k f dX
$.
}
Due to the Gagliardo-Niremberg inequality and the inequalities (\ref{biotsnabla}) and (\ref{cond1bis}), it comes
\eq{\label{E4I1}
\arr{ll}{I_1&\displaystyle \leq C\left\|\nabla U\right\|\left\|\nabla f\right\|^2_{L^4}\\
&\displaystyle \leq C\left\|\nabla U\right\|\left\|\nabla f\right\|\left\|\Delta f\right\|\\
&\displaystyle \leq \mu_2\left\|\Delta f\right\|^2+ \frac{CM\gamma \left(1-\theta\right)^6}{\mu_2}\left\|\nabla f\right\|^2,
}
}
where $\mu_2 >0$ will be chosen later.\\
\vspace{0.5cm}\\
We now estimate $I_2$ with the help of the inequality (\ref{biotsinfini}) of Lemma \ref{biots1bis}, the inequality (\ref{cond2}) and the smoothness of $G$.
\eq{\label{E4I2}
\arr{ll}{I_2 &\displaystyle \leq  \left|\eta\right|\left\|K\right\|_{L^\infty}\left\|G-\alpha_1 e^{-\tau} \Delta G\right\|\left\|\Delta f\right\|\\
&\displaystyle \leq C\left|\eta\right|\left(1+\alpha_1 e^{-\tau}\right)\left\|f\right\|^{1/2}_{H^1}\left\|f\right\|^{1/2}_{L^2(2)}\left\|\Delta f\right\|\\
&\displaystyle \leq \mu_2\left\|\Delta f\right\|^2+\frac{CM\gamma\left(1-\theta\right)^6}{\mu_2}\left(\left\|f\right\|^2+\left\|\nabla f\right\|^2+\left\|\left|X\right|^2 f\right\|^2\right).\\
}
}
We rewrite
\cent{$
I_3=I^1_3+I^2_3,
$}
where
\aligne{&I^1_3 = \eta\alpha_1 e^{-\tau} \left(\frac{\varepsilon}{\alpha_1}\Delta^2 G+\Delta G+\frac{X}{2}.\nabla \Delta G,\Delta f\right)_{L^2},\\
&I^2_3=\eta\beta e^{-2\tau}\left(\dv\left(\left|A\right|^2\nabla G\right),\Delta f\right)_{L^2}.}
Using the good regularity of $G$ and the inequality (\ref{cond2}), one can show that
\aligne{
I^1_3 &\leq C M^{1/2}\gamma^{1/2}\left(1-\theta\right)^{3} e^{-\tau}\left\|\Delta f\right\|\\
&\leq \mu_2 \left\|\Delta f\right\|^2+\frac{C M\gamma\left(1-\theta\right)^6}{\mu_2} e^{-2\tau}.
}
The estimate of $I^2_3$ is slightly more complicated. Actually, we can bound $I^2_3$ by two kinds of terms that we estimate separately. In fact, it is easy to see that
\eq{\label{E4bis}
\displaystyle I^2_3 \leq C\left|\eta\right|\beta e^{-2\tau}\int_{\R^2}\left|\nabla A\right|\left|A\right|\left|\nabla G\right|\left|\Delta f\right|dX+C\left|\eta\right|\beta e^{-2\tau}\int_{\R^2}\left|A\right|^2\left|\nabla^2 G\right|\left|\Delta f\right|dX.
}
Each term of the right hand side of (\ref{E4bis}) can be estimated in a convenient way. We use again the inequality (\ref{biotsnabla}) of the lemma \ref{biots1}, the inequality (\ref{cond1bis}) , the Hölder inequalities and the inequality (\ref{cond1bis}). We get
\aligne{
C\left|\eta\right|\beta e^{-2\tau}\int_{\R^2}\left|\nabla A\right|\left|A\right|\left|\nabla G\right|\left|\Delta f\right|dX & \leq C\left|\eta\right|\beta e^{-2\tau} \left\| \nabla^2 U\right\|\left\|\nabla G\right\|_{L^\infty}\left\|\left|A\right|\Delta f\right\|\\
&\leq \mu_1 \beta e^{-2\tau}\left\|\left|A\right|\Delta f\right\|^2 + \frac{C\left|\eta\right|^2\beta}{\mu_1} e^{-2\tau}\left\| \nabla W\right\|^2\\
&\leq \mu_1\beta e^{-2\tau}\left\|\left|A\right|\Delta f\right\|^2 + \frac{CM^2\gamma^2\left(1-\theta\right)^{12}}{\mu_1} e^{-2\tau}.
}
By the same way, we have
\cent{
$\displaystyle C\left|\eta\right|\beta e^{-2\tau}\int_{\R^2}\left|A\right|^2\left|\nabla^2 G\right|\left|\Delta f\right|dx \leq \mu_1\beta e^{-2\tau}\left\|\left|A\right|\Delta f\right\|^2 + \frac{CM^2\gamma^2\left(1-\theta\right)^{12}}{\mu_1} e^{-2\tau}$,
}
and thus we have shown
\cent{
$
\displaystyle I^2_3 \leq 2\mu_1\beta e^{-2\tau}\left\|\left|A\right|\Delta f\right\|^2+ \frac{CM^2\gamma^2\left(1-\theta\right)^{12}}{\mu_1} e^{-2\tau}.
$
}
Finally, assuming $\gamma\leq 1$, one has
\eq{\label{E4I3}
I_3 \leq \mu_2 \left\|\Delta f\right\|^2+2\mu_1\beta e^{-2\tau}\left\|\left|A\right|\Delta f\right\|^2+ \frac{CM^2\gamma\left(1-\theta\right)^{6}}{\min(\mu_1,\mu_2)} e^{-2\tau}.
}
It remains to estimate $I_4$. Recalling that $U=\eta V + K$, one has
\eq{\label{E4ter}
\arr{ll}{I_4 &\leq C\left|\eta\right|\beta e^{-2\tau}\int_{\R^2}\left|\nabla A\right|\left|A\right|\left|\nabla^2 V\right|\left|\Delta f\right|dX+C\left|\eta\right|\beta e^{-2\tau}\int_{\R^2}\left|A\right|^2\left|\nabla^3 V\right|\left|\Delta f\right|dX\\
\\
&\hspace{1.5cm}+C \beta e^{-2\tau}\int_{\R^2} \left|\nabla A\right| \left|A\right|\left|\nabla^2K\right|\left|\Delta f\right| dX+C \beta e^{-2\tau}\int_{\R^2} \left|A\right|^2 \left|\nabla^3 K\right|\left|\Delta f\right| dX .
}
}
We have to estimate each term of the right hand side of the inequality (\ref{E4ter}). The first two ones can be estimated exactly like we did for $I^2_3$. The inequality (\ref{biotsnabla}) of lemma \ref{biots1} and Gagliardo-Niremberg inequality yield
\aligne{
C\beta e^{-2\tau}\int_{\R^2} \left|\nabla A\right| \left|A\right|\left|\nabla^2 K\right|\left|\Delta f\right| dX&\leq \mu_1 \beta e^{-2\tau}\left\|\left|A\right| \Delta f\right\|^2\\
& \hspace{1cm}+ \frac{C \beta e^{-2\tau}}{\mu_1}\left\|\nabla^2 U\right\|^2_{L^4} \left\|\nabla^2 K\right\|^2_{L^4}\\
&\leq \mu_1 \beta e^{-2\tau}\left\|\left|A\right| \Delta f\right\|^2\\
& \hspace{1cm}+ \frac{C \beta e^{-2\tau}}{\mu_1}\left\|\nabla W\right\|^2_{L^4} \left\|\nabla f\right\|^2_{L^4}\\
&\leq \mu_1 \beta e^{-2\tau}\left\|\left|A\right| \Delta f\right\|^2\\
& \hspace{1cm}+ \frac{C \beta e^{-2\tau}}{\mu_1}\left\|\nabla W\right\|\left\|\Delta W\right\| \left\|\nabla f\right\|\left\|\Delta f\right\|.
}
Due to the inequality (\ref{cond1bis}), we get
\aligne{
C\beta e^{-2\tau}\int_{\R^2} \left|\nabla A\right| \left|A\right|\left|\nabla^2 K\right|\left|\Delta f\right| dX&\leq \mu_1 \beta e^{-2\tau}\left\|\left|A\right| \Delta f\right\|^2\\
&\hspace{1cm}+ \frac{C M\gamma \left(1-\theta\right)^6 e^{-\frac{3\tau}{2}}}{\mu_1}\left(\left\|\nabla f\right\|^2+\left\|\Delta f\right\|^2\right).
}
By the same method, we obtain
\aligne{
C\beta e^{-2\tau}\int_{\R^2} \left|A\right|^2 \left|\nabla^3 K\right|\left|\Delta f\right| dX &\leq \mu_1 \beta e^{-2\tau}\left\|\left|A\right| \Delta f\right\|^2 + \frac{C \beta e^{-2\tau}}{\mu_1}\left\|\nabla U\right\|^2_{L^\infty}\left\|\nabla^3 K\right\|^2\\
&\leq \mu_1 \beta e^{-2\tau}\left\|\left|A\right| \Delta f\right\|^2\\
&\hspace{1cm} + \frac{C \beta e^{-2\tau}}{\mu_1}\left\|\nabla W\right\|_{H^1}\left\|\nabla W\right\|_{L^2(2)}\left\|\Delta f\right\|^2\\
&\leq \mu_1 \beta e^{-2\tau}\left\|\left|A\right| \Delta f\right\|^2 + \frac{C M\gamma \left(1-\theta\right)^6 e^{-\tau}}{\mu_1}\left\|\Delta f\right\|^2.
}
Finally, we have shown that
\eq{\label{E4I4}
\displaystyle I_4 \leq 4\mu_1\beta e^{-2\tau}\left\|\left|A\right|\Delta f\right\|^2+ \frac{C M^2\gamma \left(1-\theta\right)^6 e^{-\tau}}{\mu_1}\left(\left\|\nabla f\right\|^2+\left\|\Delta f\right\|^2\right).
}
Going back to (\ref{E4i}) and taking into account the inequalities (\ref{E4I1}), (\ref{E4I2}), (\ref{E4I3}) and (\ref{E4I4}), we get
\cent{$\arr{l}{
\displaystyle \partial_\tau E_4 + \left(1-3\mu_2-\frac{\alpha_1}{2} e^{-\tau}\right)\left\|\Delta f\right\|^2+\left(1-7\mu_1\right)\beta e^{-2\tau}\left\|\left|A\right|\Delta f\right\|^2 \leq \left\|\nabla f\right\|^2\\
\displaystyle \hspace{1cm} + \frac{C M^2\gamma \left(1-\theta\right)^6}{\min(\mu_1,\mu_2)}\left(\left\|f\right\|^2+\left\|\nabla f\right\|^2+\left\|\Delta f\right\|^2+\left\|\left|X\right|^2 f\right\|^2\right)+ \frac{C M^2\gamma \left(1-\theta\right)^6}{\min(\mu_1,\mu_2)}e^{-2\tau}.
}$
}
Taking  for instance $\mu_1= \frac{1}{14}$, $\mu_2=\frac{1}{12}$, $\gamma$ small enough and $T = e^{\tau_0}$ large enough, we finally have
\eq{\label{E4f}
\arr{l}{\displaystyle\partial_\tau  E_4 + E_4+\frac{1}{2}\left\|\Delta f\right\|^2+\frac{\beta}{2}e^{-2\tau}\left\|\left|A\right|\Delta f\right\|^2 \leq \frac{3}{2} \left\|\nabla f\right\|^2+ C  M^2\gamma \left(1-\theta\right)^6 e^{-2\tau}\\
\displaystyle \hspace{6cm}+ C M^2\gamma \left(1-\theta\right)^6\left(\left\|f\right\|^2+\left\|\nabla f\right\|^2+\left\|\left|X\right|^2f\right\|^2\right).
}
}
\begin{flushright}
$\square$
\end{flushright}
In order to finish the $H^2-$estimate of $f$ we define a new functional $E_5$ as a linear combination of $E_3$ and $E_4$ given by
\cent{$
E_5=16E_3+E_4.
$}
From the inequalities $(\ref{E3})$ and $(\ref{E4})$, it is clear that one has
\eq{
\arr{l}{
\displaystyle \partial_\tau E_5+\theta E_5+8\left\|\left(-\Delta\right)^{-\frac{1+\theta}{4}} f\right\|^2+\frac{1}{2}\left\|\nabla f\right\|^2+\frac{1}{2}\left\|\Delta f\right\|^2  \leq CM^3\gamma \left(1-\theta\right)^2 e^{-\tau}\\
\\
\displaystyle\hspace{4cm}+CM^2\gamma \left(1-\theta\right)^2\Big(\left\|f\right\|^2+\left\|\nabla f\right\|^2+\alpha^2_1 e^{-2\tau}\left\|\Delta f\right\|^2\\
\hspace{9cm}+\left\|\left|X\right|^2 f\right\|^2+\alpha^2_1 e^{-2\tau}\left\|\left|X\right|^2\Delta f\right\|^2\Big).
}
}
Using the interpolation inequality (\ref{interpolation}) and taking $\gamma$ small enough and $\tau_0=\log(T)$ large enough, we finally obtain
\eq{\label{E5}
\arr{l}{\displaystyle \partial_\tau E_5+\theta E_5+7\left\|\left(-\Delta\right)^{-\frac{1+\theta}{4}}f\right\|^2+\frac{1}{4}\left\|\nabla f\right\|^2+\frac{1}{4}\left\|\Delta f\right\|^2\leq CM^3\gamma \left(1-\theta\right)^2 e^{-\tau}\\
\displaystyle \hspace{6cm} +CM^2\gamma \left(1-\theta\right)^2\left(\left\|\left|X\right|^2 f\right\|^2+\alpha^2_1 e^{-2\tau}\left\|\left|X\right|^2\Delta f\right\|^2\right).\\
}
}
\subsection{Estimates in  $H^2(2)$}
In order to achieve the estimate of $f$ in $H^2(2)$, it remains to perform estimates in weighted spaces. Combined with the inequality (\ref{E5}), it will give us an estimate in $H^2(2)$. To do this, we make the $L^2-$scalar product of (\ref{g3fe}) with $\displaystyle\left|X\right|^4 \left(f-\alpha_1e^{-\tau}\Delta f\right)$. We define the functional
\cent{$
\displaystyle E_6(\tau) =\frac{1}{2}\left\|\left|X\right|^2\left(f-\alpha_1 e^{-\tau}\Delta f\right)\right\|^2.
$}
Before stating the lemma which contains the estimate of $E_6$, we state a technical lemma, which gives the terms provided by the $L^2-$product of the linear terms of (\ref{g3fe}) with $\displaystyle\left|X\right|^4 \left(f-\alpha_1e^{-\tau}\Delta f\right)$.
\lem{\label{lemcalcweight}
Let $f \in C^1\left(\left(\tau_0,\tau_\varepsilon\right),H^1(2)\right)\cap C^0\left(\left(\tau_0,\tau_\varepsilon\right),H^3(2)\right)$ and $H(X , \tau, f)=\left|X\right|^4\left(f-\alpha_1 e^{-\tau}\Delta f\right)$. For all $\tau \in \left(\tau_0,\tau_\varepsilon\right)$, the next equalities hold.
\begin{enumerate}
\item $\left(-f ,H(X , \tau, f)\right)_{L^2}=-\left\|\left|X\right|^2 f\right\|^2 +8\alpha_1 e^{-\tau}\left\|\left|X\right| f \right\|^2-\alpha_1 e^{-\tau}\left\|\left|X\right|^2 \nabla f\right\|^2$.
\item $\left(-\Delta f ,H(X , \tau, f)\right)_{L^2}=\alpha_1 e^{-\tau} \left\|\left|X\right|^2\Delta f\right\|^2-8\left\|\left|X\right| f \right\|^2+\left\|\left|X\right|^2 \nabla f\right\|^2.$
\item $\arr{l}{
\left(-\frac{X}{2}.\nabla f,H(X , \tau, f)\right)_{L^2} = \frac{3}{2}\left\|\left|X\right|^2 f\right\|^2-24 \alpha_1 e^{-\tau} \left\|\left|X\right| f\right\|^2 \\
\hspace{5cm}+ 3 \alpha_1 e^{-\tau}\left\|\left|X\right|^2\nabla f\right\|^2- \frac{\alpha_1}{2} e^{-\tau} \left( X.\nabla \Delta f , \left|X\right|^4 f\right)_{L^2}.
}$
\item $\arr{l}{
\left(-\mathcal L(f),H(X , \tau, f)\right)_{L^2} = 
\frac{1}{2}\left\|\left|X\right|^2 f\right\|^2+\left(1+2\alpha_1e^{-\tau}\right)\left\|\left|X\right|^2 \nabla f\right\|^2\\
\hspace{1cm}+\alpha_1 e^{-\tau}\left\|\left|X\right|^2 \Delta f\right\|^2-\left(8+16\alpha_1 e^{-\tau}\right)\left\|\left|X\right| f\right\|^2- \frac{\alpha_1}{2} e^{-\tau} \left( X.\nabla \Delta f , \left|X\right|^4 f\right)_{L^2}.
}$
\item $\left(\alpha_1 e^{-\tau}\frac{X}{2}.\nabla \Delta f,H(X , \tau, f)\right)_{L^2}=\frac{\alpha_1}{2} e^{-\tau}\left(X.\nabla \Delta f,\left|X\right|^4 f\right)_{L^2}+\frac{3\alpha^2_1}{2} e^{-2\tau} \left\| \left|X\right|^2\Delta f\right\|^2.$
\item $\arr{l}{\varepsilon e^{-\tau} \left(\Delta^2 f,H(X , \tau, f)\right)_{L^2} =\varepsilon \alpha_1 e^{-2\tau} \left(\left\| \left|X\right|^2 \nabla \Delta f \right\|^2- 8  \left\|\left|X\right|\Delta f\right\|^2\right)\\
\hspace{2cm}+\varepsilon e^{-\tau}\left(\left\|\left|X\right|^2 \Delta f\right\|^2 -8 \left\|\left|X\right| \nabla f\right\|^2+ 32 \left\| f\right\|^2 -16 \left\|X. \nabla f\right\|^2\right).}$
\end{enumerate}
}
\textbf{Proof: }All these equalities are obtained via integrations by parts. We only show the first four ones, the others are obtained with the same method. Let us show the equality 1. Two integrations by parts imply
\aligne{
\left(-f ,\left|X\right|^4\left(f-\alpha_1 e^{-\tau}\Delta f\right)\right)_{L^2} & =-\left\|\left|X\right|^2 f\right\|^2 - \alpha_1 e^{-\tau} \left\|\left|X\right|^2 \nabla f\right\|^2\\
&\hspace{3cm} - 4\alpha_1 e^{-\tau} \sum^2_{j=1}\int_{\R^2} X_j \left|X\right|^2 f \partial_j f dX\\
& =-\left\|\left|X\right|^2 f\right\|^2 - \alpha_1 e^{-\tau} \left\|\left|X\right|^2 \nabla f\right\|^2\\
&\hspace{3cm} - 2\alpha_1 e^{-\tau} \sum^2_{j=1}\int_{\R^2} X_j \left|X\right|^2 \partial_j \left(f^2\right) dX\\
& =-\left\|\left|X\right|^2 f\right\|^2 - \alpha_1 e^{-\tau} \left\|\left|X\right|^2 \nabla f\right\|^2 + 8\alpha_1 e^{-\tau} \left\|\left|X\right| f\right\|^2.
}
The equality 2. is obtained through the same computations. We show now the third equality of this lemma. Integrating by parts, we obtain
\aligne{
\left(-\frac{X}{2}.\nabla f,\left|X\right|^4\left(f-\alpha_1 e^{-\tau} \Delta f\right)\right)_{L^2}&=-\sum^{2}_{j=1}\int_{\R^2} \frac{X_j\left|X\right|^4}{4}\partial_j\left(\left|f\right|^2\right)dX\\
&\hspace{2cm}+\alpha_1 e^{-\tau} \sum^{2}_{j=1}\int_{\R^2} \frac{X_j\left|X\right|^4}{2} \partial_j f \Delta f dX\\
&=\frac{3}{2}\left\|\left|X\right|^2 f\right\|^2+\alpha_1 e^{-\tau} \sum^{2}_{j=1}\int_{\R^2} \frac{X_j\left|X\right|^4}{2} \partial_j f \Delta f dX.\\
}
Besides, integrating several times by parts, we get
\aligne{
\alpha_1 e^{-\tau} \sum^{2}_{j=1}\int_{\R^2} \frac{X_j\left|X\right|^4}{2} \partial_j f \Delta f dX&= - \alpha_1 e^{-\tau} \sum^{2}_{j=1}\int_{\R^2}f \partial_j\left( \frac{X_j\left|X\right|^4}{2} \Delta f\right) dX\\
&= - 3\alpha_1 e^{-\tau} \int_{\R^2}  \left|X\right|^4 f \Delta f dX\\
&\hspace{3cm}- \frac{\alpha_1}{2} e^{-\tau} \left( X.\nabla \Delta f , \left|X\right|^4 f\right)_{L^2}\\
&=-24 \alpha_1 e^{-\tau} \left\|\left|X\right| f\right\|^2 + 3 \alpha_1 e^{-\tau}\left\|\left|X\right|^2\nabla f\right\|^2\\
&\hspace{4cm}- \frac{\alpha_1}{2} e^{-\tau} \left( X.\nabla \Delta f , \left|X\right|^4 f\right)_{L^2},
}
and consequently
\cent{
$\arr{l}{
\left(-\frac{X}{2}.\nabla f,\left|X\right|^4\left(f-\alpha_1 e^{-\tau} \Delta f\right)\right)_{L^2} = \frac{3}{2}\left\|\left|X\right|^2 f\right\|^2-24 \alpha_1 e^{-\tau} \left\|\left|X\right| f\right\|^2 \\
\hspace{5cm}+ 3 \alpha_1 e^{-\tau}\left\|\left|X\right|^2\nabla f\right\|^2- \frac{\alpha_1}{2} e^{-\tau} \left( X.\nabla \Delta f , \left|X\right|^4 f\right)_{L^2}.
}$
}
The fourth equality of this lemma is obtained by summing the first three ones. By the same method, we obtain easily the equalities 5. and 6. of this lemma.
\begin{flushright}
$\square$
\end{flushright}
The $H^2(2)-$estimate of $f$ is given in the following lemma.
\lem{\label{lemE6}
Let $W\in C^1\left(\left(\tau_0,\tau_\varepsilon\right),H^1(2)\right)\cap C^0\left(\left(\tau_0,\tau_\varepsilon\right),H^3(2)\right)$ be the solution of $(\ref{g3We})$ satisfying the inequality (\ref{cond1bis}) for some $\gamma>0$. There exist $\gamma_0 >0$ and $T_0\geq 1$ such that if $T\geq T_0$ and $\gamma\leq \gamma_0$, then for all $\tau \in \left[\tau_0, \tau^*_\varepsilon\right)$, $E_6$ satisfies the inequality
\eq{\label{E6}
\arr{l}{
\displaystyle \partial_\tau E_6+\theta E_6+\frac{1-\theta}{8}\left\|\left|X\right|^2 f\right\|^2+\frac{1}{4}\left\|\left|X\right|^2 \nabla f\right\|^2+\frac{\alpha_1}{4} e^{-\tau}\left\|\left|X\right|^2 \Delta f\right\|^2\displaystyle \\
\\
\displaystyle \hspace{0.5cm}\leq C M^2\gamma \left(1-\theta\right)^{6} e^{-\tau}+\frac{1024}{1-\theta}\left\|f\right\|^2+C M^2\gamma^{1/2} \left(1-\theta\right)^{3}\left(\left\|f\right\|^2+\left\|\nabla f\right\|^2+\left\|\Delta f\right\|^2\right),
}
}
where $\theta$, $0<\theta<1$ is the fixed constant introduced  at the beginning of Section \ref{secenergy}.
}
\textbf{Proof: }To show this lemma, we perform the $L^2-$product of the equality $(\ref{g3fe})$ with $\left|X\right|^4\left(f-\alpha_1 e^{-\tau}\Delta f\right)$. Applying the lemma \ref{lemcalcweight}, we obtain
\eq{\label{eqE6.0}
\arr{l}{\displaystyle\partial_\tau E_6+\frac{1}{2}\left\|\left|X\right|^2 f\right\|^2+\left(1+\alpha_1e^{-\tau}\right)\left\|\left|X\right|^2 \nabla f\right\|^2+\left(\alpha_1 e^{-\tau}+\frac{\alpha^2_1}{2}e^{-2\tau}\right)\left\|\left|X\right|^2 \Delta f\right\|^2+J\\
\\
\hspace{1cm} =\displaystyle C\varepsilon e^{-\tau} \left\|\left|X\right|\nabla f\right\|^2+C\varepsilon \alpha_1 e^{-2\tau} \left\|\left|X\right|\Delta f\right\|^2+\left(8+8\alpha_1 e^{-\tau}\right)\left\|\left|X\right| f\right\|^2\\
\hspace{10cm}+I_1+I_2+I_3+I_4+I_5,}
}
where
\aligne{
&J=-\beta e^{-2\tau}\left(\dv\left(\left|A\right|^2 \nabla f\right),\left|X\right|^4\left(f-\alpha_1 e^{-\tau}\Delta f\right)\right)_{L^2},\\
&I_1=\left(K.\nabla\left(f-\alpha_1 e^{-\tau} \Delta f\right),\left|X\right|^4\left(f-\alpha_1 e^{-\tau}\Delta f\right)\right)_{L^2},\\
&I_2=\eta\left(K.\nabla\left(G-\alpha_1 e^{-\tau} \Delta G\right),\left|X\right|^4\left(f-\alpha_1 e^{-\tau}\Delta f\right)\right)_{L^2},\\
&I_3=\eta\left(V.\nabla\left(f-\alpha_1 e^{-\tau} \Delta f\right),\left|X\right|^4\left(f-\alpha_1 e^{-\tau}\Delta f\right)\right)_{L^2},\\
&I_4=-\eta\varepsilon e^{-\tau}\left( \Delta^2 G,\left|X\right|^2\left(f-\alpha_1 e^{-\tau} \Delta f\right)\right)_{L^2}\\
&\hspace{3cm}+\eta\alpha_1 e^{-\tau} \left(\Delta G+\frac{X}{2}.\nabla \Delta G,\left|X\right|^4\left(f-\alpha_1 e^{-\tau}\Delta f\right)\right)_{L^2}\\
&\hspace{5cm}-\eta\beta e^{-2\tau}\left(\dv\left(\left|A\right|^2\nabla G\right),\left|X\right|^4\left(f-\alpha_1 e^{-\tau}\Delta f\right)\right)_{L^2},\\
&I_5=-\beta e^{-2\tau}\left(\dv\left(\nabla\left(\left|A\right|^2\right)\wedge A\right),\left|X\right|^4\left(f-\alpha_1 e^{-\tau}\Delta f\right)\right)_{L^2}.
}
We now estimate $J$. One has
\eq{\label{eqE6.3}
J=J_1+J_2,
}
where
\cent{
$
\arr{l}{
J_1 = -\beta e^{-2\tau} \left(\dv\left(\left|A\right|^2 \nabla f\right),\left|X\right|^4 f\right)_{L^2},\\
J_2 = \beta \alpha_1 e^{-3\tau} \left(\dv\left(\left|A\right|^2 \nabla f\right),\left|X\right|^4\Delta f)\right)_{L^2}.
}
$
}
We estimate $J_1$ and $J_2$ separately. Integrating by parts, we obtain
\cent{$
\displaystyle J_1 = \beta e^{-2\tau}\left\|\left|X\right|^2\left|A\right| \nabla f\right\|^2 + 4  \beta e^{-2\tau} \sum^{2}_{j=1} \int_{\R^2} X_j \left|X\right|^2 \left|A\right|^2 \partial_j f f dX.
$}
Using Hölder and Young inequalities, we obtain 
\aligne{
\left|4\beta e^{-2\tau} \sum^{2}_{j=1} \int_{\R^2} X_j \left|X\right|^2 \left|A\right|^2 \partial_j f f dX\right| &\leq \frac{\beta}{2} e^{-2\tau}\left\| \left|X\right|^2 \left|A\right| \nabla f \right\|^2\\
&\hspace{1cm} + C \beta e^{-2\tau}\left\|\left|X\right| f\right\|^2\left\| \nabla U \right\|^2_{L^\infty}.
}
Then, using the inequality (\ref{biotsinfini}) of Lemma \ref{biots1bis}, the inequality (\ref{weightLp}) of Lemma \ref{weight2} and the conditions (\ref{cond1bis}) and (\ref{cond2}), we get
\aligne{
\left|4\beta e^{-2\tau} \sum^{2}_{j=1} \int_{\R^2} X_j \left|X\right|^2 \left|A\right|^2 \partial_j f f dX\right| &\leq \frac{\beta}{2} e^{-2\tau}\left\| \left|X\right|^2 \left|A\right| \nabla f \right\|^2 \\
&\hspace{1cm}+C \beta e^{-2\tau}\left\|\left|X\right|^2 f\right\|\left\| f\right\|\left\| \nabla W \right\|_{H^1}\left\| \nabla W \right\|_{L^2(2)}\\
&\leq \frac{\beta}{2} e^{-2\tau}\left\| \left|X\right|^2 \left|A\right| \nabla f \right\|^2 + C M^2 \gamma^2\left(1-\theta\right)^{12} e^{-\tau},
}
and we conclude that
\eq{\label{E6J1}
J_1 \geq \frac{\beta}{2} e^{-2\tau}\left\| \left|X\right|^2 \left|A\right| \nabla f \right\|^2-C M^2 \gamma^2\left(1-\theta\right)^{12} e^{-\tau}.
}
By the same way, we estimate $J_2$. A short computation shows that
\cent{$
\displaystyle J_2=\beta \alpha_1 e^{-3\tau} \left\|\left|X\right|^2\left|A\right|\Delta f\right\|^2+2\beta \alpha_1 e^{-3\tau} \sum^{2}_{j=1}\int_{\R^2}\left|X\right|^4 \partial_j A : A \partial_j f \Delta f dX.
$}
We define
\cent{
$
\displaystyle I = \left|2\beta \alpha_1 e^{-3\tau} \sum^{2}_{j=1}\int_{\R^2}\left|X\right|^4 \partial_j A : A \partial_j f \Delta f dX\right|
$
}
Applying Hölder inequalities  and the continuous injection of $H^1(\R^2)$ into $L^4(\R^2)$, we obtain
\aligne{
I & \leq C\beta \alpha_1 e^{-3\tau} \left\|\left|X\right|^2 \left|A\right| \Delta f \right\| \left\|\left|X\right|^2 \nabla f \right\|_{L^4}\left\|\nabla^2 U\right\|_{L^4}\\
&\leq C\beta \alpha_1 e^{-3\tau} \left\|\left|X\right|^2 \left|A\right| \Delta f \right\| \left\|\left|X\right|^2 \nabla f \right\|_{L^4}\left\|\nabla^2 U\right\|_{H^1}.
}
Using the inequality (\ref{weightnabla2L4}) of Lemma \ref{weight2} and the inequality (\ref{biotsHs}) of Lemma \ref{biots1bis}, we get
\aligne{
I&\leq C\beta \alpha_1 e^{-3\tau} \left\|\left|X\right|^2 \left|A\right| \Delta f \right\| \left\|\left|X\right|^2 \nabla f\right\|^{1/2}\\
&\hspace{2cm}\times\left(\left\|f\right\|^{1/2}+\left\|\left|X\right|\nabla f\right\|^{1/2}+\left\|\left|X\right|^2\Delta f\right\|^{1/2}\right)\left\|W\right\|_{H^2}
.
}
Due to the Young inequality and the condition (\ref{cond1bis}), we obtain
\aligne{
I&\leq \frac{\beta }{2}\alpha_1 e^{-3\tau} \left\|\left|X\right|^2 \left|A\right| \Delta f \right\|^2 \\
&\hspace{2cm}+ C \beta \alpha_1 e^{-3\tau}\left\|W\right\|^2_{H^2}\left(\left\|f\right\|^2+\left\|\nabla f\right\|^2+\left\|\left|X\right|^2\nabla f\right\|^2+\left\|\left|X\right|^2\Delta f\right\|^2\right)\\
&\leq \frac{\beta }{2} \alpha_1 e^{-3\tau} \left\|\left|X\right|^2 \left|A\right| \Delta f \right\|^2\\
&\hspace{2cm} + CM\gamma\left(1-\theta\right)^6 e^{-2\tau}\left(\left\|f\right\|^2+\left\|\nabla f\right\|^2+\left\|\left|X\right|^2\nabla f\right\|^2+\left\|\left|X\right|^2\Delta f\right\|^2\right).
}
Thus, we can conclude that
\eq{\label{E6J2}
\arr{l}{\displaystyle J_2 \geq \frac{\beta }{2} \alpha_1 e^{-3\tau} \left\|\left|X\right|^2 \left|A\right| \Delta f \right\|^2 \\
\hspace{2cm}- CM\gamma\left(1-\theta\right)^6 e^{-2\tau}\left(\left\|f\right\|^2+\left\|\nabla f\right\|^2+\left\|\left|X\right|^2\nabla f\right\|^2+\left\|\left|X\right|^2\Delta f\right\|^2\right).
}
}
Combining the inequalities (\ref{E6J1}) and (\ref{E6J2}) and going back to (\ref{eqE6.3}), we have shown that
\eq{
\label{E6beta}\arr{l}{
\displaystyle J \geq \frac{\beta}{2} e^{-2\tau} \left(\left\|\left|X\right|^2\left|A\right| \nabla f\right\|^2+\alpha_1 e^{-\tau}\left\|\left|X\right|^2\left|A\right| \Delta f\right\|^2\right) -C M^2 \gamma^2\left(1-\theta\right)^{12} e^{-\tau}\\
\hspace{2cm}\displaystyle - CM\gamma\left(1-\theta\right)^6 e^{-2\tau}\left(\left\|f\right\|^2+\left\|\nabla f\right\|^2+\left\|\left|X\right|^2\nabla f\right\|^2+\left\|\left|X\right|^2\Delta f\right\|^2\right).
}
}
Taking into account the inequality (\ref{E6beta}), the equality (\ref{eqE6.0}) becomes
\eq{
\arr{l}{\label{eqE6.4}
\displaystyle\partial_\tau E_6+\frac{1}{2}\left\|\left|X\right|^2 f\right\|^2+\left(1+\alpha_1e^{-\tau}\right)\left\|\left|X\right|^2 \nabla f\right\|^2+\left(\alpha_1 e^{-\tau}+\frac{\alpha^2_1}{2}e^{-2\tau}\right)\left\|\left|X\right|^2 \Delta f\right\|^2\\
\hspace{5cm}\displaystyle+\frac{\beta}{2} e^{-2\tau} \left(\left\|\left|X\right|^2\left|A\right| \nabla f\right\|^2+\alpha_1 e^{-\tau}\left\|\left|X\right|^2\left|A\right| \Delta f\right\|^2\right)\leq\\
\\
\hspace{1cm} \displaystyle C\varepsilon e^{-\tau} \left\|\left|X\right|\nabla f\right\|^2+C\varepsilon \alpha_1 e^{-2\tau} \left\|\left|X\right|\Delta f\right\|^2+\left(8+8\alpha_1 e^{-\tau}\right)\left\|\left|X\right| f\right\|^2\\
\displaystyle \hspace{2cm}+ CM\gamma\left(1-\theta\right)^6 e^{-2\tau}\left(\left\|f\right\|^2+\left\|\nabla f\right\|^2+\left\|\left|X\right|^2\nabla f\right\|^2+\left\|\left|X\right|^2\Delta f\right\|^2\right)\\
\hspace{7cm}+CM^2\gamma^2\left(1-\theta\right)^{12} e^{-\tau}+I_1+I_2+I_3+I_4+I_5.
}
}
It remains to estimate every $I_i$, $i=1,...,5$. Using the divergence free property of $K$, integrating by parts and using Hölder inequalities, we get
\aligne{
I_1&= -2 \sum^{2}_{j=1} \int_{\R^2}X_j \left| X \right|^2 K_j \left|f-\alpha_1 e^{-\tau}\Delta f \right|^2 dX\\
&\leq C \left\| K\right\|_{L^\infty}\left\| \left| X \right|^2 \left(f-\alpha_1 e^{-\tau}\Delta f \right)\right\| \left\| \left| X \right|\left(f-\alpha_1 e^{-\tau}\Delta f \right)\right\|.
}
The inequalities (\ref{biotsinfini}) of Lemma \ref{biots1} and (\ref{weightLp}) of Lemma \ref{weight2}, the Young inequality $ab \leq \frac{3}{4}a^{\frac{4}{3}}+\frac{1}{4}b^4$ and the inequality (\ref{cond2}) yield
\eq{\label{E6I1}
\arr{ll}{I_1 &\displaystyle \leq C\left\|f\right\|^{1/2}_{H^1}\left\|f\right\|^{1/2}_{L^2(2)}\left\|\left| X \right|^2 \left(f-\alpha_1 e^{-\tau}\Delta f \right)\right\|^{3/2}\left\|f-\alpha_1 e^{-\tau}\Delta f \right\|^{1/2}\\
&\displaystyle \leq   CM^{1/2}\gamma^{1/2} \left(1-\theta\right)^{3} \left(\left\|\left| X \right|^2 \left(f-\alpha_1 e^{-\tau}\Delta f \right)\right\|^2+\left\|f-\alpha_1 e^{-\tau}\Delta f \right\|^2\right)\\
&\displaystyle \leq CM^{1/2}\gamma^{1/2} \left(1-\theta\right)^{3}  \left(\left\|\left| X \right|^2 f\right\|^2+\alpha^2_1e^{-2\tau} \left\|\left| X \right|^2\Delta f\right\|^2+\left\| f\right\|^2+\alpha^2_1 e^{-2\tau} \left\|\Delta f\right\|^2\right).
}
}
Using the inequality (\ref{biots4}) of Lemma \ref{biots1bis}, one can bound $I_2$ in a convenient way. Indeed, one has
\eq{\label{E6I2}
\arr{ll}{I_2 &\displaystyle \leq C \left|\eta\right|\left\|K\right\|_{L^4} \left\|\left|X\right|^2 \nabla\left(G-\alpha_1 e^{-\tau} \Delta G\right)\right\|_{L^4} \left\|\left| X \right|^2 \left(f-\alpha_1 e^{-t}\Delta f \right) \right\|\\
&\displaystyle \leq C \left|\eta\right|\left\|f\right\|_{L^2(2)} \left(\left\|\left| X \right|^2 f\right\|+\alpha_1e^{-\tau} \left\|\left| X \right|^2\Delta f\right\|\right)\\
&\displaystyle \leq CM^{1/2}\gamma^{1/2} \left(1-\theta\right)^{3}  \left(\left\|\left| X \right|^2 f\right\|^2+\alpha^2_1e^{-2\tau} \left\|\left| X \right|^2\Delta f\right\|^2+\left\| f\right\|^2\right).
}
}
Via an integration by parts, due to the facts that $V(X).X=0$ and $\dv V=0$, we show that $I_3$ vanishes. Indeed
\aligne{
I_3 &= \frac{\eta}{2}\sum^{2}_{j=1}\int_{\R^2} \left|X\right|^4 V_j\partial_j \left(\left|f-\alpha_1 e^{-\tau} \Delta f\right|^2\right) dX \\
&= -2\eta\sum^{2}_{j=1}\int_{\R^2} \left|X\right| X_j V_j\left|f-\alpha_1 e^{-\tau} \Delta f\right|^2 dX =0.
}
We rewrite $I_4=I^1_4+I^2_4$, where
\aligne{
&I^1_4 = -\eta\alpha_1 e^{-\tau}\left( \frac{\varepsilon}{\alpha_1}\Delta^2 G+\Delta G+\frac{X}{2}.\nabla \Delta G,\left|X\right|^2\left(f-\alpha_1 e^{-\tau} \Delta f\right)\right)_{L^2},\\
&I^2_4 = -\eta\beta e^{-2\tau}\left(\dv\left(\left|A\right|^2\nabla G\right),\left|X\right|^4\left(f-\alpha_1 e^{-\tau}\Delta f\right)\right)_{L^2}.
}
It is easy, using the smoothness of $G$ and the inequality (\ref{cond2}), to see that
\aligne{
I^1_4 &\leq C\left|\eta\right|e^{-\tau}\left\| G\right\|_{H^3(3)}\left(\left\| f\right\|+\alpha_1 e^{-\tau}\left\|\Delta  f\right\|\right) \\
&\leq C M\gamma\left(1-\theta\right)^{6} e^{-\tau}.
}
The term $I^3_4$ is not really harder to estimate. Due to the inequality (\ref{biotsnabla}) of Lemma \ref{biots1}, the inequality (\ref{cond2}), the continuous injection of $H^1(\R^2)$ into $L^4(\R^2)$ and the inequality (\ref{cond1bis}), we get
\aligne{
I^2_4 &\leq \left|\eta\right|\beta e^{-2\tau}\left(\left\|\nabla U\right\|^2_{L^4} \left\|\left|X\right|^4 \Delta G\right\|_{L^\infty}+\left\|\nabla U\right\|_{L^4}\left\|\nabla^2 U\right\|_{L^4} \left\|\left|X\right|^4 \nabla G\right\|_{L^\infty}\right)\\
&\hspace{10cm}\times\left(\left\| f\right\|+\alpha_1 e^{-\tau}\left\|\Delta  f\right\|\right)\\
&\leq C\left|\eta\right|\beta e^{-2\tau}\left(\left\|W\right\|^2_{L^4} +\left\|W\right\|_{L^4}\left\|\nabla W\right\|_{L^4}\right)\left(\left\| f\right\|+\alpha_1 e^{-\tau}\left\|\Delta  f\right\|\right)\\
&\leq C \left|\eta\right| e^{-2\tau}\left(\left\|W\right\|^2_{H^1}+\left\|W\right\|_{H^1}\left\|W\right\|_{H^2}\right)\left(\left\| f\right\|+\alpha_1 e^{-\tau}\left\|\Delta  f\right\|\right)\\
&\leq CM^2\gamma^2 \left(1-\theta\right)^{12} e^{-\frac{3\tau}{2}}.\\
}
Thus, assuming $\gamma \leq 1$, the following inequality holds: 
\eq{\label{E6I4}
I_4 \leq CM^2\gamma\left(1-\theta\right)^6 e^{-\tau}.
}
It remains to estimate $I_5$, which is the term that does not appear in the second grade fluids equations. We rewrite \cent{$I_5=I^1_5+I^2_5$,} where
\aligne{
&I^1_5 = -\beta e^{-2\tau}\left(\dv\left(\nabla\left(\left|A\right|^2\right)\wedge A\right),\left|X\right|^4 f\right)_{L^2},\\
&I^2_5= \beta \alpha_1 e^{-3\tau}\left(\dv\left(\nabla\left(\left|A\right|^2\right)\wedge A\right),\left|X\right|^4\Delta f\right)_{L^2}.
}
We begin by estimating $I^1_5$. After some computations, we notice that we have to estimate two kinds of terms. In fact, one has
\aligne{
I^1_5 \leq I^{1,1}_5+I^{1,2}_5, 
}
where
\aligne{
&I^{1,1}_5 =C\beta e^{-2\tau}\int_{\R^2} \left|X\right|^4 \left|\nabla^2 U\right|^2\left|\nabla U\right|\left|f\right|dX,\\
&I^{1,2}_5 = C\beta e^{-2\tau}\int_{\R^2} \left|X\right|^4 \left|\nabla^3 U\right|\left|\nabla U\right|^2\left|f\right| dX.
}
In order to simplify the notations, we set
\cent{$\delta = M \gamma \left(1-\theta\right)^6$.}
Applying the inequality (\ref{weightnabla2L4}) of Lemma \ref{weight2} and the continuous injection of $H^2(\R^2)$ into $L^\infty(\R^2)$, we obtain
\aligne{
I^{1,1}_5 &\leq \beta e^{-2\tau}\left\|\left|X\right| \nabla^2 U\right\|^2_{L^4}\left\|\nabla U\right\|_{L^\infty}\left\|\left|X\right|^2 f\right\|\\
&\leq C\beta e^{-2\tau}\left(\left\|W\right\|+\left\|\left|X\right| \nabla W\right\|\right)\left(\left\|\nabla W\right\|+\left\|\left|X\right| \Delta W\right\|\right)\left\|\nabla U\right\|_{H^2}\left\|\left|X\right|^2 f\right\|.}
Then, using the inequalities (\ref{weightLp}) of Lemma \ref{weight2} and (\ref{biotsHs}) of Lemma \ref{biots1bis} and the conditions (\ref{cond1bis}) and (\ref{cond2}), we get
\aligne{
I^{1,1}_5&\leq C\beta e^{-2\tau}\left(\left\|W\right\|+\left\|\nabla W\right\|^{1/2}\left\|\left|X\right|^2 \nabla W\right\|^{1/2}\right)\\
&\hspace{4cm}\times\left(\left\|\nabla W\right\|+\left\|\Delta W\right\|^{1/2}\left\|\left|X\right|^2 \Delta W\right\|^{1/2}\right)\left\|W\right\|_{H^2}\left\|\left|X\right|^2 f\right\|\\
&\leq C \delta e^{-3\tau/2}\left(\delta^{1/2}+\delta^{1/4}\left\|\left|X\right|^2 \nabla W\right\|^{1/2}\right)\left(\delta^{1/2}+\delta^{1/4}e^{\tau/4}\left\|\left|X\right|^2 \Delta W\right\|^{1/2}\right).
}
Then, we recall that  $W=\eta G+f$. Due to the fact that $\left|\eta\right| \leq \delta^{1/2}$ and the smoothness of $G$, we obtain
\aligne{
I^{1,1}_5 &\leq C\delta e^{-3\tau/2}\left(\delta^{1/2}+\delta^{1/4}\left\|\left|X\right|^2 \nabla f\right\|^{1/2}\right)\left(\delta^{1/2}+\delta^{1/4}e^{\tau/4}\left\|\left|X\right|^2 \Delta f\right\|^{1/2}\right)\\
&\leq C\delta^2 e^{-3\tau/2}+C\delta^{7/4} e^{-3\tau/2}\left\|\left|X\right|^2 \nabla f\right\|^{1/2}+C\delta^{7/4} e^{-5\tau/4}\left\|\left|X\right|^2 \Delta f\right\|^{1/2}\\
&\hspace{7cm}+C\delta^{3/2} e^{-5\tau/4}\left\|\left|X\right|^2 \nabla f\right\|^{1/2}\left\|\left|X\right|^2 \Delta f\right\|^{1/2}.
}
Using the Young inequalities $ab \leq \frac{1}{4}a^4+\frac{3}{4}b^{4/3}$ and $ab \leq \frac{1}{3}a^3+\frac{2}{3}b^{3/2}$, the inequality (\ref{cond2}) and assuming $\gamma\leq 1$, we finally obtain
\eq{\label{E5I5.1}
\arr{ll}{I^{1,1}_5 &\leq C\delta^2 e^{-\frac{3\tau}{2}}+ C \delta^2 \left(\left\|\left|X\right|^2 \nabla f\right\|^2+\left\|\left|X\right|^2 \Delta f\right\|^2\right) \\
&\hspace{4cm}+ C \delta^{5/3} \left(e^{-2\tau}+e^{-\frac{5\tau}{3}}\right)+C\delta^{4/3} e^{-\frac{5\tau}{3}} \left\|\left|X\right|^2 \nabla f\right\|^{2/3}\\
&\leq C\delta^2 e^{-\frac{3\tau}{2}}+ C \delta^2 \left(\left\|\left|X\right|^2 \nabla f\right\|^2+\left\|\left|X\right|^2 \Delta f\right\|^2\right)\\
&\hspace{4cm} + C \delta^{5/3} \left(e^{-2\tau}+e^{-\frac{5\tau}{3}}\right)+C\delta e^{-\frac{5\tau}{2}}+ C\delta^2\left\|\left|X\right|^2 \nabla f\right\|^2\\
&\leq CM^2 \gamma \left(1-\theta\right)^6 e^{-\frac{3\tau}{2}}+ C M^2 \gamma^2 \left(1-\theta\right)^{12} \left(\left\|\left|X\right|^2 \nabla f\right\|^2+\left\|\left|X\right|^2 \Delta f\right\|^2\right).
}
}
In order to estimate $I^{1,2}_5$, we use the Hölder inequalities and obtain
\aligne{
I^{1,2}_5 \leq C\beta e^{-2\tau} \left\|\left|X\right|^2 \nabla^3 U\right\|\left\|\left|X\right|^2 f\right\|_{L^4}\left\|\nabla U\right\|^2_{L^8}.
}
Then, using the Galiardo-Niremberg inequality, we notice that
\aligne{
\left\|\left|X\right|^2 f\right\|_{L^4} & \leq \left\|\left|X\right|^2 f\right\|^{1/2}\left\|\nabla \left(\left|X\right|^2 f\right)\right\|^{1/2}\\
&\leq C\left\|\left|X\right|^2f\right\|^{1/2}\left(\left\|\left|X\right|f\right\|+\left\|\left|X\right|^2\nabla f\right\|\right)^{1/2}.
}
The inequalities (\ref{biotsnabla}) of Lemma \ref{biots1}, (\ref{biots3nabla3}) of Lemma \ref{biots2} and the continuous injection of $H^1(\R^2)$ into $L^8(\R^2)$ imply
\aligne{
I^{1,2}_5 &\leq C\beta e^{-2\tau} \left\|\left|X\right|^2 \nabla^3 U\right\|\left\|\left|X\right|^2 f\right\|_{L^4}\left\|\nabla U\right\|^2_{L^8}\\
&\leq C\beta e^{-2\tau}\left(\left\|W\right\|+\left\|\left|X\right|\nabla W\right\|+\left\|\left|X\right|^2\Delta W\right\|\right)\\
&\hspace{3cm}\times\left\|\left|X\right|^2f\right\|^{1/2}\left(\left\|\left|X\right|f\right\|+\left\|\left|X\right|^2\nabla f\right\|\right)^{1/2}\left\|W\right\|^2_{H^1}.}
Finally, using the conditions (\ref{cond1bis}) and (\ref{cond2}) and the Young inequality $ab \leq \frac{1}{4}a^4+\frac{3}{4}b^{4/3}$, we obtain
\eq{\label{E5I5.2}
\arr{ll}{I^{1,2}_5&\leq C\delta^{7/4} e^{-\frac{3\tau}{4}}\left\|\left|X\right|^2 f\right\|^{1/2}\\
&\leq C\delta^2 e^{-\tau}+C\delta \left\|\left|X\right|^2 f\right\|^{2}\\
&\leq CM^2\gamma^2 \left(1-\theta\right)^{12} e^{-\tau}+ CM\gamma\left(1-\theta\right)^6\left\|\left|X\right|^2 f\right\|^{2}.
}
}
Thus, combining the inequalities (\ref{E5I5.1}) and (\ref{E5I5.2}), we obtain
\eq{\label{E5I15}
I^1_5 \leq CM^2 \gamma \left(1-\theta\right)^6 e^{-\tau}+ C M^2 \gamma \left(1-\theta\right)^6 \left(\left\|\left|X\right|^2 f\right\|^{2}+\left\|\left|X\right|^2 \nabla f\right\|^2+\left\|\left|X\right|^2 \Delta f\right\|^2\right).
}
It remains to estimate $I^2_5$. Like in the case of $I^1_5$, we have to consider two kinds of terms. Indeed, one can show that
\cent{
$
I^2_5 \leq I^{2,1}_5+I^{2,2}_5,
$
}
where
\aligne{
&I^{2,1}_5=C\alpha_1\beta e^{-3\tau}\displaystyle \int_{\R^2} \left|X\right|^4 \left|\nabla^2 U\right|^2 \left|\nabla U\right| \left|\Delta f\right| dX,\\
&I^{2,2}_5=C\alpha_1\beta e^{-3\tau}\displaystyle \int_{\R^2} \left|X\right|^4 \left|\nabla^3 U\right| \left|\nabla U\right|^2 \left|\Delta f\right| dX.
}
With the same tools than the ones used to estimate $I^1_5$, one can bound $I^{2,1}_5$. Due to Hölder inequalities and the continuous injection of $H^2(\R^3)$ into $L^\infty(\R^2)$, one has
\aligne{
I^{2,1}_5 &\leq C\alpha_1\beta e^{-3\tau}\left\|\left|X\right|^2 \Delta f\right\|\left\|\left|X\right| \nabla^2 U\right\|^2_{L^4}\left\|\nabla U\right\|_{L^\infty}\\
&\leq C\alpha_1\beta e^{-3\tau}\left\|\left|X\right|^2 \Delta f\right\|\left\|\left|X\right| \nabla^2 U\right\|^2_{L^4}\left\|\nabla U\right\|_{H^2}.
}
Then, using the inequality (\ref{weightnabla2L4}) of Lemma \ref{weight2} and the inequality (\ref{biotsHs}) of Lemma \ref{biots1bis}, we obtain
\aligne{
I^{2,1}_5 &\leq C\alpha_1\beta e^{-3\tau}\left\|\left|X\right|^2 \Delta f\right\|\left(\left\|W\right\|+\left\|\left|X\right| \nabla W\right\|\right)\left(\left\|\nabla W\right\|+\left\|\left|X\right| \Delta W\right\|\right)\left\|W \right\|_{H^2}}
Finally, the condition (\ref{cond1bis}) and the Young inequality imply
\eq{\label{E5I5.3}\arr{ll}{
I^{2,1}_5 &\leq C\delta^{3/2} e^{-\tau}\left\|\left|X\right|^2 \Delta f\right\|\\
&\leq C M^2 \gamma^2\left(1-\theta\right)^{12}e^{-2 \tau}+C M \gamma \left(1-\theta\right)^6 \left\|\left|X\right|^2 \Delta f\right\|^2.}
}
Likewise, using the inequality (\ref{biots3nabla3}) of Lemma \ref{biots2} and the continuous injection of $H^{\frac{3}{2}}(\R^2)$ into $L^\infty(\R^2)$, we get
\aligne{
I^{2,2}_5 &\leq C\beta \alpha_1 e^{-3\tau} \left\|\left|X\right|^2 \Delta f\right\|\left\|\left|X\right|^2 \nabla^3 U\right\|\left\|\nabla U\right\|^2_{L^\infty}\\
&\leq C\beta \alpha_1 e^{-3\tau} \left\|\left|X\right|^2 \Delta f\right\|\left(\left\|W\right\|+\left\|\left|X\right| \nabla W\right\|+\left\|\left|X\right|^2 \Delta W\right\|\right)\left\|W\right\|^2_{H^{3/2}}\\
&\leq C \delta^{1/2} e^{-2\tau}\left\|\left|X\right|^2 \Delta f\right\|\left\|W\right\|^2_{H^{3/2}}.
}
Using the well-known interpolation inequality
\cent{
$
\left\|v\right\|_{H^{3/2}} \leq C \left\|v\right\|^{1/2}_{H^1}\left\|v\right\|^{1/2}_{H^2}$,\hspace{1cm} for every $v \in H^2(\R^2)$,
}
we obtain, using again the condition (\ref{cond1bis}) and the Young inequality,
\eq{\label{E5I5.4}\arr{ll}{
I^{2,2}_5 & \leq C \delta^{1/2} e^{-2\tau}\left\|\left|X\right|^2 \Delta f\right\|\left\|W\right\|_{H^1}\left\|W\right\|_{H^2}\\
&\leq C \delta^{3/2} e^{-3\tau/2}\left\|\left|X\right|^2 \Delta f\right\|\\
&\leq C M^2 \gamma^2 \left(1-\theta\right)^{12} e^{-3\tau}+ C M \gamma \left(1-\theta\right)^6 \left\|\left|X\right|^2 \Delta f\right\|^2.
}
}
Finally, the inequalities (\ref{E5I5.3}) and (\ref{E5I5.4}) imply
\eq{\label{E5I25}
I^2_5 \leq C M^2 \gamma^2 \left(1-\theta\right)^{12} e^{-3\tau}+ C M \gamma \left(1-\theta\right)^6 \left\|\left|X\right|^2 \Delta f\right\|^2.
}
Thus, combining the inequalities (\ref{E5I15}) and (\ref{E5I25}), we get
\eq{\label{E6I5}
I_5 \leq C M^2 \gamma \left(1-\theta\right)^6 e^{-\tau}+CM^2\gamma \left(1-\theta\right)^6\left(\left\|\left|X\right|^2 f\right\|^{2}+ \left\|\left|X\right|^2 \nabla f\right\|^2+ \left\|\left|X\right|^2 \Delta f\right\|^2\right).
}
Taking into account the inequalities (\ref{E6I1}), (\ref{E6I2}), (\ref{E6I4}) and (\ref{E6I5}) and going back to (\ref{eqE6.4}), one has
\eq{
\arr{l}{\label{eqE6.5}
\displaystyle\partial_\tau E_6+\frac{1}{2}\left\|\left|X\right|^2 f\right\|^2+\left(1+\alpha_1e^{-\tau}\right)\left\|\left|X\right|^2 \nabla f\right\|^2+\left(\alpha_1 e^{-\tau}+\frac{\alpha^2_1}{2}e^{-2\tau}\right)\left\|\left|X\right|^2 \Delta f\right\|^2\\
\hspace{10cm}\displaystyle-\left(8+8\alpha_1 e^{-\tau}\right)\left\|\left|X\right| f\right\|^2 \leq\\
\\
\hspace{1cm}\displaystyle C\varepsilon e^{-\tau} \left\|\left|X\right|\nabla f\right\|^2+C\varepsilon \alpha_1 e^{-2\tau} \left\|\left|X\right|\Delta f\right\|^2+C M^2\gamma \left(1-\theta\right)^{6}e^{-\tau}\\
\hspace{3cm}\displaystyle +CM^2\gamma^{1/2} \left(1-\theta\right)^{3}\left(\left\|f\right\|^2+\left\|\nabla f\right\|^2+\left\|\Delta f\right\|^2\right)\\
\hspace{4cm}+CM^2\gamma^{1/2} \left(1-\theta\right)^{3}\left(\left\|\left|X\right|^2 f\right\|^{2}+ \left\|\left|X\right|^2 \nabla f\right\|^2+ \left\|\left|X\right|^2 \Delta f\right\|^2\right).
}
}
Via the Young inequality and the condition (\ref{cond2}), it is easy to check that
\aligne{
C\varepsilon e^{-\tau} \left\|\left|X\right|\nabla f\right\|^2 +C\varepsilon \alpha_1 e^{-2\tau} \left\|\left|X\right|\Delta f\right\|^2&\leq \varepsilon^2 \left\|\left|X\right|^2 \nabla f\right\|^2 + Ce^{-2\tau} \left\|\nabla f\right\|^2 \\
&\hspace{2cm}+ \varepsilon^2 \left\|\left|X\right|^2\Delta f\right\|^2 + C \alpha^2_1 e^{-4\tau} \left\|\Delta f\right\|^2\\
&\leq \varepsilon^2 \left\|\left|X\right|^2 \nabla f\right\|^2 +\varepsilon^2 \left\|\left|X\right|^2\Delta f\right\|^2\\
&\hspace{4cm}+ C M\left(1-\theta\right)^6 e^{-2\tau}.
}
We assume that $\displaystyle \varepsilon^2\leq \min\left(\frac{1}{2},\frac{\alpha_1e^{-\tau_0}}{2}\right)$. The inequality (\ref{eqE6.5}) becomes
\eq{
\arr{l}{
\displaystyle \partial_\tau E_6+\frac{1}{2}\left\|\left|X\right|^2 f\right\|^2+\left(\frac{1}{2}+\alpha_1e^{-\tau}\right)\left\|\left|X\right|^2 \nabla f\right\|^2\\
\displaystyle \hspace{3cm}+\left(\frac{\alpha_1}{2} e^{-\tau}+\frac{\alpha^2_1}{2}e^{-2\tau}\right)\left\|\left|X\right|^2 \Delta f\right\|^2-8\alpha_1 e^{-\tau}\left\|\left|X\right| f\right\|^2\leq\\
\\
\hspace{1cm}\displaystyle  C M^2\gamma \left(1-\theta\right)^{6} e^{-\tau}+8\left\|\left|X\right| f\right\|^2\\
\hspace{2cm}\displaystyle+C_1M^2\gamma^{1/2} \left(1-\theta\right)^{3}\left(\left\|f\right\|^2+\left\|\nabla f\right\|^2+\left\|\Delta f\right\|^2\right)\\
\hspace{3cm}\displaystyle+C_1M^2\gamma^{1/2} \left(1-\theta\right)^{3}\left(\left\|\left|X\right|^2 f\right\|^{2}+ \left\|\left|X\right|^2 \nabla f\right\|^2+ \left\|\left|X\right|^2 \Delta f\right\|^2\right),
}
}
where $C_1$ is a positive constant dependent on $\alpha_1$ and $\beta$.
\vspace{0.5cm}\\
We take now $\gamma$ sufficiently small so that $C_1M^2\gamma^{1/2} \left(1-\theta\right)^{3} \leq \frac{1-\theta}{4}$. We obtain
\eq{
\arr{l}{
\displaystyle \partial_\tau E_6+\left(\frac{\theta}{2}+\frac{1-\theta}{4}\right)\left\|\left|X\right|^2 f\right\|^2+\left(\frac{1}{4}+\alpha_1e^{-\tau}\right)\left\|\left|X\right|^2 \nabla f\right\|^2\\
\displaystyle \hspace{3cm}+\left(\frac{\alpha_1}{4} e^{-\tau}+\frac{\alpha^2_1}{2}e^{-2\tau}\right)\left\|\left|X\right|^2 \Delta f\right\|^2-8\alpha_1 e^{-\tau}\left\|\left|X\right| f\right\|^2\leq \\
\\
\hspace{4cm}\displaystyle C M^2\gamma \left(1-\theta\right)^{6} e^{-\tau}+8\left\|\left|X\right| f\right\|^2\\
\hspace{6cm}\displaystyle+C_1M^2\gamma^{1/2} \left(1-\theta\right)^{3}\left(\left\|f\right\|^2+\left\|\nabla f\right\|^2+\left\|\Delta f\right\|^2\right).
}
}
Using the inequality (\ref{weightLp}) of Lemma \ref{weight2}, one has
\cent{
$
\displaystyle 8\left\|\left|X\right|  f\right\|^2 \leq h \left\|\left|X\right|^2  f\right\|^2 +\frac{64}{h}\left\|f\right\|^2,\quad \hbox{for all} \quad h>0.
$
}
Thus, we set $h= \frac{1-\theta}{8}$ and obtain
\eq{
\arr{l}{\label{eq6.6}
\displaystyle\partial_\tau E_6+\left(\frac{\theta}{2}+\frac{1-\theta}{8}\right)\left\|\left|X\right|^2 f\right\|^2+\left(\frac{1}{4}+\alpha_1e^{-\tau}\right)\left\|\left|X\right|^2 \nabla f\right\|^2\\
\displaystyle \hspace{3cm}+\left(\frac{\alpha_1}{4} e^{-\tau}+\frac{\alpha^2_1}{2}e^{-2\tau}\right)\left\|\left|X\right|^2 \Delta f\right\|^2-8\alpha_1 e^{-\tau}\left\|\left|X\right| f\right\|^2 \leq \\
\\
\displaystyle \hspace{3cm} C M^2\gamma \left(1-\theta\right)^{6} e^{-\tau}+\frac{1024}{1-\theta}\left\|f\right\|^2\\
\hspace{6cm}\displaystyle+C_1M^2\gamma^{1/2} \left(1-\theta\right)^{3}\left(\left\|f\right\|^2+\left\|\nabla f\right\|^2+\left\|\Delta f\right\|^2\right).
}
}
Integrating several times by parts, we notice that
\cent{$
E_6 = \frac{1}{2}\left\|\left|X\right|^2 f\right\|^2+\alpha_1 e^{-\tau}\left\|\left|X\right|^2\nabla f\right\|^2+\frac{\alpha^2_1}{2}e^{-2\tau}\left\|\left|X\right|^2 \Delta f\right\|^2-8\alpha_1e^{-\tau}\left\|\left|X\right|  f\right\|^2.
$}
Consequently, the inequality (\ref{eq6.6}) can be written
\eq{
\arr{l}{
\displaystyle \partial_\tau E_6+\theta E_6+\frac{1-\theta}{8}\left\|\left|X\right|^2 f\right\|^2+\frac{1}{4}\left\|\left|X\right|^2 \nabla f\right\|^2+\frac{\alpha_1}{4} e^{-\tau}\left\|\left|X\right|^2 \Delta f\right\|^2\displaystyle \\
\\
\displaystyle \hspace{0.5cm}\leq C M^2\gamma \left(1-\theta\right)^{6} e^{-\tau}+\frac{1024}{1-\theta}\left\|f\right\|^2+C M^2\gamma^{1/2} \left(1-\theta\right)^{3}\left(\left\|f\right\|^2+\left\|\nabla f\right\|^2+\left\|\Delta f\right\|^2\right).
}
}
\begin{flushright}
$\square$
\end{flushright}
\section{\label{secdem}Proof of Theorem \ref{theo1}}
In this section, we consider the solution $W_\varepsilon$ of (\ref{g3We}) with initial data $W_0$ satisfying the condition (\ref{cond1}) for some $\gamma>0$ and we take advantage of the energy estimates obtained in Section \ref{secenergy} to show that $W_\varepsilon$ satisfies the inequality (\ref{inetheo1}). Then, we pass to the limit when $\varepsilon$ tends to $0$ and show that $W_\varepsilon$ converges, up to a subsequence, to a weak solution of (\ref{g3W}) which satisfies also the inequality (\ref{inetheo1}). We recall that
\cent{$
W_\varepsilon = \eta G + f_\varepsilon,
$}
where $G$ is the Oseen vortex sheet given by (\ref{oseen}), $\displaystyle \eta =\int_{\R^2} W_0(X)dX$ and $f_\varepsilon$ satisfies the equality (\ref{g3fe}). We define the functional 
\cent{$\displaystyle E_7= \frac{K}{1-\theta} E_5 +E_6$,} where $K$ is a large positive constant that will be made more precise later and $E_5$ and $E_6$ are the energy functionals defined in Section \ref{secenergy}.
\vspace{0.5cm}\\
If $K$ is large enough, this energy is suitable to estimate the $H^2(2)$ norm of $f_\varepsilon$, as it is shown by the next lemma.
\lem{\label{E7H22}
Let $f_\varepsilon\in C^1\left(\left(\tau_0,\tau_\varepsilon\right),H^1(2)\right)\cap C^0\left(\left(\tau_0,\tau_\varepsilon\right),H^3(2)\right)$. If $K$ is large enough, there exist two positive constants $C_1$ and $C_2$ such that, for all $\tau \in \left(\tau_0,\tau_\varepsilon\right)$,
\aligne{
E_7 \leq \frac{C_1}{1-\theta} \left(\left\|f_\varepsilon\right\|^2_{H^1}+\alpha_1 e^{-\tau}\left\|\Delta f_\varepsilon\right\|^2+\left\|\left|X\right|^2 f_\varepsilon\right\|^2+\alpha^2_1  e^{-2\tau}\left\|\left|X\right|^2 \Delta f_\varepsilon\right\|^2\right),\\
C_2\left(\left\|f_\varepsilon\right\|^2_{H^1}+\alpha_1 e^{-\tau}\left\|\Delta f_\varepsilon \right\|^2+\left\|\left|X\right|^2 f_\varepsilon\right\|^2+\alpha^2_1  e^{-2\tau}\left\|\left|X\right|^2 \Delta f_\varepsilon\right\|^2\right) \leq E_7.\\
}
}
\textbf{Proof: }The first inequality of this lemma comes directly from the definition of $E_7$. To prove the second one, we notice that
\cent{$
\displaystyle E_7\geq \frac{C K}{1-\theta}\left(\left\|f_\varepsilon\right\|^2_{H^1}+\alpha_1 e^{-\tau}\left\|\Delta f_\varepsilon \right\|^2\right)+\frac{1}{2}\left\|\left|X\right|^2\left( f_\varepsilon-\alpha_1 e^{-\tau} \Delta f_\varepsilon\right)\right\|^2.
$}
Furthermore, we have already shown that
\cent{$
\arr{l}{\displaystyle \left\|\left|X\right|^2\left( f_\varepsilon-\alpha_1 e^{-\tau} \Delta f_\varepsilon\right)\right\|^2=\left\|\left|X\right|^2 f_\varepsilon\right\|^2+2\alpha_1  e^{-\tau}\left\|\left|X\right|^2 \nabla f_\varepsilon\right\|^2\\
\hspace{7cm}+\alpha^2_1  e^{-2\tau}\left\|\left|X\right|^2 \Delta f_\varepsilon\right\|^2- 16 \left\|\left|X\right|f_\varepsilon\right\|^2.
}$}
Via the Hölder and Young inequalities, we get
\cent{$
\arr{l}{\displaystyle\left\|\left|X\right|^2\left( f_\varepsilon-\alpha_1 e^{-\tau} \Delta f_\varepsilon\right)\right\|^2 \geq \left\|\left|X\right|^2 f_\varepsilon\right\|^2+2\alpha_1  e^{-\tau}\left\|\left|X\right|^2 \nabla f_\varepsilon\right\|^2\\
\hspace{5cm}\displaystyle+\alpha^2_1  e^{-2\tau}\left\|\left|X\right|^2 \Delta f_\varepsilon\right\|^2- \frac{1}{2} \left\|\left|X\right|^2 f_\varepsilon\right\|^2-128 \left\| f_\varepsilon\right\|^2.}$
}
Consequently, one has
\cent{$
\arr{l}{\displaystyle E_7 \geq \frac{C K}{1-\theta}\left(\left\|f_\varepsilon\right\|^2_{H^1}+\alpha_1 e^{-\tau}\left\|\Delta f_\varepsilon \right\|^2\right)+\frac{1}{4}\left\|\left|X\right|^2 f_\varepsilon\right\|^2+\alpha_1  e^{-\tau}\left\|\left|X\right|^2 \nabla f_\varepsilon\right\|^2\\
\hspace{9cm}\displaystyle +\frac{\alpha^2_1}{2}  e^{-2\tau}\left\|\left|X\right|^2 \Delta f\right\|^2 -64 \left\| f_\varepsilon\right\|^2.}
$}
Thus, if $K$ is big enough, we get the second inequality of this lemma.
\begin{flushright}
$\square$
\end{flushright}
\lem{\label{lemE7}
Let $W_\varepsilon\in C^0\left(\left[\tau_0,\tau_\varepsilon\right),H^3(2)\right)$ be a solution of $(\ref{g3We})$ satisfying the inequality (\ref{cond1bis}) for some $\gamma>0$. There exist $T_0>0$ and $\gamma_0 >0$  such that if $T=e^{\tau_0} \geq T_0$ and $\gamma\leq \gamma_0$, then, for all $\tau\in \left[\tau_0,\tau^*_\varepsilon\right)$, $E_7$ satisfies the inequality
\eq{\label{E7}
\partial_\tau E_7+\theta E_7 \leq CM^3 \gamma \left(1-\theta\right) e^{-\tau}.
}
}
\textbf{Proof: }
We take $\gamma_0$ and $T_0$ respectively as small and large as necessary to satisfy the conditions of the lemmas \ref{lemE1} to \ref{lemE6}. According to the inequalities $(\ref{E5})$ and $(\ref{E6})$, one has
\aligne{
&\partial_\tau E_7+\theta E_7+\frac{K}{1-\theta}\left(7\left\|\left(-\Delta\right)^{-\frac{1+\theta}{4}}f_\varepsilon\right\|^2+\frac{1}{4}\left\|\nabla f_\varepsilon\right\|^2+\frac{1}{4}\left\|\Delta f_\varepsilon\right\|^2\right)\\
&\hspace{1cm} + \frac{1-\theta}{8}\left\|\left|X\right|^2 f_\varepsilon\right\|^2+\frac{\alpha_1}{4}e^{-\tau}\left\|\left|X\right|^2\Delta  f_\varepsilon\right\|^2 \leq CM^3 \gamma \left(1-\theta\right) e^{-\tau}+\frac{1024}{1-\theta}\left\| f_\varepsilon\right\|^2\\
&\hspace{5cm}+C M^2 \left(1-\theta\right) \gamma^{1/2}K\left(\left\|\left|X\right|^2 f_\varepsilon\right\|^2+\alpha^2_1 e^{-2\tau}\left\|\left|X\right|^2\Delta f_\varepsilon\right\|^2\right)\\
&\hspace{7cm}+C M^2 \left(1-\theta\right) \gamma^{1/2}\left(\left\|f_\varepsilon\right\|^2+\left\|\nabla f_\varepsilon\right\|^2+\left\|\Delta f_\varepsilon\right\|^2\right).\\
}
Using the interpolation inequality (\ref{interpolation}) of $\left\|f_\varepsilon\right\|^2$ between $\left\|\left(-\Delta\right)^{-\frac{1+\theta}{4}} f_\varepsilon\right\|^2$ and $\left\|\nabla f_\varepsilon\right\|^2$ and taking $K$ large enough and $\gamma$ small enough, we get
\eq{
\partial_\tau E_7+\theta E_7 \leq CM^3 \gamma \left(1-\theta\right) e^{-\tau}.
}
\begin{flushright}
$\square$
\end{flushright}
\begin{remark}
We can see in the proofs of the lemmas \ref{lemE1} to \ref{lemE7} that $\gamma_0$ does not depend on $\theta$, but only on $\alpha_1$, $\beta$ and $M$.
\end{remark}
\subsection{Regularized problem}
Before proving Theorem \ref{theo1}, we show an intermediate theorem. This one gives the same result than Theorem \ref{theo1}, but for the solutions of the regularized system (\ref{g3We}).
\theo{\label{theoe} Let $\theta$ be a constant such that $0<\theta \leq 1$. There exist $\varepsilon_0=\varepsilon_0(\alpha_1,\beta)>0$, $\gamma_0=\gamma_0(\alpha_1,\beta)>0$ and $T_0=T_0(\alpha_1,\beta)\geq 0$ such that, for all $\varepsilon \leq \varepsilon_0$, $T=e^{\tau_0} \geq T_0$ and $W_0\in H^2(2)$ satisfying the condition $(\ref{cond1})$ with $\gamma \leq \gamma_0$, there exist a unique global solution $W_\varepsilon\in C^1\left(\left(\tau_0,+\infty\right),H^1(2)\right)\cap C^0\left(\left(\tau_0,+\infty\right),H^3(2)\right)$ of $(\ref{g3We})$ and a positive constant $C=C(\alpha_1,\beta,\theta)>0$ such that, for all $\tau\geq \tau_0$,
\eq{\label{th3}
\left\|\left(1-\alpha_1e^{-\tau}\Delta\right)\left(W_\varepsilon(\tau)-\eta G\right)\right\|^2_{L^2(2)}\leq C\gamma e^{-\theta \tau},
}
where $\eta=\displaystyle \int_{\R^2}W_0(x) dx$ and the parameters $\alpha_1$ and $\beta$ are fixed and given in (\ref{g3u}).
}
\paragraph{Proof of Theorem \ref{theoe}:}
Let $W_0 \in H^2(2)$ satisfying the condition $(\ref{cond1})$ with $0\leq \gamma \leq \gamma_0$ and $0\leq T_0 \leq T$, where $\gamma_0$ and $T_0$ will be made more precise later. By theorem \ref{theo2}, there exist $\tau_\varepsilon > \tau_0=\log(T)$ and a solution $W_\varepsilon$ to the system $(\ref{g3We})$ which belongs to $C^1\left(\left(\tau_0,\tau_\varepsilon\right),H^1(2)\right)\cap C^0\left(\left(\tau_0,\tau_\varepsilon\right),H^3(2)\right)$. Let $\eta =\displaystyle \int_{\R^2} W_0(X) dX$, and $f_\varepsilon$ defined by the equality
\eq{\label{decomp2}
W_\varepsilon = \eta G + f_\varepsilon.
}
Let $M>2$ be a positive constant that will be set later and $\tau^*_\varepsilon \in \left[\tau_0,\tau_\varepsilon\right)$ be the highest positive time such that the inequality (\ref{cond1bis}) holds. As shown at the beginning of Section \ref{secenergy}, the inequality (\ref{cond2}) holds on $\left[\tau_0,\tau^*_\varepsilon\right)$. We take $T_0$ sufficiently large and $\gamma_0$ and $\varepsilon$ sufficiently small so that the results of the lemmas \ref{lemE1} to \ref{lemE7} occur. Consequently, there exists $C=C(\alpha_1,\beta)>0$ such that, for all $\tau \in \left[\tau_0,\tau^*_\varepsilon\right)$,
\eq{
\partial_\tau \left(E_7e^{\theta \tau}\right)\leq CM^3 \gamma \left(1-\theta\right) e^{-\left(1-\theta\right)\tau}.
}
Integrating this inequality in time between $\tau_0$ and $\tau \in \left[\tau_0,\tau^*_\varepsilon\right)$, we obtain
\eq{\label{ineqdemo1}
\displaystyle E_7 (\tau) \leq E_7 (\tau_0) e^{-\theta\left(\tau-\tau_0\right)}+CM^3 \gamma \left(e^{-\left(1-\theta\right) \tau_0}e^{-\theta \tau}-e^{-\tau}\right).
}
Due to the decomposition (\ref{decomp2}) and the lemma \ref{E7H22}, for every $\tau \in \left[\tau_0,\tau^*_\varepsilon\right)$, one has
\aligne{
\left\|W_\varepsilon(\tau) \right\|^2_{H^1}+ \left\|\left|X\right|^2 W_\varepsilon(\tau)\right\|^2 + \alpha_1 e^{-\tau} \left\|\Delta W_\varepsilon(\tau)\right\|^2 + \alpha^2_1 e^{-2\tau} \left\|\left|X\right|^2 \Delta W_\varepsilon(\tau)\right\|^2 &\leq\\
&\hspace{-1cm} C\eta^2 + C E_7(\tau).
}
Since $f_\varepsilon$ satisfies the inequality (\ref{cond2}), one has $\eta^2 \leq C \gamma \left(1-\theta\right)^6$. Taking into account the inequality (\ref{ineqdemo1}), it comes
\eq{\label{ineqdemo2}
\arr{l}{\displaystyle\left\|W_\varepsilon(\tau) \right\|^2_{H^1}+ \left\|\left|X\right|^2 W_\varepsilon(\tau)\right\|^2 + \alpha_1 e^{-\tau} \left\|\Delta W_\varepsilon(\tau)\right\|^2 + \alpha^2_1 e^{-2\tau} \left\|\left|X\right|^2 \Delta W_\varepsilon(\tau)\right\|^2 \leq \\
\hspace{7.5cm} \displaystyle C\gamma \left(1-\theta\right)^6 + E_7 (\tau_0) e^{-\theta\left(\tau-\tau_0\right)}+CM^3 \gamma e^{-\tau_0}.
}
}
Using again the lemma \ref{E7H22} and arguing like for the establishment of the inequality (\ref{cond2}), we can show that
\cent{
$
\arr{ll}{\displaystyle E_7 (\tau_0) &\displaystyle\leq \frac{C}{1-\theta} \Big(\left\|f_\varepsilon(\tau_0)\right\|^2_{H^1}+\alpha_1 e^{-\tau_0}\left\|\Delta f_\varepsilon(\tau_0)\right\|^2\\
&\displaystyle \hspace{3cm}+\left\|\left|X\right|^2 f_\varepsilon(\tau_0)\right\|^2+\alpha^2_1  e^{-2\tau_0}\left\|\left|X\right|^2 \Delta f_\varepsilon(\tau_0)\right\|^2\Big)\\
& \displaystyle \leq  C\gamma \left(1-\theta\right)^5.
}$
}
Consequently, the inequality (\ref{ineqdemo2}) becomes
\eq{
\arr{l}{\displaystyle\left\|W_\varepsilon(\tau) \right\|^2_{H^1}+ \left\|\left|X\right|^2 W_\varepsilon(\tau)\right\|^2 + \alpha_1 e^{-\tau} \left\|\Delta W_\varepsilon(\tau)\right\|^2 + \alpha^2_1 e^{-2\tau} \left\|\left|X\right|^2 \Delta W_\varepsilon(\tau)\right\|^2 \leq \\
\hspace{10cm} \displaystyle C_1\gamma \left(1-\theta\right)^5 +C_2 M^3 \gamma e^{-\tau_0},
}
}
where $C_1$ and $C_2$ are two positive constants independent of $W_0$ and $\theta$.
\vspace{0.5cm}\\
We set $M = \frac{4 C_1}{1-\theta}$, and we get
\eq{
\arr{l}{\displaystyle\left\|W_\varepsilon(\tau) \right\|^2_{H^1}+ \left\|\left|X\right|^2 W_\varepsilon(\tau)\right\|^2 + \alpha_1 e^{-\tau} \left\|\Delta W_\varepsilon(\tau)\right\|^2 + \alpha^2_1 e^{-2\tau} \left\|\left|X\right|^2 \Delta W_\varepsilon(\tau)\right\|^2 \leq \\
\hspace{10cm} \displaystyle \frac{M\gamma \left(1-\theta\right)^6}{4} +C_2 M^3 \gamma e^{-\tau_0}.
}
}
Finally, taking $T_0$ sufficiently large so that $\displaystyle C_2 M^3 \gamma e^{-\tau_0} \leq \frac{M\gamma \left(1-\theta\right)^6}{4}$, we obtain, for all $\tau \in \left[\tau_0, \tau^*_\varepsilon\right)$,
\eq{\label{ineqdemo3}
\displaystyle\left\|W_\varepsilon(\tau) \right\|^2_{H^1}+ \left\|\left|X\right|^2 W_\varepsilon(\tau)\right\|^2 + \alpha_1 e^{-\tau} \left\|\Delta W_\varepsilon(\tau)\right\|^2 + \alpha^2_1 e^{-2\tau} \left\|\left|X\right|^2 \Delta W_\varepsilon(\tau)\right\|^2 \leq \frac{M\gamma \left(1-\theta\right)^6}{2}.
}
This inequality shows in particular that $\tau^*_\varepsilon = \tau_\varepsilon$ and thus (\ref{ineqdemo3}) holds for all $\tau \in \left[\tau_0,\tau_\varepsilon\right)$. From the inequality (\ref{ineqdemo3}), we deduce also that $\tau_\varepsilon = +\infty$. Indeed, if $\tau_\varepsilon<+\infty$, the boundedness of $W_\varepsilon$ in $H^2(2)$ on $\left[\tau_0,\tau_\varepsilon\right)$ given by (\ref{ineqdemo3}) is a contradiction to the finiteness of $\tau_\varepsilon$.
\vspace{0.5cm}\\
In particular, the inequality (\ref{ineqdemo1}) occurs on $\left[\tau_0,+\infty\right)$. Applying the lemma \ref{E7H22} in the inequality (\ref{ineqdemo1}), we finally obtain the inequality (\ref{th3}).
\subsection{\label{secexis}Existence of weak solutions in $H^2(2)$}
Now, we show that under the hypotheses of Theorem \ref{theoe}, there exists a global weak solution $W$ of (\ref{g3W}) which belongs to $ C^0\left(\left[\tau_0,+\infty\right),H^2(2)\right)$, and that this solution converges to the Oseen vortex sheet $G$ when $\tau$ goes to infinity. To this end, we pass to the limit in the system $(\ref{g3We})$ when $\varepsilon$ tends to $0$ and show that, up to a subsequence, $W_\varepsilon$ converges in some sense to a solution of the system $(\ref{g3W})$ which satisfies the inequality $(\ref{th3})$. Let $(\varepsilon_n)_{n\in \mathbb N}$ be a sequence of positive numbers tending to $0$. We consider the solution $W_{\varepsilon_n} \in C^1\left(\left(\tau_0,+\infty\right),H^1(2)\right)\cap C^0\left(\left(\tau_0,+\infty\right),H^3(2)\right)$ of $(\ref{g3We})$ which satisfies the conditions of Theorem \ref{theoe}. Due to technical reasons linked to the compactness properties of Sobolev spaces, it is more convenient to establish the convergence of $W_{\varepsilon_n}$ to $W$ in every bounded regular domain of $\R^2$. Let $\Omega$ be a bounded regular domain of $\R^2$ and $\tau_1$ be a fixed positive time such that $\tau_0<\tau_1<+\infty$. In what follows, $H^s(\Omega)$, $s\geq 0$, denotes the restrictions to $\Omega$ of the functions of the Sobolev space $H^s(\R^2)$. From Theorem \ref{theoe}, we know that $W_{\varepsilon_n}$ is bounded in $L^\infty \left(\left[\tau_0,+\infty\right),H^2(2)\right)$ uniformly with respect to $n$. Consequently, there exists $W \in L^\infty\left(\left[\tau_0,\tau_1\right],H^2(2)\right)$ such that
\cent{$
W_{\varepsilon_n} \rightharpoonup W \quad \hbox{weakly in}\quad L^p\left(\left[\tau_0,\tau_1\right],H^2(\Omega)\right),
$ for all $p\geq 2$.}
Looking at the system (\ref{g3We}), we can see that $\partial_\tau W_{\varepsilon_n}$ is bounded in $L^\infty\left(\left[\tau_0,\tau_1\right],H^1(\Omega)\right)$ uniformly with respect to $n$. This implies that $W_{\varepsilon_n}$ is equicontinuous in $H^1(\Omega)$. Indeed, for $\sigma_1,\sigma_2 \in \left[\tau_0,\tau_1\right]$, $ \sigma_2 \geq \sigma_1$, we have
\aligne{
\left\|W_{\varepsilon_n}(\sigma_2) - W_{\varepsilon_n}(\sigma_1)\right\|_{H^1(\Omega)} &\displaystyle= \left\|\int^{\sigma_2}_{\sigma_1} \partial_\tau W_{\varepsilon_n}(s)ds\right\|_{H^1(\Omega)}\\
&\leq \left(\sigma_2-\sigma_1\right) \left\|\partial_\tau W_{\varepsilon_n}(s)\right\|_{L^\infty\left(\left[\tau_0,\tau_1\right],H^1(\Omega)\right)}.\\
} 
Furthermore, for every $\tau\in \left[\tau_0,\tau_1\right]$, the set $\bigcup\limits_{n \in \mathbb N} f_{\epsilon_n}(\tau)$ is bounded in $H^2(\Omega)$ and thus compact in $H^1(\Omega)$. Using the Arzela-Ascoli theorem, we get
\cent{$
W_{\varepsilon_n} \rightarrow W \quad \hbox{ strongly in} \quad C^0\left(\left[\tau_0,\tau_1\right],H^1(\Omega)\right)
$.}
By interpolation, we can show that 
\eq{\label{conv1}W_{\varepsilon_n} \rightarrow W\quad \hbox{in}\quad C^0\left(\left[\tau_0,\tau_1\right],H^s(\Omega)\right), \quad \hbox{for all}\quad s<2.}
This is enough to pass to the limit in the system $(\ref{g3We})$ in the sense of the distributions on $\left[\tau_0,\tau_1\right]\times\Omega$ and to show that $W$ is a weak solution of the system $(\ref{g3W})$. Since most of the terms of the equation (\ref{g3We}) have already been studied in \cite{jaffal11}, we will just show that the convergence holds for the term $\quad -\dv \curl\left(\left|A_{\varepsilon_n}\right|^2 A_{\varepsilon_n} \right)\quad$ which does not appear in the second grade fluids equations.\\
\vspace{0.5cm}\\
We consider $\varphi\in C^\infty_0 \left(\left[\tau_0,\tau_1\right]\times \Omega\right)$. For all $\tau\in \left[\tau_0,\tau_1\right]$, we want to show that
\eq{\label{conv3}
\arr{l}{\displaystyle \int^\tau_{\tau_0} \int_{\Omega} \left|A_{\varepsilon_n}(\tau,X)\right|^2  A_{\varepsilon_n}(\tau,X) \diamond \nabla^2 \varphi (\tau,X) dX d\tau \longrightarrow \\
\hspace{5cm}\displaystyle\int^\tau_{\tau_0} \int_{\Omega} \left|A(\tau,X)\right|^2 A(\tau,X) \diamond \nabla^2 \varphi(\tau,X) dX d\tau,
}
}
when $n$ tends to infinity, where, for $A,B \in \mathcal M_2(\R)$, we use the notation
\cent{
$
\displaystyle A \diamond B = \sum^2_{j=1} \left(A_{1,j} B_{2,j}- A_{2,j}B_{1,j}\right).
$
}
The term of the right hand side of (\ref{conv3}) appears naturally via two integrations by parts, when performing the $L^2-$scalar product of $-\dv \curl \left(\left|A\right|^2 A\right)$ with $\varphi$. The strong convergence of $W_{\varepsilon_n}$ to $W$ in $C^0\left(\left[\tau_0,\tau_1\right],H^1(\Omega)\right)$ implies directly the identity (\ref{conv3}). Indeed, due to the continuous injection of $H^1(\Omega)$ into $L^3(\Omega)$, $W_{\varepsilon_n}$ converges to $W$ in $C^0\left(\left[\tau_0,\tau_1\right],L^3(\Omega)\right)$. Furthermore, the inequality (\ref{biotsnabla}) implies
\cent{
$
\left\|A_{\varepsilon_n}-A\right\|_{L^3} \leq \left\|W_{\varepsilon_n}-W\right\|_{L^3},
$
}
and consequently $A_{\varepsilon_n}$ converges to $A$ strongly in $C^0\left(\left[\tau_0,\tau_1\right],L^3(\Omega)\right)$. This fact suffices to show that the identity (\ref{conv3}) occurs. Thus $W$ is a global weak solution of (\ref{g3W}) which belongs to $C^0\left(\left[\tau_0,+\infty\right),H^2(2)\right)$.
\vspace{0.5cm}\\
The fact that $W$ satisfies the inequality (\ref{inetheo1}) is a direct consequence of the weak convergence of $W_{\varepsilon_n}$ to $W$. Indeed, for all $\tau \in \left[\tau_0,+\infty\right)$,  $W_{\varepsilon_n}(\tau)$ is bounded in $H^2(2)$ uniformly with respect to $n$ and consequently we have
\cent{
$
W_{\varepsilon_n}(\tau) \rightharpoonup W(\tau)
$, weakly in $H^2(2)$, for all $\tau \in \left[\tau_0,+\infty\right)$.
}
Since $W_{\varepsilon_n}$ satisfies the inequality (\ref{inetheo1}), it implies that $W$ also satisfies (\ref{inetheo1}).
\subsection{Uniqueness}
The aim of this part is to prove that the solution $w$ of the system $(\ref{g3w})$ obtained in Section \ref{secexis} is unique in $L^2(2)$. Let $w_1$ and $w_2$ be two solutions of $(\ref{g3w})$ with the same initial data $w_0 \in H^2(2)$. Let $u_1$ and $u_2$ be the divergence free vector fields obtained via the Biot-Savart law respectively from $w_1$ and $w_2$. We also define $A_i = \nabla u_i+\left(\nabla u_i\right)^t$. Applying the Biot-Savart law to the system (\ref{g3u}), we can see that, for $i=1,2$, the divergence free vector field $u_i$ satisfies the system
\eq{\label{g3diffui}
\arr{l}{
\partial_t\left(u_i-\alpha_1 \Delta u_i\right)-\Delta u_i+\curl\left(u_i-\alpha_1 \Delta u_i\right)\wedge u_i-\beta \dv \left(\left|A_i\right|^2 A_i\right)+\nabla p_i = 0,\\
\dv u_i=0,\\
u_{i\left|t=0\right.} = u_0,
}
}
where $u_0$ is obtained from $w_0$ via the Biot-Savart law.
\vspace{0.5cm}\\
Notice that since $w_i$ belongs to $L^\infty_{loc}\left(\R^+,H^2(2)\right)$ and $\partial_t w_i$ belongs to $L^\infty_{loc}\left(\R^+,H^1(\R^2)\right)$, the inequalities (\ref{biotspq}) and (\ref{biotsnabla}) imply in particular
\cent{
$
\arr{l}{
u_i \in L^\infty_{loc}\left(\R^+,L^p(\R^2)^2\right), \quad \hbox{for all} \quad p>2,\\
\nabla u_i \in L^\infty_{loc}\left(\R^+,H^2(\R^2)^4\right),\\
\partial_t u_i \in L^\infty_{loc}\left(\R^+,L^p(\R^2)^2\right), \quad \hbox{for all} \quad p>2,\\
\partial_t \Delta u_i \in L^\infty_{loc}\left(\R^+,L^2(\R^2)^2\right).
}
$
}
Consequently, the system (\ref{g3diffui}) has a meaning in the sense of distributions.
\vspace{0.5cm}\\
We note $w=w_1-w_2$, $u=u_1-u_2$, $L=L_1-L_2$ and $A=A_1-A_2$. A short computation shows that $u$ satisfies the system 
\eq{\label{g3diffu}
\arr{l}{
\partial_t\left(u-\alpha_1 \Delta u\right)-\Delta u+\curl \left(u-\alpha_1\Delta u\right)\wedge u_1+\curl \left(u_2-\alpha_1\Delta u_2\right)\wedge u\\
\hspace{6cm}+\beta \dv \left(\left|A_2\right|^2 A_2\right)-\beta \dv \left(\left|A_1\right|^2 A_1\right)+\nabla q = 0,\\
\dv u=0,\\
u_{\left|t=0\right.} = 0.
}
}
Notice that, although $u_1$ and $u_2$ do not belong to $L^2(\R^2)$, the divergence free vector field $u$ does. Indeed, since $w_1$ and $w_2$ have the same initial data, for all $t\geq 0$, we have
\cent{
$
\displaystyle \int_{\R^2} w(t,x) dx =0.
$
}
By application of the lemma \ref{biots2}, this fact implies that $u$ belongs to $L^2(\R^2)$. Let $t_0>0$ be a fixed positive time. We notice that both $w_1$ and $w_2$ are bounded in $L^\infty\left(\left[0,t_0\right],H^2(2)\right)$. More precisely, one has
\cent{
$
\sup\limits_{t\in \left[0,t_0\right]}\left(\left\|w_1(t)\right\|_{H^2(2)}+\left\|w_2(t)\right\|_{H^2(2)}\right) \leq C.$
}
Applying the lemma \ref{biots1}, it implies in particular
\cent{
$
\sup\limits_{t\in \left[0,t_0\right]}\left(\left\| u_i(t)\right\|_{L^4}+\left\|\nabla u_i(t)\right\|_{L\infty}+\left\|\Delta u_i(t)\right\|_{L^4}\right) \leq C, \quad \hbox{for} \quad i=1,2.
$
}
In order to show that $u\equiv 0$, we now perform estimates on the $H^1-$norm of $u$. The uniqueness of the solutions of (\ref{g3diffui}) has been shown in \cite{busuiociftimie04} for solutions with initial data in $H^2(\R^2)$. In our case, the proof is slightly simpler, because the vector field $u$ belongs to $H^3(\R^2)^2$. We consider the $L^2-$inner product of (\ref{g3diffu}) with $u$. First of all, integrating by parts, we notice that
\aligne{
\beta \left(\dv\left(\left|A_2\right|^2 A_2-\left|A_1\right|^2 A_1\right),u\right)_{L^2} & = \frac{\beta}{2} \left(\left|A_1\right|^2 A_1-\left|A_2\right|^2 A_2,A\right)_{L^2} \\
& = \frac{\beta}{4} \int_{\R^2}\left(\left|A_1\right|^2+\left|A_2\right|^2\right)\left|A\right|^2 dx\\
&\hspace{1cm}+ \frac{\beta}{4} \int_{\R^2}\left(\left|A_1\right|^2-\left|A_2\right|^2\right)\left(A_1+A_2\right):A dx \\
& = \frac{\beta}{4} \int_{\R^2}\left(\left|A_1\right|^2+\left|A_2\right|^2\right)\left|A\right|^2 dx\\
&\hspace{2cm}+ \frac{\beta}{4} \int_{\R^2}\left(\left|A_1\right|^2-\left|A_2\right|^2\right)^2dx.
}
Thus, using integrations by parts and the divergence free property of $u$, we have
\eq{\label{diffH1}
\arr{l}{
\displaystyle \frac{1}{2}\partial_t \left(\left\|u\right\|^2+\alpha \left\|\nabla u\right\|^2\right)+\left\|\nabla u\right\|^2+\frac{\beta}{4} \int_{\R^2}\left(\left|A_1\right|^2+\left|A_2\right|^2\right)\left|A\right|dx \\
\hspace{7cm} \displaystyle + \frac{\beta}{4} \int_{\R^2}\left(\left|A_1\right|^2-\left|A_2\right|^2\right)^2dx = I_1+I_2,
}
}
where
\cent{
$
\arr{l}{I_1=\left(\curl \left(u_2-\alpha_1\Delta u_2\right)\wedge u,u\right)_{L^2},\\
I_2 = \left(\curl u\wedge u_1,u\right)_{L^2},\\
I_3 = -\alpha_1 \left(\curl \Delta u \wedge u_1,u\right)_{L^2}.
}
$}
A short computation shows that $I_1$ vanishes. Indeed, we set $\omega = u_2-\alpha_1\Delta u_2$ and we recall the notation $u = \left(u^1,u^2,0\right)$ and $\curl \omega = \left(0,0,\partial_1 \omega_2-\partial_2 \omega_1\right)$. We have
\aligne{
I_1 &= \left(\curl \omega \wedge u, u\right)_{L^2}\\
&= -\left(\left(\partial_1 \omega_2 - \partial_2 \omega_1\right) u^2,u^1\right)_{L^2}+\left(\left(\partial_1 \omega_2 - \partial_2 \omega_1\right) u^1,u^2\right)_{L^2}\\
&=0.
}
Due to the boundedness of $u_1$ in $L^4(\R^2)$, applying Hölder inequalities we obtain
\aligne{
I_2 & \leq \left\|u_1\right\|_{L^4} \left\|\nabla u\right\|\left\|u\right\|\\
&\leq C(\alpha_1) \left(\left\|u\right\|^2+\alpha_1\left\|\nabla u\right\|^2\right).
}
Using \cite[Lemma A.1]{paicuraugelrekalo12}, we check that
\cent{
$
\displaystyle I_3 \leq C\alpha_1 \int_{\R^2} \left|\Delta u_1\right|\left|\nabla u\right|\left| u\right| dx + C\alpha_1\int_{\R^2} \left|\nabla u_1\right|\left|\nabla u\right|^2 dx.
$
}
Using Hölder inequalities, the Gagliardo-Nirenberg inequality and the Young inequality $a b \leq \frac{1}{4}a^4+\frac{3}{4}b^{4/3}$, we obtain
\aligne{
I_3 & \leq C\alpha_1 \left\|u\right\|_{L^4} \left\|\Delta u_2\right\|_{L^4}\left\|\nabla u\right\|+C\alpha_1\left\|\nabla u_1\right\|_{L^\infty}\left\|\nabla u\right\|^2\\
& \leq C \alpha_1\left\|\nabla u\right\|^{3/2}\left\| u\right\|^{1/2}+C\alpha_1\left\|\nabla u\right\|^2\\
& \leq C (\alpha_1) \left( \left\| u\right\|^2+\alpha_1\left\|\nabla u\right\|^2\right).
}
Going back to (\ref{diffH1}), we get
\eq{\label{diffH1.2}
\arr{l}{
\displaystyle \frac{1}{2}\partial_t \left(\left\|u\right\|^2+\alpha \left\|\nabla u\right\|^2\right)\leq C(\alpha)\left(\left\|u\right\|^{2}+\alpha \left\|\nabla u\right\|^{2}\right).
}
}
Integrating in time this inequality between $0$ and $t \in \left[0,t_0\right]$ and applying the Gronwall lemma, we finally obtain
\cent{
$
\left\|u(t)\right\|^2+\alpha \left\|\nabla u(t)\right\|^2 = 0, \quad \hbox{for all} \quad t \in \left[0,t_0\right].
$
}
Since $t_0$ is arbitrary, we conclude that $u\equiv 0$ on $\R^+$. Consequently $u$ is unique and so is $w$. Thus, the system (\ref{g3w}) has a unique global solution in the space $C^0\left(\R^+,H^2(2)\right)$.
\vspace{0.5cm}\\
\textbf{Acknowledgement:} I would like to express my gratitude to Genevieve RAUGEL and Marius PAICU, whose advice were very helpful for this work.

\end{document}